\documentclass[11pt, leqno]{amsart}
\usepackage{amsthm, amsmath}
\usepackage{amssymb} 
\usepackage{amscd}
\usepackage[curve,matrix,arrow,frame,tips]{xy}
\usepackage{verbatim}
\usepackage{graphicx}

\newtheorem{thm}{Theorem}[section]
\newtheorem{prop}[thm]{Proposition}
\newtheorem{lemma}[thm]{Lemma}

\newtheorem{Remarknumb}[thm]{Remark}
\newtheorem{Remarks}[thm]{Remarks}
\newtheorem{Remark}[thm]{Remark}

\newtheorem{conjecture}[thm]{Conjecture}
\newtheorem{cor}[thm]{Corollary}

\newcounter{ex}[section]

\newcommand{\cal}{\mathcal}

\newcommand{\ds}{\displaystyle}

\newcommand{\A}{{\mathcal A}}
\newcommand{\G}{{\bf G}}

\newcommand{\C}{{\bf C}}
\newcommand{\R}{{\bf R}}
\newcommand{\Q}{{\bf Q}}

\newcommand{\ep}{\epsilon}

\newcommand{\Hom}{{\rm Hom}}

\newcommand{\Gm}{{{\bf G}_m}}
\newcommand{\Z}{{\bf Z}}
\newcommand{\F}{{\mathcal F}}
\newcommand{\ti}{\tilde}
\newcommand{\Spec}{{\rm Spec }\, }
 \renewcommand{\O}{{\mathcal O}}

\newcommand{\und}{\underline}

\newcommand{\Lfr}{{\mathfrak L}}

\newcommand{\M}{{\mathcal M}}
\renewcommand{\L}{{\mathcal L}}

\newcommand{\Res}{{\rm Res}}

\newcommand{\Tr}{{\rm Tr}}

\newcommand{\LL}{{\mathcal L}}
\newcommand{\fp}{{{\bf F}_p}}

\newcommand{\DD}{{\mathcal D}}

\newcommand{\fG}{{\mathfrak g}}
\newcommand{\fF}{{\mathfrak f}}
\newcommand{\fP}{{\mathfrak p}}

\def\thfill{\null\nobreak\hfill}

\def\endproof{\thfill\vbox{\hrule
  \hbox{\vrule\hbox to 5pt{\vbox to 5pt{\vfil}\hfil}\vrule}\hrule}}

\DeclareMathOperator{\abs}{abs}
\DeclareMathOperator{\der}{der}
\DeclareMathOperator{\ord}{ord}
\DeclareMathOperator{\val}{val}

\DeclareMathOperator{\Ker}{Ker}

\renewcommand{\to}{\longrightarrow}
\newcommand{\mto}{\longmapsto}

\newcommand{\la}{\langle}
\newcommand{\ra}{\rangle}

\newcommand{\mR}{{\mathbb R}}
\newcommand{\mZ}{{\mathbb Z}}

\newenvironment{eq}{\addtocounter{thm}{0}\begin{equation}}{\end{equation}}

\begin{document}

\title[Local Models III]{Local Models in the ramified case. \hfill
\center III.  Unitary groups}
\author[G. Pappas]{G. Pappas*}
\thanks{*Partially supported by  NSF Grant DMS05-01409.}
\address{Dept. of
Mathematics\\
Michigan State
University\\
E. Lansing\\
MI 48824-1027\\
USA}
\email{pappas@math.msu.edu}
\author[M. Rapoport]{M. Rapoport}
\address{Mathematisches Institut der Universit\"at Bonn,  
Beringstrasse 1\\ 53115 Bonn\\ Germany.}
\email{rapoport@math.uni-bonn.de}

\date{\today}
 
\maketitle


\medskip

\centerline{\sc Contents}
\medskip

\noindent  \ \ \  \ \ Introduction\\
\S 1.  Unitary Shimura varieties and moduli problems\\
\S 2.  Affine Weyl groups and affine flag varieties; the $\mu$-admissible set\\
\S 3.  Affine Flag varieties\\
\S 4. The structure of local models\\
\S 5. Special parahorics\\
\S 6. The local models of Picard surfaces\\
\S 7. Functor description; the spin condition\\
\S 8. Remarks on the case of (even) orthogonal groups
\bigskip
\medskip

\section*{Introduction}

One of the basic problems in the arithmetic theory of Shimura varieties is the construction of natural models over the ring of integers $\O_E$ of the reflex field $E$. One would like to construct models which are flat and have mild singularities. Such models should be useful in various arithmetic applications.
In particular, they should allow one to use the Lefschetz fixed point formula for the complex of nearby cycles to calculate the semi-simple zeta function at a non-archimedean prime of $E$. 

When the Shimura varieties is the moduli space over $\Spec (E)$ of abelian varieties with additional polarization, endomorphisms and level structure (a Shimura variety of PEL type), versions of the moduli problem sometimes make sense over $\Spec (\O_E)$. In this case the corresponding moduli schemes define models over $\O_E$, and one may ask whether they satisfy the requirements spelled out above. 

Let us fix a prime  number $p$ and consider a level structure of parahoric type in $p$. In this case one can expect that the problem of defining a natural integral model has a positive solution. When the parahoric is hyperspecial, and the group $G$ defining the Shimura variety has as simple factors only groups of type $A$ or $C$, Kottwitz [Ko] has shown that the model defined by the natural extension of the moduli problem has good reduction (at least when $p\neq 2$). When the parahoric is defined in an elementary way as the stabilizer of a selfdual periodic lattice chain in a $p$-adic vector space, integral models of the above type were defined in [RZ]. When the group $G$ only involves types $A$ and $C$ and splits over an unramified extension of $\Q_p$, all parahorics can be described in this way and G\"ortz [G\"o1, G\"o2] has shown in this case that the models defined in [RZ] are indeed flat with reduced special fiber and with irreducible components which are normal and with rational singularities. This follows earlier work by Deligne and Pappas [DP], by Chai and Norman [CN], and by de Jong [J]. In a few cases, these singularities can be resolved to construct semi-stable models or at least models with toroidal singularities. Here we mention work by Genestier [Ge1], by Faltings [Fa1-2], and by G\"ortz [G\"o3].  When the group $G$ is of type $A$ or $C$ and is split over $\Q_p$, Haines and Ngo
[HN] have calculated the semi-simple trace of Frobenius on the sheaves of nearby cycles at any point in the reduction of the natural model  at any prime ideal of $\O_E$ of residue characteristic $p$. Indeed, they give a group-theoretical expression for this semi-simple trace that had been conjectured earlier by Kottwitz. 

In [P1] the case of a unitary group  that corresponds at $p$ to a {\sl ramified} quadratic extension of $\Q_p$ was considered. The parahoric subgroup considered in [P1] is the stabilizer of a selfdual lattice. In this case, it was shown in loc. cit. that the models proposed in [RZ] are not flat in general. Furthermore, a closed subscheme of this model was defined and it was conjectured that this closed subscheme is flat and has other good properties. This conjecture is still open, although it is proved in [P1] for signature type $(r, 1)$. 
In the intervening years it became clear that more generally, when the group $G$ splits over a ramified extension of $\Q_p$, the models defined in [RZ], which were subsequently renamed {\it naive models},  are not the correct ones. Also, as was pointed out by Genestier to one of the authors several years ago, the naive models have pathological properties in the case of even orthogonal groups, even those that split over $\Q_p$. 

There are two ways to overcome the shortcomings of the naive models. The first one is to force flatness  by taking the flat closure of the generic fiber in the naive model, and to investigate the properties of the models obtained in this way. The second one is to strengthen the naive formulation of the moduli problem in order to obtain closed subschemes of the naive models and to show that these models have good properties, like flatness. In the case that the parahoric subgroup is the stabilizer of a self-dual periodic lattice chain, there is a well-known procedure ([DP], [RZ]) to reduce these questions to problems of the corresponding {\it naive local models}. The advantage of this approach is that we are then dealing with varieties that can be defined in terms of linear algebra. Furthermore, the second approach sometimes leads in this way to very explicit problems on the structure of varieties given by matrix identities. 

In [PR1] and [PR2] we considered groups $G$ which after localization over $\Q_p$ are of the form ${\rm Res}_{F/\Q_p}G'$, where
$G'$ is the general linear group or the group of symplectic similitudes, and where $F$ is a ramified extension of $\Q_p$. We gave three ways of defining good models in this case, which are all in the spirit of the first approach. The first one is in terms of a {\it splitting model} which in turn is defined in terms of the naive local model for the group $G'$. The second one is by flat closure. The third one is by flat closure for the maximal parahorics, and then by taking intersections inside the naive local model for $G$. It was shown in [PR2] that all three methods lead to the same models. The  approach through a strengthening of the moduli problem  turned out to be very complicated in the case of $G'=GL_n$. In the case $G'=GSp_{2n}$, we conjectured that the naive model is flat. This conjecture is still open, but G\"ortz [G\"o4]  has at least proved that the naive model is topologically flat in this case. Furthermore, we gave in [PR2] a partial calculation of the semi-simple trace of Frobenius on the sheaves of nearby cycles for these cases, by reduction to the theorem of Haines and Ngo (here `partial' means that we only compute these traces after a ramified base change). Note that in the cases considered in [PR2], all parahorics are  stabilizers of selfdual periodic lattice chains. 

In the present paper we deal with the other typical case of a group $G$ which splits over a ramified extension of $\Q_p$, namely the group of unitary similitudes corresponding to a quadratic extension of $\Q$ which is ramified at $p$. At the end we also comment on the case of an even orthogonal group since this case is closely intertwined with the case of the ramified unitary group when one approaches the construction problem in the second way. 

To explain our results, we need to introduce some notation. We consider the group $G$ of unitary similitudes for a hermitian vector space $(W, \phi)$ of dimension $n\geq 3$ over an imaginary quadratic  field $K\subset \C$, and fix a conjugacy class of  homomorphisms $h: {\rm Res}_{\C/\R} \G_m\to G_{\R}$ corresponding to a Shimura datum $(G, X_h)$ of signature $(r, s)$ with  $s\leq r$. We assume that $K/\Q$ is ramified over  $p$ and that $p\neq 2$. Let $F=K\otimes\Q_p$ and $V=W\otimes_{\Q}\Q_p$. We assume that the hermitian form on $V$ is split, i.e that there is a basis $e_1, \ldots,e_n$ such that 
$$
\phi(e_i, e_{n-j+1})=\delta_{i j}\ ,\ \forall i,j=1,\ldots,n \ .
$$
We fix a square root $\pi$ of $p$. For $i=0,\ldots,n-1$, set 
$$
\Lambda_i={\rm span}_{\O_F}\{\pi^{-1}e_1,\ldots, \pi^{-1}e_i, e_{i+1},\ldots, e_n\} \ .
$$
We complete this into a selfdual periodic lattice chain by setting $\Lambda_{i+kn}=\pi^{-k}\Lambda_i$. Let $n=2m+1$ when $n$ is odd and $n=2m$ when $n$ is even. 

It turns out that when $n=2m+1$ is odd, the stabilizer of a partial selfdual periodic lattice chain is always a parahoric. The conjugacy classes of parahoric subgroups of $G(\Q_p)$ correspond in this way to non-empty subsets $I$ of $\{0,\ldots,m\}$ (stabilizer of the lattices $\Lambda_j$, where $j=\pm i+kn$ for some $i\in I$ and some $k\in \bf Z$). When $n=2m$ is even, the situation is a little more complicated, since the stabilizer groups sometimes contain a parahoric subgroup with index $2$. Also, the conjugacy classes of parahoric subgroups correspond in this case to subsets $I$ of $\{0,\ldots,m\}$ with the property that if $m-1\in I$ then also $m\in I$. Let us denote by $M^{\rm naive}_I$ the {\it naive local
model} in the sense of [RZ], cf. further below. This is a projective scheme over $\Spec(\O_E)$. We also introduce the {\it local model}, the flat closure $M^{\rm loc}_I$ of the generic fiber in $M^{\rm naive}_I$.  In all cases we construct (naive resp. flat) models of the Shimura variety $Sh_C(G, X_h)$, where $C$ is an open compact subroup of the finite adele group with parahoric $p$-component,  which is \'etale locally around each point isomorphic to the naive local model resp. the
local model for suitable $I$. To understand the structure of these models of the Shimura variety it therefore suffices to study their local versions.  

The main tool  towards this goal is, as already in G\"ortz' paper [G\"o1], the embedding of the geometric special fiber  of the naive local model in a partial affine flag variety over the algebraic closure $k$ of the residue field. 
In our case, we have to employ the  affine flag varieties 
for non-split groups whose
theory is developed in [PR3].  More precisely, let $H$ be the quasi-split unitary group   over $k((t))$ which corresponds to a ramified quadratic extension of $k((t))$ and let $\F$ be the affine flag variety of $H$.   We then obtain a closed embedding of the special fiber of $M^{\rm naive}$ in $\F$. Here $ M^{\rm naive}$ denotes the naive local model for $I=\{0,\ldots,m\}$. This embedding is equivariant for the action of the Iwahori subgroup,
and hence its image is a union of Schubert varieties. These Schubert varieties are enumerated by certain elements of the affine Weyl group $W_a$ of $H$. To describe this subset of $W_a$, recall that the Shimura datum $(G, X_h)$ defines a minuscule  coweight  $\mu=\mu_{r, s}$ of $H$, and that to a coweight $\mu$ there is associated a finite subset ${\rm Adm}(\mu)$ of $W_a$, the {\it $\mu$-admissible set} [R]. One of the main results of the present paper is the following theorem.
\begin{thm}\label{thmA}
The union of Schubert varieties over the $\mu$-admissible set $$
\A(\mu)=\bigcup\nolimits_{w\in {\rm Adm}(\mu)}S_w
$$
 is contained in 
$\bar M^{\rm loc}$, the special fiber of the local model. If the coherence conjecture of {\rm [PR3]} is true, then this containment is in fact an equality and $\bar M^{\rm loc}$ is reduced and all its irreducible components are normal and with rational singularities.
\end{thm}
The coherence conjecture of [PR3] is an explicit formula for the 
dimension of the space of global sections on $\A(\mu)$ of the natural ample line bundle $\L$ on $\F$. The conjecture is proved in the cases of $GL_n$ and of $GSp_{2n}$, but is open in general for ramified unitary groups. 

Something analogous is proved  here also for local models $ M^{\rm loc}_I$ for proper subsets $I$ of $\{0,\ldots,m\}$. In various special cases we can however dispense with the coherence conjecture of [PR3]. We prove the following result.
\begin{thm}\label{thmB}
Let $I=\{0\}$ if $n$ is odd, and $I=\{m\}$ if $n=2m$ is even. The special fiber of the local model $M_I^{\rm loc}$ is irreducible and reduced and is normal, Frobenius split and with only rational singularities. 
\end{thm}

We note that the parahoric subgroups corresponding to the subsets $I$ in this theorem are {\it special} in the sense of Bruhat-Tits theory,
but that when $n$ is odd, there are special parahorics which are not conjugate to the parahoric subgroup for $I=\{0\}$. 

For the case of Picard surfaces, we have a complete result for all parahorics.
\bigskip
\vfill\eject

\begin{thm}\label{picard} Let $n=3$.  

a) $M^{\rm naive}_{\{0\}}$ is normal and Cohen-Macaulay. Furthermore,  $M^{\rm naive}_{\{0\}}$  is flat over $ \Spec(\O_E)$ and is smooth outside a single point  of the special fiber. Blowing up this  point yields a semi-stable model with special fiber consisting of two smooth surfaces meeting transversely along a smooth curve.

b) $M^{\rm loc}_{\{1\}}$ is smooth over $\Spec(\O_E)$, with special fiber isomorphic to 
${\bf P}^2$.

c)  $M^{\rm loc}_{\{0, 1\}}$ is normal and Cohen-Macaulay, with  reduced special fiber. Its special fiber  has two irreducible components which are normal and with only rational singularities.  These two irreducible 
 components meet along two smooth curves which intersect transversally in a single point.
\end{thm}

It is an interesting question whether local models can be defined by first forming the flat closure for maximal parahorics, i.e for subsets $I$ consisting of a single element, and then taking intersections of their inverse images in the naive local model. This leads to the combinatorial question whether the $\mu$-admissible set can be defined vertex-wise.  Figures 1 and 2 show that this is indeed the case for the cases of relative rank $2$.

 Let us briefly explain these figures. The ambient tesselation by alcoves corresponds to the affine Weyl group associated to a finite root system $\Sigma$ of type $B_2$ in the case of figure 1, resp. of type $C_2$ in the case of figure 2. The dots in these figures indicate the translation subgroup $Q(\Sigma^\vee)$ of $W_a$; the vertices of the simplices correspond to the elements in the coweight lattice $P(\Sigma^\vee)$. The base alcove is marked in bold face, and the extreme alcoves in ${\rm Adm}(\mu)$ are shaded in darker gray. They are translates of the fundamental alcove under the four
translation elements $\lambda_{r,s}$ in the Iwahori Weyl group associated to the coweight $\mu=\mu_{r,s}$. The other elements of ${\rm Adm}(\mu)$ are shaded in lighter gray.  
\bigskip

\begin{figure}[h]
	\begin{center}
		\includegraphics[width=12cm]{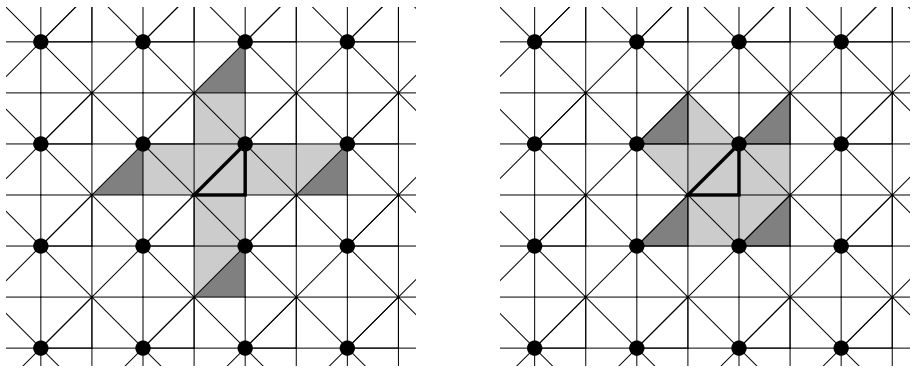}
		\caption{The admissible sets for $U(2,2)$ (left) and $U(3,1)$ (right)}
	\end{center}
	\label{fig:U4}
\end{figure}

\begin{figure}[h]
\includegraphics[width=12cm]{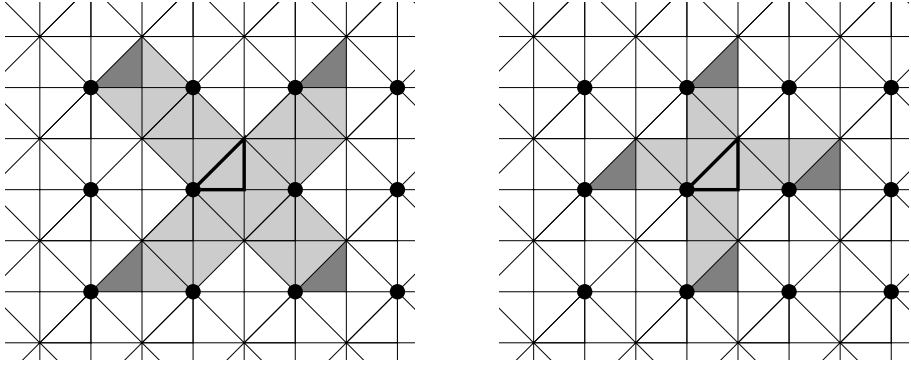}
\caption{The admissible sets for $U(3,2)$ (left) and $U(4,1)$ (right)}
\end{figure}

\vfill\eject

In this paper we also pursue the   approach of defining flat models through a strengthening of the moduli problem.   Recall from [RZ] the definition of the naive local model $M^{\rm naive}_I$. We change the notation slightly by now denoting by $E$  the localization at $p$ of the Shimura field, i.e  $E=\Q_p$ if $r=s$, and $E=F=K\otimes \Q_p$ if $r\neq s$.  The naive local model represents a moduli problem over $\Spec (\O_E)$.  A point of $ 
M^{\rm naive}_I$ with values in an $\mathcal O_E$-scheme $S$
is given by a $\mathcal O_{F}\otimes\mathcal O_S$-submodule
$$
\mathcal F_j\subset \Lambda_j\otimes\mathcal O_S
$$ 
for each $j\in\bf  Z$ of the form $j=\pm i+kn$ for some $k\in \bf Z$ and some $i\in I$. The following
conditions a)--d) are imposed.
\begin{itemize}
\item[a)] As an $\mathcal O_S$-module, $\mathcal F_j$ is
locally on $S$ a direct summand of rank $n$.
\item[b)] For each $j<j'$, there is a commutative diagram
$$
\begin{matrix}
\Lambda_j\otimes\mathcal O_S
&
\longrightarrow
&
\Lambda_{j'}\otimes\mathcal O_S
\\
\cup
&&
\cup
\\
\mathcal F_j & \longrightarrow & \mathcal F_{j'}
\end{matrix}
$$
where the top horizontal map is induced by the inclusion $\Lambda_j\subset \Lambda_{j'}$,
and for each $j$, the isomorphism $\pi: \Lambda_j\to\Lambda_{j-n}$ induces an isomorphism of  $\mathcal F_j$ with $\mathcal F_{j-n}$.

\item[c)] We have $\mathcal F_{-j}=\mathcal F_j^\bot$ where $\mathcal F_j^\bot$
is the orthogonal complement of $\mathcal F_j$ under the natural
perfect pairing 
$$
(\Lambda_{-j}\otimes\mathcal O_S)\times
(\Lambda_j\otimes\mathcal
O_S)\longrightarrow \mathcal O_S\ .
$$
\item[d)] For each $j$, the characteristic polynomial
equals
$$
\det((T\cdot {\rm id}- \pi)\mid \mathcal F_j)= (T-\pi)^s\cdot
(T+\pi)^r\in \mathcal O_E[T]\ \ .
$$
\end{itemize}
Here $\pi$ denotes a uniformizer of $F$ with $\pi^2=p$. 

In [P1] a closed subscheme of $M^{\rm naive}_I$ was defined by imposing an additional condition.

 \begin{itemize}
 \item[e)]  If $r\neq s$,  we have 
\begin{equation}
\wedge^{r+1}(\pi -\sqrt{p}\ |\ \F_j)=0\ ,
\end{equation}
\begin{equation}
\wedge^{s+1}(\pi+\sqrt{p} \ |\ \F_j)=0\ .
\end{equation}
\end{itemize}

Denote  by $M^{\wedge}_I$ the corresponding  closed subscheme of $M^{\rm naive}_I$.  It has the same generic fiber.
In [P1] it is conjectured that $M^{\wedge}_I$ is flat over $\Spec (\O_E)$ when $I=\{0\}$, i.e that condition e) above cuts out the local model in $M^{\rm naive}_{\{0\}}$ . The conjecture for $I=\{0\}$ is then reduced to a question about schemes 
defined by explicit matrix identities. 
It turns out that something similar can be done for  $I=\{m\}$ when $n=2m$. 
Let us state here these questions.

\smallskip\noindent
{\it Consider the scheme of matrices
$X$ in $M_{n\times n}$ over $\Spec(k)$ described by
$$
X^2=0, \quad X^t=HXH, \quad
 {\rm char}_X(T)=T^n\ ,\quad \wedge^{s+1}X=0,\quad  \wedge^{r+1}X=0\ ,
$$
 resp. 
$$
X^2=0, \quad X^t=-JXJ, \quad
 {\rm char}_X(T)=T^n , \ (\text{and } \wedge^{s+1}X=0,\quad  \wedge^{r+1}X=0 , \text{ when } r\neq s) ,
$$
if both $n$ and $s$ are even. 

Is this scheme  reduced (in which case  it is normal, with rational singularities)? }
\smallskip

Here $H=H_n$ denotes the antidiagonal unit matrix, and $J=J_{2m}$ the skew-symmetric matrix 
with  square blocks $0_m$ on the diagonal and $H_m$, resp. $-H_m$, above the diagonal, resp. below the diagonal. Varieties given by similar matrix equations were also considered by Faltings [Fa1]. However, it does not seem
that his results can be used   to answer the above question. Indeed,
the schemes he considers are the special fibers of the (normalizations)
of the flat closures of mixed characteristic schemes closely related to $M_I^{\rm naive}$. We can only see, after the fact, that if the answer to our question is positive,
then our schemes coincide with the ones considered in [Fa1].
(Regardless, the results of [Fa1] establish a connection between 
the singularities of the normalization of the local models and the complete symmetric varieties of de Concini-Procesi; see also [P2].)

We give examples which show that for more general $I$ the condition e) is not enough to define the local model. To treat a general index set $I$  we propose in the present paper an additional  condition. 

\begin{itemize}
\item[ \ f)]  \ ({\sl Spin condition}) \ The line  $\wedge^n \F_j\subset   \wedge^n(\Lambda_j\otimes\O_S)$ is contained in 
the subspace   $(\wedge^n\Lambda_j\otimes\O_E)_{\pm}\otimes_{\O_E}\O_S$ with $\pm=(-1)^s$.
\end{itemize}

We refer to the body of the text for the definition of $(\wedge^n\Lambda_j\otimes\O_E)_{\pm}\otimes_{\O_E}\O_S$. We denote by $M_I$ the closed subscheme of $M^{\rm naive}_I$ defined by the conditions a)--f). Then $M_I$ has again the same generic fiber and we conjecture that $M_I$ is flat over $\O_E$. Unfortunately, this seems very difficult to prove, although the computer evidence seems quite convincing. It is remarkable that, in contrast to the unramified cases of type $A$ or $C$, one seems to need in these ramified cases and also for type $D$  higher tensors to describe the local model. 

We now give an overview of the paper. In \S 1 we define the unitary Shimura varieties and reduce our structure problem to the case of local models. In \S 2 we recall some facts about affine Weyl groups and about the $\mu$-admissible set and make these explicit in the case of unitary groups. In \S 3 we recall results from [PR3] on affine flag varieties, and state the coherence conjecture. We also prove in this section that the points in Schubert varieties corresponding to elements in ${\rm Adm}(\mu)$ can be lifted to the generic fiber of $M^{\rm naive}$, i.e the first statement of Theorem \ref{thmA}. In \S 4 we give our general results. As pointed out above, they depend on some conjectures of a combinatorial nature. In \S 5 we consider cases we can treat unconditionally. In particular, we relate the special fibers in the cases addressed in Theorem \ref{thmB} to nilpotent orbits for the symmetric pairs $({\mathfrak gl}(n), {\mathfrak o}(n))$ resp. $({\mathfrak gl}(n), {\mathfrak sp}(n))$, and recover in this way some of Ohta's results about them [Oh], even in positive characteristics. We also treat in \S 6 the case of Picard surfaces. In \S 7 we explain the spin condition and the evidence we have that it is sufficient to cut out the flat closure. In \S 8 we give remarks on local models for even orthogonal groups and explain on examples how the spin condition also seems to eliminate the pathological points in the naive models for them found by Genestier. 

It should be pointed out that local models for ramified unitary groups are still quite mysterious. Many questions on them remain open. This also explains why in  the present paper we have made abstract concepts as explicit  as possible for the case at hand and why we have given many concrete examples. We hope that this will be useful for future work. One major open problem is that of calculating the semi-simple trace of the Frobenius on the complex of nearby cycles. Besides the special case treated by Kr\"amer [Kr] nothing is known. 

The present paper should be viewed as a follow-up to [PR3]. The results in [PR3] on affine flag varieties for non-split groups over $k((t))$  are crucial to several of our results in this paper.

We thank A. Genestier, U. G\"ortz   and T. Wedhorn for helpful discussions. We are especially indebted to  C. Kaiser and J.-L.  Waldspurger for their explanations on Bruhat-Tits theory. We are grateful to the Institute for Advanced Study and the Universit\'e de Paris-Sud for their support during our collaboration.  
\medskip


\section{Unitary Shimura varieties and moduli problems}
\setcounter{equation}{0}

\subsection{Unitary Shimura varieties}\label{uniShimura} Let $K$ be an imaginary quadratic field 
with an embedding $\ep: K\hookrightarrow \C$. Denote by $\O$ the ring of integers of $K$. Denote by $a\mapsto \bar a$ the non trivial automorphism of $K$.
Let $W=K^n$ be a $n$-dimensional
$K$-vector space, and suppose that $\phi: W\times W\to K$ is a non-degenerate hermitian form.
We assume $n\geq 3$.
Now set $W_\C=W\otimes_{K, \ep}\C$. Choosing a suitable isomorphism $W_\C\simeq \C^n$, we
may write $\phi$ on $W_\C$ in a normal form $\phi(w_1, w_2)=\ {^t\bar w_1H w_2}$
where
$$
H={\rm diag}(-1,\ldots, -1, 1,\ldots ,1)
$$
We denote by $s$ (resp. $r$) the number of places,
where $-1$, (resp. $1$) appears in $H$. We will say that 
$\phi$ has signature $(r,s)$. By replacing $\phi$ by $-\phi$ if needed, we can make sure that $s\leq r$.
We will assume that $s\leq r$ throughout the paper.  
Let $J: W_\C\to W_\C$ be the endomorphism given by
the matrix $-\sqrt{-1}H$. 
 We have $J^2=-id$ and so the endomorphism $J$ gives  
 an ${\bf R}$-algebra homomorphism $h_0: \C\to {\rm End}_{\bf R}(W\otimes_\Q\R)$
with $h_0(\sqrt{-1})=J$ and hence a complex structure on $W\otimes_\Q\R=W_\C$. For this complex structure we have
$$
{\rm Tr}_{\C}(a ; W\otimes_\Q\R)=s\cdot \ep(a)+r\cdot \bar \ep(a), \ a\in K\, .
$$
Denote by $E$ the subfield of $\C$ which is generated by the traces above (the `reflex field'). This is equal to $\Q$ if $r=s$ and to $K$ if $r\neq s$.
The representation of $K$ on $W\otimes_\Q\R$ with the above trace is
defined over $E$, i.e there is an $n$-dimensional
$E$-vector space $W_0$ on which $K$ acts such that
\begin{eq} \label{traeq}
{\rm Tr}_{E}(a; W_0)=s\cdot a+r\cdot \bar a
\end{eq}
and such that $W_0\otimes_{E}\C$ together with the
above $K$-action is isomorphic to $W\otimes_\Q\R$
with the $K$-action induced by $\epsilon: K\hookrightarrow \C$ and the above complex structure.

Now let let us fix a
non-zero element $\alpha\in K$
with $\bar \alpha=-\alpha$.
Set
\begin{equation}
\psi(x, y)={\rm Tr}_{K/\Q}(\alpha^{-1}\phi(x,y))
\end{equation}
which is a non-degenerate alternating form $W\otimes_\Q W\to \Q$.
This satisfies
\begin{equation} \label{alt}
\psi(av, w)=\psi(v, \bar a w), \quad \hbox{for all\ }
a\in K, \ v, w\in W.
\end{equation}
By replacing $\alpha$ by $-\alpha$, we can make sure that
the symmetric $\R$-bilinear form on $W_\C$ given by
$\psi(x, Jy)$ for $x$, $y\in W_\C$ is positive definite.

Let $G=GU(\phi)$ be the 
unitary
similitude group of the form $\phi$, 
$$
G=GU(\phi)=\{ g\in GL_K(W)\ | \ \phi(gv, gw)=c(g)\phi(v, w), \ c(g)\in \G_m \}\ .
$$

Set
$$
GU(r, s):=\{A\in GL_n(\C)\ | \ ^t\bar AHA=c(A)H, \ c(A)\in \R^{\times}\}
$$
By the above discussion, the embedding
$\ep: K\hookrightarrow \C$ induces an isomorphism
$G(\R)\simeq GU(r, s)$. The group $G$ is a reductive group over $\Q$ which is also given by
$$
G(\Q)=\{ g\in GL_K(W)\ | \ \psi(gv, gw)=c(g)\psi(v, w), \ c(g)\in \Q^{\times} \}\ .
$$
 
We define a homomorphism $h: {{\Res}}_{\C/\R}{{\bf G}}_{{m, \C}}\to G_\R$
by restricting $h_0$ to $\C^{\times}$. Then $h(a)$ for $a\in \R^{\times}$ acts on $W\otimes_\Q\R$ by
multiplication by $a $ and   $h(\sqrt{-1})$ acts as $J $. Consider $h_\C(z,1): \C^{\times}\to G(\C)\simeq {\rm GL}_n(\C)\times \C^{\times}$. Up to conjugation $h_\C(z,1)$
is given by
\begin{equation}
\mu_{r,s}(z)=({\rm diag}(\,   z^{(s)}, \, 1^{(r)}\, ), \, z);
\end{equation}
this is a cocharacter of $G$ defined over the number field $E$.
(Here and in the rest of the paper, we write $a^{(m)}$ to denote a list of $m$ copies of $a$).

Denote by $X_h=X_{r,s}$ the conjugation orbit of $h(i)$ under $G({\bf R})$. 
The pair $(G, h)$
gives rise to a Shimura variety $Sh(G, h)$ which is defined over the reflex field $E$.
In particular, if $C=\prod_{v} C_v\subset G({\bf A}^f_{\Q})$, with 
$C_v\subset G(\Q_v)$,  is an open compact subgroup of the finite adelic points of $G$, we can consider $Sh_C(G, h)$;
this is a quasi-projective variety over $E$ whose set of complex points is identified with
$$
Sh_C(G, h)(\C)=G(\Q)\backslash X_{r,s}\times G({\bf A}^f_\Q)/C\ .
$$
Suppose now that $p$ is an {\sl odd} rational prime which ramifies 
in $K$.  We denote by $\O_{E_w}$ the 
ring of integers of the completion of $E$ at the unique place above $(p)$.
We are interested in the construction  of  models of $Sh(G, h)_C$
over  $\O_{E_w}$ and in the reductions of these models
when the ``level subgroup at $p$",  $C_p\subset G(\Q_p)$, is a parahoric subgroup of $G(\Q_p)$.
(We will also assume that the level ``away from $p$", i.e the subgroup $C^p=\prod_{v\neq p}C_v$,
is sufficiently small, i.e it is contained in the principal congruence subgroup for some $N\geq 3$
relatively prime to the discriminant of $K$).
\medskip

\subsection{Parahoric subgroups of the unitary similitude group.}\label{parahoric}

It turns out that the  parahoric subgroups are the neutral components of the subgroups
of the unitary similitude group over the local field that stabilize certain sets of lattices.
In this paragraph, we will explain this statement in a slightly more general context.

\subsubsection{} Let $F_0$ be a complete discretely valued field with ring of integers 
$\O=\mathcal O_{F_0}$ and perfect residue field $k$ of characteristic $\neq
2$ and uniformizer $\pi_0$. Let $F/F_0$ be a ramified quadratic extension and let
$\pi\in F$ be a uniformizer with $\pi^2=\pi_0$, so that $\overline{\pi}=-\pi$. Set $\Gamma={\rm Gal}(F/F_0)$. Let
$V$ be a $F$-vector space of dimension $n\geq 3$ and let
$$
\phi:\  V\times V\longrightarrow F
$$
be a $F/F_0$-hermitian form.  We assume that
$\phi$ is split. This means that there
exists a basis $e_1,\ldots, e_n$ of $V$ such that
$$
\phi(e_i,e_{n-j+1})= \delta_{ij}\ \ ,\ \ \forall\  i,j=1,\ldots, n\  .
$$
Set
$$
GU(V, \phi)=\{g\in GL_F(V)\ |\ \phi(gx, gy)=c(g)\phi(x, y),\ \ c(g)\in F_0^{\times}\}\, .
$$
We have an exact sequence 
of algebraic groups over $F_0$
\begin{equation}
1\to SU(V, \phi)\to GU(V, \phi)\to D\to 1\, .
\end{equation}
Here $SU(V, \phi)$ is also the derived group of $GU(V,\phi)$ and  $D$ is the  torus $D=T/(T\cap SU(V, \phi) )$ with $T$ the standard (diagonal) maximal torus of $GU(V,\phi)$.

\subsubsection{}\label{latt} We have two associated $F_0$-bilinear forms,
$$
(x,y)=\frac{1}{2}\Tr_{F/F_0}(\phi(x,y))\ \ ,\ \ \langle x,y\rangle = \frac{1}{2}\Tr_{F/F_0} (\pi^{-1}\cdot\phi (x,y))\ \ .
$$
The form $(\ ,\  )$ is symmetric while $\langle \ , \ \rangle$ is alternating. They satisfy the identities,
\begin{equation}
(x,\pi y)=-(\pi x, y), \quad \langle x,\pi y\rangle=-\langle \pi x,  y\rangle \ .
\end{equation}
For any $\mathcal O_{F}$-lattice $\Lambda$ in $V$ we set,
$$
\hat\Lambda =\{ v\in V\mid \phi (v,\Lambda)\subset
\mathcal O_{F}\} =\{ v\in V\mid \langle v,\Lambda\rangle 
\subset \O\} \ .
$$
Similarly, 
we set
$$
\hat\Lambda^s =\{ v\in V\mid (v,\Lambda)\subset \O\}\  ,
$$
so that $\hat\Lambda^s =\pi^{-1}\cdot\hat\Lambda$. For
$i=0,\ldots, n-1$, set
$$
\Lambda_i={\rm span}_{\mathcal O_{F}} \{ \pi^{- 1}
e_1,\ldots, \pi^{- 1} e_i, e_{i+1},\ldots, e_n\}\ .
$$ 
The lattice $\Lambda_0$ is self-dual for the alternating form $\langle \ , \ \rangle$.
\medskip

\subsubsection{}\label{oddeven}
 We now distinguish two cases:

a) $n=2m+1\geq 3$ is odd. Then we have 
\begin{equation*}
T={\rm diag}\left(a_1,\ldots ,a_{m}, a,  {a\bar a}{\bar a_{m}^{-1}},\ldots ,  {a\bar a}{\bar a_{1}^{-1}}\right),
\end{equation*}
\begin{equation*}
T\cap SU(V,\phi)=\{{\rm diag}\left(a_1,\ldots ,a_{m}, {\bar a}a^{-1},  \bar a_{m}^{-1},\ldots , \bar a_{1}^{-1}\right)\ |\ 
a=a_1\cdots a_m \} \ .
\end{equation*}
We can see that 
\begin{equation}
D=T/(T\cap SU(V,\phi))\xrightarrow{\sim} {\rm Res}_{F/F_0}(\Gm)\ .
\end{equation}
with the isomorphism given by sending an element of $T$ as above to 
$ a^{-1}\cdot  (a_1 \bar a_1^{-1}) \cdots  ( a_m \bar a_m^{-1})$.

Now let $I$ be a non-empty
subset of $\{0,\ldots, m\}$ and  consider the subgroup
\begin{equation*}
P_I=\{g\in GU(V, \phi)\ |\ g\cdot \Lambda_i= \Lambda_i,\ \forall i\in I\}\ .
\end{equation*}
of $GU(V, \phi)$ that preserves the lattice set $\Lambda_i$, $i\in I$. 
The following statement can be shown as in [PR3]. Notice that in this case 
$X_*(D)_\Gamma\simeq \Z$ and so the Kottwitz invariant of each element of $P_I$
is trivial.
\medskip

{\it The subgroup $P_I$ is a parahoric subgroup of $GU(V,\phi)$. Any parahoric subgroup of $GU(V,\phi)$ 
is conjugate to a subgroup $P_I$ for a unique subset $I$.  
The sets $I=\{0\}$ and $I=\{m\}$ correspond to the special maximal parahoric subgroups.} 
\medskip

In fact, $P_{\{0, \ldots , m\}}$ is an Iwahori subgroup 
and its choice allows us to identify $\{0,\ldots, m\}$
with the local Dynkin diagram $\Delta(GU(V, \phi))$; the index $i$ corresponds to the vertex associated to 
$P_{\{i\}}$. 
\medskip
 
b) $n=2m\geq 4$. Then we have 
\begin{equation*}
T=\{{\rm diag}\left(a_1,\ldots ,a_{m}, c{\bar a_{m}^{-1}},\ldots ,  c{\bar a_{1}^{-1}}\right)\ |\  c\in F_0^{\times}\}
\end{equation*}
\begin{equation*}
T\cap SU(V,\phi)=\{{\rm diag}\left(a_1,\ldots ,a_{m},  \bar a_{m}^{-1},\ldots , \bar a_{1}^{-1}\right)\ |\ 
a_1\cdots a_m\in F_0^{\times} \} \ .
\end{equation*}
We can see that 
\begin{equation}
D=T/(T\cap SU(V,\phi))\xrightarrow{\sim} \Gm\times {\rm ker}\big({\rm Res}_{F/F_0}(\Gm)\xrightarrow{\rm Norm} \Gm\big )
\end{equation}
with the isomorphism given by sending an element of $T$ as above to  
$$
(c, \, a_1\cdots a_m\cdot  \bar a_1^{-1}\cdots \bar a_m^{-1})\ .
$$
(Here of course we use Hilbert's theorem 90.) We can now see that $X_*(D)_\Gamma\simeq \Z\oplus \Z/2\Z$
and that the Kottwitz invariant of an element $t\in T(F_0)$ as above is given
by 
\begin{equation}
\kappa (t)=({\rm val}_{F_0}(c), {\rm val}_{F}(a_1\cdots a_m)\, {\rm mod}\, 2) \ .
\end{equation} 
This shows that in this case the Kottwitz invariant 
$\kappa (g)$ of $g\in GU(V, \phi)(F_0)$
can be obtained as follows: Consider $d(g):=c(g)^{-m} \det_F(g)\in F^{\times}$; this has norm $1$ and so we can write
$d(g)=x\cdot \bar x^{-1}$. Then $\kappa (g)=({\rm val}_{F_0}(c(g)), {\rm val}_{F}(x)\, {\rm mod}\, 2)$.
Notice that if $g$ stabilizes a maximal $\phi$-isotropic $F$-subspace $L\subset V$, then the second component of $\kappa (g)$ 
is the valuation modulo $2$ of the determinant $\det(g\, | \, L)$. Also note that if $g$ preserves 
a lattice, then the first component of $\kappa (g)$ is zero. 

Now consider non-empty subsets $I\subset \{0,\ldots ,m\}$.   
As above, consider the subgroup
\begin{equation*}
P_I=\{g\in GU(V, \phi)\ |\ g\Lambda_i=\Lambda_{i}, \ \forall i\in I\}\ .
\end{equation*}
of $GU(V, \phi)$ that preserves the  lattices $\Lambda_i$, $i\in I$. 
We also consider the kernel of the Kottwitz homomorphism, i.e.,
\begin{equation*}
P^0_I=\{g\in P_I\ |\ \kappa_{H}(g)=1\}\ .
\end{equation*} 
We can see that if $m-1\in I$, then $P^0_I=P^0_{I\cup \{m\}}$.

In fact, in this case the following statement holds:
\medskip

{\it The subgroup $P^0_I$ is a parahoric subgroup of $GU(V,\phi)$. 
Any parahoric subgroup of $GU(V,\phi)$ is conjugate to a subgroup $P^0_I$
for a unique subset $I$ with the property that if $m-1$ is in $I$, then $m$
is also in $I$. For such a subset $I$, we have $P^0_I=P_I$ exactly when $I$ contains 
$m$. The set $I=\{m\}$ corresponds to a special 
maximal parahoric subgroup.}
\medskip

This follows from the results on parahoric subgroups of $SU(V, \phi)$
in [PR3]. To explain this,  we also introduce  
\begin{equation*}
\Lambda_{m'}={\rm span}_{\O_{F}}\{ \pi^{-1}e_1,\ldots, \pi^{-1}e_{m-1}, e_{m},  \pi^{-1}e_{m+1}, e_{m+2},\ldots , e_n\}\ .
\end{equation*}
Then both $\Lambda_m$ and $\Lambda_{m'}$ are self-dual for the symmetric form
$(\ ,\  )$. Now consider non-empty subsets $J$ of   $\{0,\ldots, m-2, m, m'\}$. Consider the subgroups
\begin{equation*}
P_J=\{g\in GU(V, \phi)\ |\ g\Lambda_j=\Lambda_{j}, \ \forall j\in J\}\ ,\quad  P^0_J=\{g\in P_J\ |\ \kappa_{H}(g)=1\}\ .
\end{equation*}
Notice that if both $m$ and $m'$ are in $J$ then $\Lambda_j$, $j\in J$, is not a lattice chain.
 
As in [PR3] we can see that $P_J$ is parahoric exactly when $J$ contains at least one
of the two elements $\{m, m'\}$. Then $P_J^0=P_J$. When  $J$  contains neither
$m$ nor $m'$, the kernel $P^0_J$ is a parahoric subgroup and $P_J/P^0_J\simeq\Z/2\Z$. 
Now recall that [PR3] gives a description of the parahoric subgroups of $SU(V,\phi)$
via $P'_J=P_J\cap SU(V,\phi)$. In fact, $P_{\{0, \ldots , m-2,m,m'\}}$ is an Iwahori subgroup 
and its choice allows us to identify  $\{0,\ldots, m-2, m, m'\}$
with the set of vertices of the local Dynkin diagram 
$\Delta:=\Delta(GU(V,\phi))=\Delta(SU(V,\phi))$; the index $j$ corresponds to the vertex 
associated to the subgroup $P_{\{j\}}$. 

For $J\subset \{0,\ldots ,m-2, m, m'\}$ let $J^*$ be the subset obtained 
by replacing $m$ by $m'$ and vice versa. Observe that if $\tau$ is the unitary automorphism
defined by $e_m\mapsto e_{m+1}$, $e_{m+1}\to e_m$, $e_i\mapsto e_i$, for $0\leq i\leq m-2$,
then $\tau\cdot \Lambda_m=\Lambda_{m'}$, $\tau\cdot \Lambda_{m'}=\Lambda_{m}$, $\tau\cdot \Lambda_i=\Lambda_{i}$,
for $0\leq i\leq m-2$. This shows that $\tau P_J\tau^{-1}=P_{J^*}$.
By [T] 2.5 we can see that the action of $GU(V,\phi)$ on the local Dynkin diagram 
$\Delta =\{0,\ldots, m-2, m, m'\}$ factors through the Kottwitz homomorphism 
$\kappa: GU(V,\phi)\to \Z/2\Z$; now $\kappa(\tau)=-1$ and as we have seen $\tau$ 
fixes $0\leq i\leq m-2$ and exchanges $m$ and $m'$. It follows that the conjugacy classes of
parahoric subgroups of $GU(V,\phi)$ are parametrized by the orbits of $J\mapsto J^*$ on the set 
of non-empty subsets $J\subset \{0,\ldots , m-2, m, m'\}$. These are in turn parametrized 
by non-empty subsets $I\subset \{0,\ldots, m-1, m\}$ with the property that if $m-1$ is in $I$, then $m$
is also in $I$: The set $I=J^\sharp=(J^*)^\sharp$ that corresponds to $\{J, J^*\}$ is obtained by the following recipe:
If both $m$ and $m'$ belong to $J$, put  $J^\sharp$ to be the set which is obtained from $J$ by replacing $m'$ by $m-1$.  
If $m'\not\in J$,  set $J^\sharp=J$. Finally, if $m'\in J$ but $m\not\in J$, let $J^\sharp =J^*$.
Observe now that for all $J$, we have $P_J=P_{J^\sharp}$ if $m\in J$ and $P_J=\tau P_{J^*}\tau^{-1} =\tau P_{J^\sharp}\tau^{-1}$
if $m'\in J$. Therefore our statements on the parahoric subgroups of $GU(V, \phi)$ now follow from these observations and the results on parahoric subgroups of $SU(V, \phi)$ in \S4 of [PR3].
\medskip

\subsection{Reduction to level subgroups that are lattice chain stabilizers.}

We return to the set up of \S \ref{uniShimura}, so $G=GU(W,\phi)$ is 
a unitary similitude group over $\Q$, and $X=X_{r,s}$. 
We will assume $s>0$. (When $s=0$ the corresponding Shimura varieties
are zero-dimensional.) Let $p$ be an odd rational prime which ramifies in $K$. Assume that the form $\phi$ 
is split on $V=W\otimes_{\Q}\Q_p$ so that the set-up of \S \ref{parahoric} applies. 
Let $C=C_p\cdot C^p$ with $C_p=P^0_I\subset G(\Q_p)$ a parahoric subgroup as in \S \ref{parahoric}
and where $C^p\subset G({\bf A}^{f,p})$ contains the principal
congruence subgroup for some $N\geq 3$ relatively prime to the discriminant of $K$.
Set $C'_p=P_I\subset G(\Q_p)$ for the corresponding stabilizer of the set
of lattices; we have either $P^0_I=P_I$ or $P_I/P^0_I\simeq \Z/2\Z$.
Set $C=C'_p\cdot C^p$. Our first observation is that the Shimura varieties
${\rm Sh}_{C}(G, X)$ and  ${\rm Sh}_{C'}(G, X)$ have isomorphic geometric connected
components. This follows essentially from the fact that $C\cap G_{\rm der}({\bf A}^f)=C'\cap G_{\rm der}({\bf A}^f)$
(see below). Therefore, from the point of view of constructing reasonable integral models over $\O_{E_w}$,
we may restrict our attention to ${\rm Sh}_{C'}(G, X)$; since $C'_p$ corresponds to a lattice set stabilizer, this 
Shimura variety is given by a simpler moduli problem.

Denote by 
$$
\nu: G\xrightarrow{\ }D=G/G_{\rm der}
$$
the maximal
torus quotient of $G$. 
We can see ([Ko], \S 7) that
\begin{equation}\label{dgroup}
\ D\simeq {\rm Res}_{K/\Q}(\Gm),\quad D\simeq \Gm\times {\rm ker}\big({\rm Res}_{K/\Q}(\Gm)\rightarrow \Gm\big ),\  \hbox{if $n$ is odd, resp. even} .
\end{equation}
 Since $G_{\rm der}$ is simply connected, the set of connected components of the complex Shimura variety
${\rm Sh}_C(G, X) (\C)$ can be identified with the double coset
\begin{equation*}
D(\Q)\backslash Y\times D({\bf A}^f)/\nu(C)
\end{equation*}
where $Y=D(\R)/{\rm Im}(Z(\R)\to D(\R))$, provided that $r>0, s>0$, comp. [M], \S 5, p.311. 
(Notice that $D(\R)=\C^{\times}$ if $n$ is odd, resp. $D(\R)=\R^{\times}\times U_1$ if $n$ is even, with $U_1$ the complex unit circle. The image of the real points of the center $Z(\R)$ is $\C^{\times}$, resp. 
$\R_{>0}\times U_1$. Therefore, $Y=\{1\}$ if $n$ is odd, and $Y=\{\pm 1\}$ if $n$ is even.)
The action
of the Galois group ${\rm Gal}(\bar E/E)$ on the set of connected components factors through ${\rm Gal}( E^{\rm ab}/E)$
and is given as follows: Consider the composition $\nu\circ \mu_{r,s}: \Gm_E\to G_E\to D_E$ which is defined over $E$, and set
$$
\rho={\rm Norm}_{E/\Q} \circ \nu\circ \mu_{r,s}\ :\ \Gm_E\to D_\Q\ .
$$
If $\sigma\in {\rm Gal}( E^{\rm ab}/E)$ corresponds to the idele 
$x_{\sigma}\in {\bf A}_E^{\times}$ via Artin reciprocity, then 
\begin{equation}
\sigma\cdot [y, d]=\rho(x_\sigma)\cdot [y, d]\ ,
\end{equation}
with $\rho: {\bf A}_E^{\times}\to D({\bf A}_\Q)$   given as  above. 
(Here we normalize the reciprocity isomorphism by asking that the local uniformizer corresponds to
the inverse of the Frobenius).

The conjugacy class of the cocharacter $\mu_{r,s}$ of $G_E$ 
defines a conjugacy class of a local cocharacter $\mu_{r,s}: \Gm_{E_w}\to G_{E_w}$.
Suppose first that $r=s$. Then $E_w=\Q_p$ and we can assume that 
$\mu_{r,s}$ is given by
\begin{equation}
a\, \mapsto\,  {\rm diag}(a^{(s)}, 1^{(r)}) \ \in\ T(\Q_p)\subset G(\Q_p)
\end{equation}
in the notation of \S \ref{parahoric}. If $r\neq s$
then $E_w=K_v$ and 
and we can assume that 
$\mu_{r,s}$ is given by
\begin{equation}
a\, \mapsto\,  ({\rm diag}(  a^{(s)}, 1^{(r)}), a)) \ \in\  
 GL_n(K_v)\times K_v^{\times}\simeq G(K_v)\, .
\end{equation}
Note that the isomorphism $GL_n(K_v)\times K_v^{\times}\simeq G(K_v)$
takes the conjugation action on $G(K_v)$ to the involution 
$(A, c)\mapsto (\bar c (A^*)^{-1}, \bar c )$ where $(A^*)^{-1}$ is the inverse of the
hermitian adjoint.
Therefore, we have 
\begin{equation}\label{normmu}
{\rm Norm}_{E_w/\Q_p}(\mu_{r,s})(a)={\rm diag}((a\bar a)^{(s)}, \bar a^{(r-s)}, 1^{(s)}) \in T(\Q_p)\subset G(\Q_p)\, .
\end{equation}
Under the identification of $D(\Q_p)$ as in \S\ref{oddeven} this gives 
\begin{equation}\label{rhoformula1}
\rho_{E_w}(a)=\bar a^{-1} \left(\frac{\bar a}{  a}\right)^{m-s}, \ \ \hbox{if $n$ is odd},
\end{equation}
\begin{equation}\label{rhoformula2}
 \rho_{E_w}(a)=\left( a\bar a , \left(\frac{\bar a}{ a}\right)^{m-s}\right), 
\ \ \hbox{if $n$ is even and $r\neq s$}\ , 
\end{equation}
\begin{equation}\label{rhoformula3}
\rho_{E_w}(a)=(a  , 1), 
\ \ \hbox{if $n$ is even and $r=s$.}
\end{equation}

Consider now the morphism between the 
corresponding Shimura varieties,
\begin{equation*}
\pi_{C, C'}: {\rm Sh}_{C}(G, X)\to {\rm Sh}_{C'}(G, X)\ .
\end{equation*}
Recall that if $n$ is odd, we always have $C=C'$ and so we may restrict our attention to the case when
$n=2m$ is even.

\begin{prop}
\, Assume $C\neq C'$. Then the morphism $\pi_{C,C'}$ is an \'etale $\Z/2\Z$-cover which 
 splits after base changing by an extension 
$E'/E$ which is of degree $1$ or $2$ and is unramified over $p$. If $m-s$ is odd, then 
the cover is not trivial, $[E':E]=2$ and $w$ remains prime in $E'$.
If either $r=s$ or more generally $m-s$ is even, then either $E'=E$ or $[E':E]=2$ and $w$ splits in $E'$.
\end{prop}

\begin{Proof} 
It is clear  that $\pi_{C, C'}$ is an \'etale cover of degree $1$ or $2$. 
Let us consider the corresponding map between geometric connected components:
\begin{equation}\label{concomp}
D(\Q)\backslash Y\times D({\bf A}^f)/\nu(C)\xrightarrow{\ \ } D(\Q)\backslash Y\times D({\bf A}^f)/\nu(C')\ .
\end{equation}
Obviously, this is a surjective group homomorphism with kernel either trivial or of order $2$. 
We will show that the kernel is actually of order $2$; this will imply  that the cover $\pi_{C, C'}$ becomes trivial after base changing to the algebraic closure $\bar E$. Consider the element 
$\alpha$ in $Y\times D({\bf A}^f)$ which is $1$ at each place except at $p$ where it is 
equal to $\alpha_p=(1, -1)\in \Q_p^{\times}\times {\rm ker}({\rm Norm}(K^{\times}_v\to \Q^{\times}_p))$. 
From \S \ref{parahoric}, we can see that $\alpha_p$ is in $\nu(C'_p)$ but not in $\nu(C_p)$. 
We claim that $\alpha$ gives a non-trivial element of $D(\Q)\backslash Y\times D({\bf A}^f)/\nu(C)$
whose image under (\ref{concomp}) is trivial:
Indeed, suppose there is $d\in D(\Q)$ such that $d\cdot \alpha\in \nu(C)$. Since $\alpha$ is in $\nu(C')$ we
obtain that $d$ is in the intersection $D(\Q)\cap \nu(C')$; this is given by pairs 
of a unit of $\Q$ and a unit of $K$ which are both congruent to $1$ modulo $N$. We can see that the only such
pair is $(1,1)$; this implies $\alpha\in \nu(C)$ which contradicts our choice. 
It remains to show that the cover becomes trivial
after base changing by an extension $E'/E$ as in the statement. By the formula
for the Galois action on the set of connected components and the above we see that
the cover $\pi_{C,C'}\otimes_E E_w$ is described via local class field theory by the map
\begin{equation*}
\rho_{E_w}\, {\rm mod}\, \nu(C'_p) :\  E^{\times}_w\xrightarrow{\ } D(\Q_p)/\nu(C_p)\to \Z/2\Z\ .
\end{equation*}
Here the last map is given by $(a , b/\bar b)\mapsto {\rm val}(b)\, {\rm mod}\, 2$.
The result now follows using (\ref{rhoformula2}), (\ref{rhoformula3}). \endproof

\end{Proof}
\begin{cor}
Assume $C\neq C'$. Let $\M_{C'}$ be a model of $ {\rm Sh}_{C'}(G, X)$ over $\Spec(\O_{E_w})$. Then there exists a unique model $\M_C$ of  ${\rm Sh}_{C}(G, X)$ over $\Spec(\O_{E_w})$, 
such that the following diagram is commutative

$$
\begin{matrix}
 {\rm Sh}_{C}(G, X)\otimes_EE_w & \hookrightarrow &
\mathcal M_{ C}
\\
\big\downarrow && \big\downarrow
\\
 {\rm Sh}_{C'}(G, X)\otimes_EE_w  & \hookrightarrow &
\mathcal M_{ C'}
\end{matrix}\ \ 
$$
in which the vertical arrows are finite etale. 
\qed
\end{cor}
\medskip

\subsection{Unitary moduli problems}\label{unimoduli}

Here we follow [RZ] to define moduli schemes over $\O_{E_w}$
whose generic fiber agrees with ${\rm Sh}_{C'}(G,X)\otimes_EE_w$
when $C'=C'_p\cdot C^p$ with $C'_p$ one of the lattice set stabilizer 
subgroups of \S \ref{parahoric}. More precisely, we consider 
non-empty subsets $I\subset\{0,\ldots , m\}$ where $m=[n/2]$ with
the requirement that for $n=2m$ even, if $m-1$ is in $I$, then $m$ is in $I$ too.
We will use the notations of \S \ref{parahoric} with $F_0=\Q_p$, $F=K_v$, $V=W\otimes_KK_v$.
We can extend $\Lambda_i$, $i\in I$, given as in \ref{latt}, to a periodic self-dual lattice chain 
by first including the duals
$\hat\Lambda^s_i=\Lambda_{n-i}$ for $i\neq 0$, and then all the $\pi$-multiples of our lattices:
For $j\in \Z$ of the form $j=k\cdot n\pm i$ with $i\in I$
we put
\begin{equation*}
\Lambda_j=\pi^{-k}\cdot \Lambda_i\ .
\end{equation*}
Then $\{\Lambda_j\}_j$ form a periodic lattice chain $\Lambda_I$
(with $\pi\cdot \Lambda_j=\Lambda_{j-n}$) which satisfies
$\hat\Lambda_j=\Lambda_{-j}$. 

We set $C'_p=P_I$. We will use the construction
of Chapter 6 of [RZ] applied to the current situation.
In the notation of loc. cit. we take   
$B=K$,  $*=$ the conjugation of $K/\Q$, $V=K^n=\Q^{2n}$,
$( x , y )=\psi(x,y)$. Then $G$ is the unitary similitude group as before.
We take  the selfdual multichain $\LL$
of lattices to be $\Lambda_j$, $j=kn\pm i$ with $i\in I$.  
For simplicity, set $R=\O_{E_w}$.

We consider the moduli functor $\A_{C'}$
over $\Spec(R)$ given in [RZ] Definition 6.9:

A point of $\A_{C'}$
with values in the $\Spec(R)$-scheme $S$ is the isomorphism class
of the following set of data $(A, \bar\lambda, \bar \eta)$:

1. An $\LL$-set of abelian schemes $A=\{A_j\}$, $j=kn\pm i$ with $i\in I$, over $S$
(terminology of loc. cit.). By definition, this amounts to the data of abelian schemes $A_j$ over $S$
up to prime-to-$p$-isogeny for each $j$, isogenies $A_j\to A_{j'}$ of height $\log_p(\Lambda_{j'}/\Lambda_j)$
for each pair $j\leq j'$, and periodicity isomorphisms $\theta^a: A_j\simeq A_{j-n\cdot {\rm val}(a)}$
for each $a\in \O_K\otimes\Z_{(p)}$. These should satisfy the conditions of loc. cit.

2. A $\Q$-homogeneous principal polarization $\bar\lambda$
of the $\LL $-set $A$;

3. A $C^p$-level structure
$$
\bar\eta: H_1(A, {\bf A}^p_f)=\Big(\prod_{l\neq p}T_l(A_j)\Big)\otimes\Q\simeq W\otimes {{\bf A}^p_f}\ {\rm mod}\ C^p
$$
that respects the bilinear forms of both sides up to a constant
in $({\bf A}^p_f)^{\times}$ (see loc. cit. for details).
\smallskip

The set $A$ should satisfy the determinant
condition (i) of loc. cit.
\medskip

Recall that we  assume that $C^p$ is sufficiently small, i.e
that it is contained in the principal
congruence subgroup of level $N$ for some $N\geq 3$ relatively prime to
the discriminant of $K$.
Then $\A_{C'}$ is representable by a quasi-projective
scheme over $\Spec(R)$ which we will also denote by $\A_{C'}$.
Since the Hasse principle is satisfied for the unitary group, we 
can see as in loc. cit. that there is 
a natural isomorphism 
\begin{equation}
\A_{C'}\otimes_R E_w= {\rm Sh}_{C'}(G,X)\otimes_EE_w\ .
\end{equation}
 \medskip

\subsection{Local models for $GU$.}\label{naive}

We fix non-negative integers $r,s$ with $n=r+s$
and consider $F/F_0$ as in \S \ref{parahoric}. We set
$E=F$ if $r\neq s$ and $E=F_0$ if $r=s$ (this is the reflex
field of the local model we are about to define). 
As in \S \ref{unimoduli}, we will consider  non-empty subsets $I\subset\{0,\ldots , m\}$ where $m=[n/2]$ with
the requirement that for $n$ even, if $m-1$ is in $I$, then $m$ is in $I$ too.

\subsubsection{} First consider $I=\{0,\ldots, m\}$ which gives as above a complete lattice chain 
$\Lambda_j$, $j\in \Z$. We define as follows a functor $  M^{\rm naive}$ on the
category of $\mathcal O_E$-schemes. A point of $ 
M^{\rm naive}$ with values in an $\mathcal O_E$-scheme $S$
is given by a $\mathcal O_{F}\otimes_{\mathcal
O_{F_0}}\mathcal O_S$-submodule
$$
\mathcal F_j\subset \Lambda_j\otimes_{\mathcal
O_{F_0}}\mathcal O_S
$$ 
for each $j\in\bf  Z$. The following
conditions are imposed.
\begin{itemize}
\item[a)] as an $\mathcal O_S$-module, $\mathcal F_j$ is
locally on $S$ a direct summand of rank $n$.
\item[b)] for each $j<j'$, there is a commutative diagram
$$
\begin{matrix}
\Lambda_j\otimes_{\mathcal O_{F_0}}\mathcal O_S
&
\longrightarrow
&
\Lambda_{j'}\otimes_{\mathcal O_{F_0}}\mathcal O_S
\\
\cup
&&
\cup
\\
\mathcal F_j & \longrightarrow & \mathcal F_{j'}
\end{matrix}
$$
where the top horizontal map is induced by the inclusion $\Lambda_j\subset \Lambda_{j'}$, and for each $j$, the isomorphism $\pi: \Lambda_j\to\Lambda_{j-n}$ induces an isomorphism of  $\mathcal F_j$ with $\mathcal F_{j-n}$.

\item[c)] we have $\mathcal F_{-j}=\mathcal F_j^\bot$ where $\mathcal F_j^\bot$
is the orthogonal complement of $\mathcal F_j$ under the
perfect pairing
$$
(\Lambda_{-j}\otimes_{\mathcal O_{F_0}}\mathcal O_S)\times
(\Lambda_j\otimes_{\mathcal O_{F_0}}\mathcal
O_S)\longrightarrow \mathcal O_S
$$ induced by $\langle\ ,\
\rangle\otimes_{\mathcal O_{F_0}}\mathcal O_S$.
\end{itemize}

 Next note that $\mathcal F_j$ is an $\mathcal
O_{F}\otimes_{\mathcal O_{F_0}}\mathcal O_S$-module, hence
$\mathcal O_{F}$ and $\mathcal O_E$  act on it. We require
further that
\begin{itemize}
\item[d)] for each $j$, the characteristic polynomial
equals
$$
\det((T\cdot {\rm id}- \pi)\mid \mathcal F_j)= (T-\pi)^s\cdot
(T+\pi)^r\in \mathcal O_E[T]\ \ .
$$
\end{itemize}

This concludes the definition of the functor $M^{\rm naive}$, which is obviously representable by a
projective scheme over $\Spec (\mathcal O_E)$. We call
$M^{\rm naive}$ the {\it naive} local model
associated to the group $GU(V,\phi)$, the signature type
$(r,s)$ and for the complete lattice chain $\Lambda_j$,
$j\in\bf  Z$.  

\subsubsection{} We can generalize this definition to
incomplete selfdual periodic lattice chains,
i.e to all subsets $I$ as above. For each such $I$,
we obtain a functor $  M_I^{\rm
naive}$ by only giving the submodules $\mathcal F_j$ for
$j\in \Z$ of the form $j=k\cdot n\pm i$ with $i\in I$. Therefore $ 
M^{\rm naive}= 
M^{\rm naive}_{\{0,\dots ,m\}}$. Denote by  ${\mathcal P}_I$ the (smooth) group scheme 
of automorphisms (up to similitude) of the polarized chain $\L$
over $\O_{F_0}$; then ${\mathcal P}_I(\O_{F_0})=P_I$;  the group scheme ${\mathcal P}_I$ has $GU(V,\phi)$ as its generic fiber. Then ${\mathcal P}_I\times_{\O_{F_0}}R$ acts on $M^{\rm naive}_I$. We have forgetful
morphisms of projective schemes
\begin{equation}\label{11.1}
 M^{\rm naive} \xrightarrow{\ }   M_I^{\rm
naive}\ ,\ \mbox{and}\ \   M_{I'}^{\rm
naive}\xrightarrow{\ }  M_I^{\rm naive}\ \mbox{for}\
I'\supset I\ \ .
\end{equation}

\subsubsection{}\label{1e3}
Note that the map
$
\{\F_j\}_j\mapsto {\rm ker}(\pi-\sqrt{\pi_0}\ |\ \F_j)
$
gives an isomorphism   between $M^{\rm naive}_I\otimes_{\O_E}F$ and the
Grassmannian ${\rm Gr}(s, n)_F$ of $s$-dimensional spaces in the $n$-dimensional space
$V_0={\rm ker}(\pi-\sqrt{\pi_0}\ |\ V\otimes_{F_0}F)$.

\subsubsection{} As is explained in [RZ], [P1], [PR2], when $C'_p=P_I$, for $F_0=\Q_p$ and $F=K_p$, the (naive) local model is connected to the moduli scheme
via a diagram
\begin{equation}
\A_{C'}\xleftarrow{\ \pi\ } \tilde \A_{C'}\xrightarrow{\ \phi\ } M^{\rm naive}_I \ ,
\end{equation}
where the morphism $\pi$ is a ${\cal P}_I\otimes_{\Z_p}R$-torsor and $\phi$ is a smooth morphism of relative dimension ${\rm dim}(G)$ which is ${\cal P}_I\otimes_{\Z_p}R$-equivariant.
Therefore, there is a relatively representable  morphism of algebraic stacks
\begin{equation}
\A_{C'}\xrightarrow{ \ } [M^{\rm naive}_I/({\cal P}_I\otimes_{\Z_p}R)]
\end{equation}
which is smooth of relative dimension ${\rm dim}(G)$. (See [P1], [PR2] \S 15 for some more details.)

\subsubsection{} As was observed in [P1], the schemes $M^{\rm naive}_I$ are almost never flat over $R$;
by the above, the same is true for $\A_{C'}$. In the present paper, we will examine the flat closure $\A^{\rm flat}_{C'}$ of 
${\rm Sh}_{C'}(G, X)_{E_w}=\A_{C'}\otimes_RE_w$ in $\A_{C'}$. By the above,  we can 
reduce certain local questions about $\A^{\rm flat}_{C'}$ to similar questions
about the flat closure of $M^{\rm naive}_I\otimes_R E_w$ in $M^{\rm naive}$. This last flat closure is,  
by definition, the {\sl local model} $M^{\rm loc}_I$. We can see that $M^{\rm loc}_I$ supports an action of 
${\cal P}_I\otimes_{\Z_p}R$ and there is a relatively representable smooth morphism of relative dimension ${\rm dim}(G)$, 
\begin{equation}
\A^{\rm flat}_{C'}\xrightarrow{ \ } [M^{\rm loc}_I/({\cal P}_I\otimes_{\Z_p}R)]\, .
\end{equation}
This of course implies that each closed point  of $\A^{\rm flat}_{C'}$ 
has an \'etale  neighborhood which is isomorphic to an \'etale neighborhood of 
some corresponding point of $M^{\rm loc}_I$. 

Since the generic fibers
$M_I^{\rm naive}\otimes_R{E_w}$ are, for the various choices of $I$,  all
identical, we obtain for $I'\supset I$ commutative diagrams
of projective morphisms, resp.\ closed embeddings,
$$
\begin{matrix}
M_{I'}^{\rm loc} & \hookrightarrow & 
M_{I'}^{\rm naive}
\\
\big\downarrow && \big\downarrow
\\
M_I^{\rm loc} & \hookrightarrow & 
M_I^{\rm naive}
\end{matrix}
\  ,
$$
whose generic fibers are all the identity morphism. 
We can see that these correspond 
to similar diagrams between the schemes $\A^{\rm flat}_{C'}$ and $\A_{C'}$
for the choices $C'_p=P_{I'}$ or $P_I$.

\subsubsection{}\label{1e4} Let again $F/F_0$ be as in \S\ref{parahoric}. We can also define a subfunctor 
$M^{\wedge}_I$ of $M^{\rm naive}_I$ by specifying that a point of $M^\wedge_I$ with values in an $\O_E$-scheme $S$ is given by
an $\O_F\otimes_{\O_{F_0}}\O_S$-submodule $\F_j\subset \Lambda_j\otimes _{\O_{F_0}}\O_S$ for each $j=k\cdot n\pm i$, $i\in I$,
that in addition to the conditions (a)--(d) above also satisfies:
\smallskip

e) If $r\neq s$, for each $j=k\cdot n\pm i$, $i\in I$, we have 
\begin{equation}
\wedge^{r+1}(\pi -\sqrt{\pi_0}\ |\ \F_j)=0\ ,
\end{equation}
\begin{equation}
\wedge^{s+1}(\pi+\sqrt{\pi_0} \ |\ \F_j)=0\ ,
\end{equation}
where we have  set $\pi=\pi\otimes 1$, $\sqrt{\pi_0}=1\otimes\pi\in \O_F\otimes_{\O_{F_0}}\O_E=\O_F\otimes_{\O_{F_0}}\O_F$.
\smallskip

Denote also by $M^{\wedge}_I$ the corresponding moduli scheme which is a closed subscheme of $M^{\rm naive}_I$.
Suppose that $S$ is an $E$-scheme. Then conditions (d) and (e) are  equivalent since the action of $\pi$ on 
$\F_j$ is semisimple. Therefore, the generic fibers $M^{\wedge}_I\otimes_{\O_E}E$ and $M^{\rm naive}_I\otimes_{\O_E}E$ agree.
It turns out that $M^\wedge_I$ is in some cases flat; then $M^{\rm loc}_I=M^\wedge_I$. This was shown when $I=\{0\}$
and $n= 3$ in [P1]; when $I=\{0\}$ it is conjectured in [P1] that $M^\wedge_{\{0\}}$ is flat. However, 
$M^{\wedge}_{I}$ is not flat in general (see Remark \ref{nonflat}.)

\smallskip
In the sequel we will pursue two goals. The first goal is to understand the structure of $M^{\rm loc}_I$. The second goal is to define $M^{\rm loc}_I$ as a closed subscheme of $M^{\rm naive}_I$ by imposing conditions similar to those defining $M^\wedge_I$. 
\medskip

\section{Affine Weyl groups and affine flag varieties; the $\mu$-admissible set.}
\setcounter{equation}{0}

In this section, $G$ will be a (connected) reductive group over a local field $L$
with perfect residue field.  We assume that $G$ is residually split and hence quasi-split
([T1] 1.10). 
 
\subsection{Affine Weyl groups} We start by recalling some  facts on affine Weyl groups ([HR], [R]).
\subsubsection{} Let $S$ be a maximal split torus in $G$ and let $T$ be its centralizer.
Since $G$ is quasi-split, $T$ is a
maximal torus in $G$. Let $N=N(T)$ be the normalizer of $T$; denote by
$T(L)_1 $ the kernel of the Kottwitz homomorphism
$\kappa_T: T(L)\to X_*(T)_I$; then $T(L)_1={\mathcal T}^0(\O_L)$
where ${\mathcal T}^0$ the connected Neron model of the torus $T$ over $\O_L$.
By definition, the {\sl Iwahori-Weyl group associated to $S$} is the
quotient group
\begin{equation*}
\ti W=N(L)/T(L)_1\ .
\end{equation*}
Since $\kappa_T$ is surjective, the Iwahori-Weyl group $\ti W$
is an extension of the relative Weyl group $W_0=N(L)/T(L)$ by $X_*(T)_I$:
\begin{equation}\label{exactWeyl}
0\to X_*(T)_I\to \ti W\to W_0\to 1 .
\end{equation}
We have ([HR], comp.~also [R], \S 2)

\begin{prop}\label{relpos}
Let $B_0$ be the Iwahori subgroup of $G(L)$ associated to an alcove contained
in the apartment associated to the maximal split torus $S$. Then $G(L)=B_0\cdot N(L)\cdot B_0$
and the map $B_0\cdot n\cdot B_0\mapsto n\in \ti W$ induces a bijection
\begin{equation*}
B_0\backslash G(L)/B_0\xrightarrow{\sim} \ti W\ .
\end{equation*}
If $P$ is the parahoric subgroup of $G(L)$ associated to a facet contained in the apartment corresponding to $S$, then
\begin{equation}
P \backslash G(L)/P \xrightarrow{\sim}  W^P\backslash \ti W/ W^P, \ \ \hbox{\rm where\ \ \ }  W^P:=(N(L)\cap P)/T(L)_1\ .
\end{equation}
\end{prop}

In fact, if $P$ is the (special) parahoric subgroup $P_{\und x}$ that corresponds to a special vertex $\und x$ in
the apartment corresponding to $S$, then the subgroup $ W^P\subset \ti W$ maps isomorphically to $W_0$ under
the quotient $\ti W\to W_0$ and the exact sequence (\ref{exactWeyl}) represents the Iwahori-Weyl group as
a semidirect product
\begin{equation}\label{semidirect}
\ti W=W_0\ltimes X_*(T)_I\ ,
\end{equation}
see [HR].

\subsubsection{} Now let $S_{{\rm sc}}$,   $T_{{\rm sc}}$, resp. $N_{{\rm sc}}$ be the inverse images
of $S\cap G_{{\rm der}}$, $T\cap G_{{\rm der}}$, resp. $N\cap G_{{\rm der}}$ in the simply connected covering
$G_{{\rm sc}}$ of the derived group $G_{{\rm der}}$. Then $S_{{\rm sc}}$ is a maximal split torus of $G_{{\rm sc}}$  and $T_{{\rm sc}}$, resp. $N_{{\rm sc}}$ is its centralizer, resp. normalizer. Hence
\begin{equation*}
W_a:=N_{{\rm sc}}(L)/T_{{\rm sc}}(L)_1
\end{equation*}
is the Iwahori-Weyl group of $G_{{\rm sc}}$. This group is also called the {\sl affine Weyl group
associated to $S$} and is a Coxeter group.
Indeed, we can  recover $W_a$ in the following way:
Let $N(L)_1$ be the intersection of $N(L)$ with the kernel $G(L)_1$ of the Kottwitz homomorphism
$\kappa_G: G(L)\to \pi_1(G)_I$. Then one can see ([HR]) that the natural homomorphism
\begin{equation}
W_a=N_{{\rm sc}}(L)/T_{{\rm sc}}(L)_1\xrightarrow{\sim} N(L)_1/T(L)_1
\end{equation}
is an isomorphism and that there is  an exact sequence
\begin{equation}
1\to W_a\to \ti W\xrightarrow{\kappa_G}  \pi_1(G)_I \to 1\ ,
\end{equation}
where $\pi(G)_I=X_*(T)_I/X_*(T_{{\rm sc}})_I$. 
Now let $B_0$ be the Iwahori subgroup of $G(L)$ associated to an alcove  $C$ in the apartment corresponding to $S$
and let ${\bf S}$ be the set of reflections about the walls of $C$. Then by [BTII], 5.2.12 the quadruple
$(G(L)_1, B_0, N(L)_1, {\bf S})$ is a double Tits system and $W_a=N(L)_1/T(L)_1$
is the affine Weyl group of the affine root system $\Phi_a$ of $S$. The affine Weyl group $W_a$
acts simply transitively on the set of alcoves in the apartment of $S$. Since $\ti W$ acts transitively
on the set of these chambers, $\ti W$ is the semi-direct product of $W_a$ with the normalizer $\Omega$ of
the base alcove $C$, i.e., the subgroup of $\ti W$ which preserves the alcove,
\begin{equation}\label{semi}
\ti W=W_a\rtimes \Omega\ .
\end{equation}
We can identify $\Omega$ with $\pi_1(G)_I$.

Let us write
${\bf S}=\{s_i\}_{i\in I}\subset W_a$ for the finite set of reflections about the walls of
$C$ that generate the Coxeter group $W_a$.
For each $w\in W_a$
its length $l(w)$ is the minimal number of factors in a product of $s_i$'s representing $w$.
Any such product realizing the minimum is called a reduced decomposition of $w$.
We will denote by $\leq $ the corresponding Bruhat order. Recall its definition.  Fix a reduced decomposition
of $w\in W_a$. The elements $w'\leq w$
are obtained by replacing some factors in it by $1$. (This set of such $w'$'s
is independent of the choice of the reduced decomposition of $w$.) We extend the Bruhat order from $W_a$ to $\ti W$ using the semi-direct product decomposition (\ref{semi}): for $w=w_1\cdot \tau$, with $w_1\in W_a$ and $\tau \in \pi_1(G)_I$, the elements smaller than $w$ are the $w'$ of the form $w'=w'_1\cdot \tau$ with $w'_1 \leq w_1 \in W_a$.

Let us denote by $\alpha_i\in \Phi_a$ the
unique affine root with  corresponding affine reflection
equal to $s_i$ (since the group is residually split $\frac{1}{2}\alpha_i\notin \Phi_a$, cf. [T1], 1.8). We will denote by $\Delta=\Delta_G$ the (local) Dynkin diagram of the affine root system
$\Phi_a$ (this can be obtained from the set $\{\alpha_i\}_{i\in I}$;
see [BTI], 1.4 and [T1], 1.8). 
 For  a  subset $Y\subset {\bf S}$, we denote by $W_Y\subset W_a$ the
subgroup generated by $s_i$ with $i\in Y$; we set $P_Y=B_0\cdot W_Y\cdot B_0 $.
By  general properties of Tits systems these are subgroups of $G(L)_1\subset G(L)$; by [BTII], 5.2.12 (i)
they are the parahoric subgroups of $G(L)$ that contain $B$. Using [BTI], 1.3.5 we see that
we can identify $P_Y$ with the parahoric subgroup $P_{C_Y}$, where $C_Y$
is the facet consisting of $a\in \overline C$ for which $Y$ is exactly
the set of reflections $s\in {\bf S}$ which fix $a$.

Finally, let us recall that there exists a reduced root system $\Sigma$ such that the semi-direct
product (\ref{semidirect}) (for $G_{{\rm sc}}$ instead of $G$) presents $W_a$ as the affine Weyl group
associated (in the sense of Bourbaki) to $\Sigma$,
\begin{equation}\label{semidirect2}
W_a=W(\Sigma)\ltimes Q^\vee (\Sigma) ,
\end{equation}
(cf. [BTI],~1.3.8,  [T],~1.7, 1.9). In other words, we have   identifications $W_0\simeq W(\Sigma)$,
$X_*(T_{{\rm sc}})_I\simeq Q^\vee(\Sigma)$ compatible with the semidirect product
decompositions (\ref{semidirect}) and (\ref{semidirect2}).

\subsection{The $\mu$-admissible set} We next recall the definition of the $\mu$-admissible set. 
\subsubsection{}
 To $\mu\in X_\ast(T)$ we attach its image $\lambda$ in the
coinvariants $X_\ast(T)_I$. By (\ref{exactWeyl}) we can consider $\lambda$ as an element in the Iwahori-Weyl group $\tilde{W}$ of $G$. The {\it admissible subset of $\tilde{W}$ associated to the coweight} $\mu$ is defined as 
\begin{equation}
{\rm Adm} (\mu) = \{ w\in\tilde{W}\mid  w\leq w_0(\lambda)\, \text{ for some}\, w_0\in W_0 \}.
\end{equation}
Note that all elements of ${\rm Adm} (\mu)$ have the same image in $\tilde{W}/W_a = \Omega$, namely the image of $\mu$ in $\pi_1(G)_I$.
Furthermore, the set ${\rm Adm} (\mu)$ only depends on the geometric conjugacy class of the one-parameter
subgroup $\mu$, cf.~[R], \S 3.

\subsubsection{}\label{admissible}
In the present paper we find it more convenient to pass to the adjoint group. Let $T_{\rm ad}$ denote the image of $T$ in the adjoint group $G_{\rm ad}$. To $\mu$ we associate
its image $\mu_{\rm ad}$ in $X_\ast(T_{\rm ad})$ and its image $\lambda_{\rm ad}$ in
$X_\ast(T_{\rm ad})_I$.  Now the set ${\rm Adm} (\mu)$ is mapped bijectively to ${\rm Adm} (\mu_{\rm ad})$ under the homomorphism $\ti W\to \ti W_{\rm ad}$, so that it suffices to consider the latter set, which we sometimes also denote by ${\rm Adm} (\mu)$. 
We have a commutative diagram
\begin{equation}
\begin{matrix}
W_a &\xrightarrow{\sim} &W(\Sigma)\ltimes Q^\vee(\Sigma) \\
\downarrow && \downarrow \\
\tilde{W}_{\rm ad} & \xrightarrow{} & W(\Sigma)\ltimes P^\vee(\Sigma)
\end{matrix}
\end{equation}
Here all arrows are injective and $Q^\vee(\Sigma)$, resp.
$P^\vee(\Sigma)$, is the group of coroots, resp. of coweights, of the finite root system $\Sigma$, cf.~[R], \S 3.
We have $\lambda_{\rm ad} \in P^\vee(\Sigma)$, cf. the proof of Lemma 3.1 in [R]. 
In the sequel we write $W_0$ for the finite Weyl group $W(\Sigma)$.

 Denoting by $\tau_{\rm ad}$ the common image of all elements of 
${\rm Adm} (\mu)$ or ${\rm Adm} (\mu_{\rm ad})$ in $\Omega_{\rm ad}=\tilde{W}_{\rm ad}/W_a$, we can define the subset ${\rm Adm} (\mu)^{\circ}$ of
$W_a$ by 
$${\rm Adm}^\circ (\mu)=\{ w \in W_a\mid w \cdot \tau_{\rm ad} \in {\rm Adm} (\mu_{\rm ad}) \}.$$

For a non-empty subset $Y$ of the set of simple affine roots, let $Y^{\circ}\subset {\bf S}$ be the subset that corresponds 
to the set of simple reflections of the form $\{ \tau_{\rm ad}\cdot s_i\cdot \tau^{-1}_{\rm ad}\ |\ i\in Y\}$ where  $s_i$ is the 
reflection corresponding to the simple affine root parametrized by $i\in Y$.

We may define the subset  ${\rm Adm}^Y (\mu)$ of $\ti W_{\rm ad}$ by 
$W^Y\cdot {\rm Adm} (\mu_{\rm ad})\cdot W^{Y}$ and
 the subset  ${\rm Adm}^Y (\mu)^{\circ}$ of
$W_a$ by
$$
{\rm Adm}^Y (\mu)^{\circ}:=W^Y\cdot {\rm Adm} (\mu)^{\circ}\cdot W^{Y^{\circ}}\ .
$$
Since $\tau_{\rm ad}\cdot W^Y\cdot \tau^{-1}_{\rm ad}=W^{Y^\circ}$ this is also equal to $
(W^Y\cdot {\rm Adm} (\mu_{\rm ad})\cdot W^Y)\cdot \tau^{-1}_{\rm ad}$.

\subsection{Calculation of the affine root system}

The following recipe for obtaining the affine root system $\Phi_a$ was explained 
to us by Waldspurger. We may assume for this that the derived group $G_{\der}$ is
simply connected. Let $L' / L$ be a finite Galois extension which splits $G$. Since $G$
is residually split, the extension $L'/L$ is totally ramified. We
extend the normalized valuation $\val : L\to\mZ$ to a valuation 
$\val : L'\to\frac{1}{e}\mZ$, where $e = |L':L|$ is the ramification degree of
$L' / L$. We identify $X_\ast(S)$ with $X_\ast (T)^I$ and $X^{\ast} (S)$ with 
$X^\ast (T)_I$/torsion, and denote by $\varphi$ the natural map,
\begin{equation}
\varphi : X^{\ast} (T)\to X^{\ast} (T)_I \text{/torsion} = X^{\ast} (S) \ .
\end{equation}
Let $\Phi^{\abs} \subset X^{\ast} (T)$ be the set of (absolute) roots. The
relative root system $\Phi$ is the image of $\Phi^{\abs}$ under $\varphi$. Let 
$\Phi_a$ be the affine root system, which is a set of affine functions on 
$V = X_\ast (S)\otimes\mR$ of the form $\beta + r$, for $\beta\in\Phi$ and 
$r\in\mR$. There is a unique action of $N(L)$ on $V$ by affine transformations
such that the elements $t\in T(L)$ act by translations by $\ord (t)\in V$.
Here $\ord (t)$ is uniquely defined by
\begin{equation*}
\langle\ord (t), \chi\rangle = -\val\chi (t) , \quad \chi\in X^\ast (T) \ .
\end{equation*}
Let
\begin{equation}
l_\beta = |\varphi^{-1}(\beta)\cap\Phi^{\abs} | \ , \ \beta\in\Phi \ .
\end{equation}
Let for $\beta\in\Phi,$
\begin{equation}
R_\beta =  \begin{cases}\ds{ \frac{1}{l_\beta}\mZ\, ,} & \text{if \ }  \ds{\frac{\beta}{2}}\not\in\Phi\\ 
    \\                                          \ds{\frac{1}{2l_\beta}+\frac{1}{l_\beta}\mZ\, ,} & \text{if \ }  \ds{\frac{\beta}{2}}\in\Phi\, .\end{cases}
\end{equation}
\begin{prop}
The affine root system $\Phi_a$ equals
\begin{equation*}
\Phi_a = \{\beta + r\mid\beta\in\Phi, r\in R_\beta\} \ .
\end{equation*}
\end{prop}

\begin{proof}
For any $\alpha\in\Phi^{\abs}$, we denote by $U_\alpha$ the corresponding root
subgroup and by ${\alpha}^{\vee}\in X_\ast (T)$ the corresponding coroot.
Similarly, for $\beta\in\Phi$ we have $U_\beta$ and ${\beta}^{\vee}\in X_\ast (S)$
(in case $2\beta\in\Phi$, we have $U_{2\beta}\subset U_\beta$).

Let $\beta\in\Phi$. For $u\in U_\beta (L)\backslash\{1\}$, the intersection
$U_{-\beta} (L)\cdot u\cdot U_{-\beta} (L)\cap N(L)$ consists of a single element,
denoted by $n(u)$, with image in $W_0$ equal to the reflection $s_\beta$.
The action of $n(u)$ on $V$ is an affine reflection about a hyperplane parallel
to $\Ker (\beta)$. Let $r_\beta (u)\in\mR$ be the real number such that this hyperplane
is defined by the equation
\begin{equation*}
\beta (v) + r_\beta (u) = 0\quad, \quad v\in V \ .
\end{equation*}
We also set $r_\beta (1) = \infty$. For $r\in\mR$, define
\begin{equation}
U_{\beta +r} = \{u\in U_\beta (L)\mid r_\beta (u)\geq r\},\quad U_{\beta + r^+} = \{u\in U_\beta (L)\mid r_\beta (u) > r\} \ .
\end{equation}
Then the affine roots in $\Phi$ are defined to be the affine functions
$\beta +r$ such that
\begin{equation}
U_{\beta +r} / U_{\beta +r^+}\cdot U_{2\beta +2r}\neq \{1\} \ ,
\end{equation}
cf. [T], 1.6. Here we set $U_{2\beta +2r} =\{1\}$, if $2\beta\not\in\Phi$. 

Let
$\beta\in\Phi$. Assume first that $2\beta\not\in\Phi$ and $\frac{\beta}{2}\not\in\Phi$,
(i.e., $\beta$ of type I). There exists a sub-extension $L\subset L_\beta\subset L'$
such that $l_\beta =|L_\beta : L|$ and such that the subgroup of $G(L)$ generated
by $U_\beta (L)$ and $U_{-\beta} (L)$ is isomorphic to $SL_2 (L_\beta)$ (here the
fact that $G_{\der}$  is simply connected, enters). We choose this isomorphism such
that $U_\beta (L)$ becomes identified with the upper unipotent matrices, and
$U_{-\beta} (L)$ with the lower unipotent matrices and the intersection of $S(L)$
with this subgroup with the diagonal torus.  For $b\in L_\beta\backslash\{0\}$
and $u =  \begin{pmatrix} 1 & b\\   0 & 1\end{pmatrix}$, one calculates ([BT II], p. 80)
\begin{equation*}
n(u) = \ \begin{pmatrix} 1 & 0\\ -b^{-1} & 1\end{pmatrix}\cdot \begin{pmatrix} 1 & b\\ 0 & 1\end{pmatrix}\cdot \begin{pmatrix} 1 & 0\\ -b^{-1} & 1\end{pmatrix} =  \begin{pmatrix} b & 0\\ 0 & b^{-1}\end{pmatrix}\cdot w \  ,  \ \ w =\begin{pmatrix}0&1\\-1&0\end{pmatrix}.
\end{equation*}
The action of $n(u)$ on $V$ is given by
\begin{equation}
v\mto s_\beta (v) - \val (b)\cdot\beta^\vee (v) \ .
\end{equation}
Here $s_\beta$ denotes the reflection about the hyperplane $\Ker\beta$. It follows that
$r_\beta (u) = \val (b)$. Hence $\beta + r\in\Phi$ if and only if $r\in\val (L_\beta)$. Since
$L_\beta / L$ is totally ramified, $\val (L_\beta) = \frac{1}{l_\beta}\cdot\mZ$. This proves
the claim for affine functions of type $v\mto\beta (v)+r$, for $\beta$ of Type I.

Now let $\beta\in\Phi$ be of type II, i.e., such that $2\beta\in\Phi$ or
$\frac{\beta}{2}\in\Phi$. Then there is a tower of subextensions
$L\subset L_\beta\subset L'_\beta\subset L'$ such that $|L'_\beta : L_\beta| = 2$ and such
that
\begin{equation*}
|L_\beta : L| =  \begin{cases} \displaystyle{\frac{l_\beta}{2}},  &\mbox{if }  2\beta\in\Phi\\ 
                                                            l_\beta, & \mbox{if }  \displaystyle{\frac{\beta}{2}}\in\Phi\, .  \end{cases}\ 
\end{equation*}
The subgroup of $G(L)$ generated by $U_\beta (L)$ and $U_{\frac{\beta}{2}} (L)$
($=\{1\}$ if $\frac{\beta}{2}\not\in\Phi$) and $U_{-\beta} (L)$ and $U_{-\frac{\beta}{2}} (L)$
is isomorphic to $SU_3 (L_\beta)$, the special unitary group relative to the quadratic
extension $L'_\beta / L_\beta$ and the hermitian form with anti-diagonal unit matrix.
We again choose this isomorphism in the evident way. Let $\beta_0$ resp. $\beta_0^\vee$
the root resp. coroot of $SU_3 (L_\beta)$, given by
\begin{equation*}
 \begin{pmatrix} t & \quad & \quad\\ 
                                  \quad & 1 & \quad\\
                                  \quad & \quad & t^{-1}\end{pmatrix} \mto t\quad,\quad t\mto
 \begin{pmatrix} t^2 & \quad & \quad\\ 
                                  \quad & 1 & \quad\\ \quad & \quad & t^{-2}\end{pmatrix}\, .
\end{equation*}
The identification is chosen such that if $2\beta\in\Phi$, then $U_\beta$ is identified
with
\begin{equation*}
U_{\beta_0} = \left\{\begin{pmatrix} 1 & b & c\\ 0 & 1 & \overline{b}\\
0 & 0 & 1\end{pmatrix}\
 \mid\ c+\overline{c} = b \ \overline{b}\right\}\, ,
\end{equation*}
and if $\frac{\beta}{2}\in\Phi$, then $U_\beta$ is identified with
\begin{equation*}
U_{2 \beta_0} = \left\{\begin{pmatrix} 1 & 0 & c\\ 0 & 1 & 0\\0 & 0 & 1\end{pmatrix} \ \mid\ c+\overline{c} = 0 \right\}\, .
\end{equation*}
Here $x\mapsto\overline{x}$ denotes Galois automorphism of $L'_\beta / L_\beta$. For
$u\in U_{\beta_0}$, written as above, one has, [BT II], p. 84,
\begin{equation*}
n(u) = \begin{pmatrix} c & \quad & \quad\\ \quad &   \overline{c}{c^{-1}} & \quad\\
                                              \quad & \quad & \overline{c}^{-1}\end{pmatrix} w \,,
\end{equation*}
where
\begin{equation*}
w =  \begin{pmatrix} \quad & \quad & 1\\ \quad & -1 & \quad\\
                                           1 & \quad & \quad\end{pmatrix}\,.
\end{equation*}
Hence the action of $n(u)$ on $V$ is given by
\begin{equation}
v\mto s_\beta (v) - \frac{\val (c)}{2}\cdot{\beta}^{\vee}_0 \ .
\end{equation}
It follows that
\begin{equation}
r_\beta (u) =  \begin{cases} \ds{\frac{\val (c)}{2}}\,, & \mbox{\ if\ } \  2\beta\in\Phi\\ 
                                                       \val (c)\,, & \mbox{\ if\ }  \ds{\frac{\beta}{2}}\in\Phi\, .   \end{cases}
\end{equation}
Hence the condition that $U_{\beta +r} / U_{\beta +r^+}\cdot U_{2\beta +2r}\neq \{1\}$ is
equivalent to

\begin{equation}
 \begin{cases} \ \exists\ b \in L'_\beta \  \text{with}\   r = \frac{1}{2}\displaystyle{\rm max}\{\val (c)\mid c+\overline c=b\overline b\} , &
                                                                        \   \text{ \, if }2\beta\in\Phi\\ 
                                   \ \exists\ c \in L'_\beta \  \mbox{with}\  c+\overline{c} = 0, \text{\ and } r = \val (c)\,, & \ \text{ \,  if } \displaystyle{\frac{\beta}{2}}\in\Phi\,. \end{cases}
\end{equation}
Since $p\neq 2$, and $L' / L$ is totally ramified, this is equivalent to
\begin{equation}
 \begin{cases}\ r\in\displaystyle{\frac{1}{l_\beta}\mZ}\,, & \text{\ if }\,  2\beta\in\Phi\\ 
                          \ \displaystyle{ r\in\frac{1}{2l_\beta} + \frac{1}{l_\beta}\mZ}\,, & \text{\ if }\,   \displaystyle{\frac{\beta}{2}}\in\Phi\,. \end{cases}
\end{equation}
\end{proof}

\subsection{The  case of a unitary group}    We now take $G$ to be the (quasi-split) unitary similitude group 
over $L=F_0$ corresponding to a tame ramified 
quadratic extension $F/F_0$ (as in \ref{parahoric}). 
Set $\Gamma={\rm Gal}(F/F_0)=\langle \tau\rangle$
and notice that we can take $L'=F$ in the notation of the previous paragraph.
As in \S \ref{oddeven} we denote by $T$ the standard maximal torus of $G$ and by $S$ the maximal split subtorus of $T$. 
Our first purpose is to calculate explicitly the inclusion
$$
W_a\subset\tilde W_{\rm ad}\ .
$$
Here $W_a\cong X_*(T_{\rm sc})_\Gamma\rtimes W_0$ 
(as in (\ref{semidirect}), after the choice of a special vertex) 
is the Coxeter group generated by the reflections in the root hyperplanes in
$\Phi_a$. The group $\tilde W_{\rm ad}$ is obtained from $W_a$ by the
push-out via the inclusion $X_*(T_{\rm sc})_\Gamma\hookrightarrow X_*(T_{\rm
ad})_\Gamma$.

Recall the definition of the cocharacter $\mu_{r,s}\in X_*(T)$ in \S \ref{traeq}. Using $G(F)\simeq GL_n(F)\times F^{\times}$ we identify $X_*(T)$ with $\Z^n\times\Z$\,\,; then  $\mu_{r,s}=((1^{(s)}, 0^{(r-s)}), 1)$.
Our second purpose is to determine the image $\lambda=\lambda_{r,s}$ of $\mu=\mu_{r,s}\in
X_*(T)$ in $X_*(T_{\rm ad})_\Gamma$ and the set
$$
{\rm Adm}\ (\mu)=\{ w\in\tilde W_{\rm ad}\mid w\leq w_0(\lambda),\ \mbox{some}\
w_0\in W_0\}\ \ .
$$
Let ${\rm Adm}_0(\mu)$ be the image of $W_0\cdot {\rm Adm}(\mu)$ in
$W_0\setminus \tilde W_{\rm ad}/W_0$. Our third purpose is to determine ${\rm
Adm}_0(\mu)$.

\subsubsection{The case $A_{2m-1}^{(2)}$ (corresponding to $GU_{2m}(F/F_0)$) for $m\geq 2
$}\label{caseEv}

Let $n=2m$. In this case,  $X_\ast(S)$ can be identified with the subgroup of $X_*(T)\simeq\Z^n\times\Z$ formed by the elements of the form  $(x_1,\ldots,x_m,-x_m,\ldots, -x_1;y)$. The relative roots $\beta$ are of the form $\pm x_i\pm x_j$ for $i\neq j$ and of the form $\pm 2x_i$, for $i, j\in \{1,\ldots , m\}$. For the first kind  of $\beta$ we have $l_\beta=2$, and for the second kind $l_\beta=1$. Therefore the affine root system $\Phi_a$ consists of the affine functions $\pm x_i\pm x_j +\frac{1}{2} \Z$ and $\pm 2x_i+\Z$. Hence the set of root hyperplanes is the zero sets of the affine functions 
$\{\pm x_i\pm x_j+\frac{1}{2}\Z;\  \pm x_i+\frac{1}{2}\Z\}.$  It follows 
that  the  root system $\Sigma$  is of type $B_m$ (i.e. $\Sigma ^\vee$ is of type $C_m$) and we have
$$
W_a= Q^\vee\rtimes W_0\subset\tilde W_{\rm ad}=P^\vee\rtimes W_0\  ,
$$
with
\begin{eqnarray*}
W_0
&=&
S_m\rtimes \{ \pm 1\}^m
\\
Q^\vee
&=&
\{ x\in\mathbb Z^m\mid \Sigma(x)\equiv 0(2)\}
\\
P^\vee
&=&
\mathbb Z^m\ \ .
\end{eqnarray*}
The map
$$
X_*(T)\to X_*(T_{\rm ad})_\Gamma\buildrel\sim\over\longrightarrow \mathbb Z^m
$$
is given by
$$
(x_1,\ldots, x_{2m}, y)\longmapsto (x_1-x_{2m}, x_2-x_{2m-1},\ldots, x_m-x_{m+1})\ \
.
$$
The image $\lambda=\lambda_s$ of $\mu_{r, s}=(1^{(s)}, 0^{(r-s)}, 1)$ in $X_*(T_{\rm ad})_\Gamma=\mathbb Z^m$ is equal to
$(1^{(s)},0^{(m-s)})$.

Let us identify $W_0\setminus \tilde W_{\rm ad} / W_0=P^\vee / W_0$ with
$$P^\vee\cap
\bar C=\{ x\in\mathbb Z^m\mid x_1\geq x_2\geq \ldots \geq x_m\geq 0\}\ \ .$$
The induced order from the Bruhat order is the dominance order. The positive coroots
are (Bourbaki, table $C_m$)
$$e_i-e_j,\ e_i+e_j,\ 2e_i\quad (i<j)\ \ .$$
The only possibility of subtracting a positive coroot from $\lambda_s$ and stay inside
$\bar C$ is by $e_{s-1}+e_s$. Hence we see inductively that
$${\rm Adm}_0(\mu_{r, s})=\{ \lambda=\lambda_s> \lambda_{s-2}>\ldots \}\ \ ,$$
with the chain ending in $\lambda_1$ if $s$ is odd and in $\lambda_0$ if $s$ is even.

Note that $\lambda_s$ is minuscule if and only if $s=0$ or $s=1$.

\subsubsection{The case $A_{2m}^{(2)}$ (corresponding to $GU_{2m+1}(F/F_0)$) for $m\geq 1$}\label{caseOdd}
Set $n=2m+1$. In this case,  $X_\ast(S)$ can be identified with the subgroup of $X_*(T)\simeq\Z^n\times\Z$ formed by the elements of the form  $(x_1,\ldots, x_m, 0, -x_m,\ldots, -x_1;y)$. The relative roots $\beta$ are of the form $\pm x_i\pm x_j$ for $i\neq j$ and of the form $\pm x_i$ and of the form $\pm 2x_i$, for $i, j\in \{1,\ldots , m\}$. For the first and second kind  of $\beta$ we have $l_\beta=2$, and for the third kind $l_\beta=1$. Therefore the affine root system $\Phi_a$ consists of the affine functions $\pm x_i\pm x_j +\frac{1}{2} \Z$ and $\pm x_i+\frac{1}{2}\Z$ and $\pm 2x_i+\frac{1}{2}+\Z$. 
Hence the set of root hyperplanes is the zero sets of the affine functions 
$\{\pm x_i\pm x_j+\frac{1}{2}\Z;\  \pm 2x_i+\frac{1}{2}\Z\}.$  In this case, the  root system $\Sigma$ is of type $C_m$, and we obtain
$$W_a=Q^\vee \rtimes W_0=\tilde W_{\rm ad}\ \ .$$
Here
\begin{eqnarray*}
W_0
&=&
S_m\rtimes \{ \pm 1\}^m
\\
Q^\vee
&=&
\mathbb Z^m
\end{eqnarray*}
(whereas $P^\vee =\mathbb Z^m+\mathbb Z\left( \frac{1}{2} \sum\limits_{i=1}^m
e_i\right)$).
The map
$$
X_*(T)\to X_*(T_{\rm ad})_\Gamma\buildrel\sim\over\longrightarrow \mathbb Z^m$$
is given by
$$(x_1,\ldots, x_{2m+1}, y)\longmapsto (x_1-x_{2m+1}, x_2-x_{2m},\ldots, x_m-x_{m+2})\ \
.$$
The image of $\mu_{r,s}$ is $\lambda_s= (1^{(s)},0^{(m-s)})$. Again we identify
$W_0\setminus W_a /W_0$ with
$$Q^\vee\cap \bar C=\{ x\in \mathbb Z^m\mid x_1\geq x_2\geq \ldots \geq x_m\geq 0\}\
\ .$$
In this case the positive coroots are
$$e_i - e_j\ ,\ e_i+e_j\ ,\ e_i\quad (i<j)\ \ ,$$
and
$${\rm Adm}_0(\mu_{r, s})= \{\lambda= \lambda_s > \lambda_{s-1} >\ldots > \lambda_1
>\lambda_0=0\}\ \ .$$
Note that $\lambda_s$ is minuscule if and only if $s=0$.

\subsection{Cases of small rank} Below we give the complete information for the cases of small rank.  Note that here, in order to simplify the notation,  we are switching to the unitary group from the unitary similitude group; the group over $F$ is $GL_n$ etc.  

\subsubsection{The case $A_3^{(2)}$ (corresponds to $U_4(F/F_0)$)}
\begin{itemize}
\item[a)] {\bf Apartment over $F$:} $(x_1,x_2,x_3,x_4)\in \mathbb R^4$
\item[] {\bf Affine roots:} $\alpha+\frac{1}{2}\mathbb Z$ 
\item[] {\bf Fundamental alcove:} $x_4+\frac{1}{2} \geq x_1\geq x_2\geq x_3\geq x_4$
\item[] {\bf Simple affine roots:} $\alpha_1,\alpha_2,\alpha_3$ and $\alpha_0(x)=
x_4-x_1+\frac{1}{2} =-\theta +\frac{1}{2}$ with $\theta=\alpha_1+\alpha_2+\alpha_3$

\medskip
\item[b)] {\bf Action of $\tau$ on apartment:} $(x_1,x_2,x_3,x_4)\mapsto
(-x_4,-x_3,-x_2,-x_1)$
\item[] {\bf Apartment over $F_0$:} $(x_1,x_2,-x_2,-x_1)$
\item[] {\bf Positive roots:}
$\alpha_1,\alpha_2,\alpha_3,\alpha_1+\alpha_2,\alpha_2+\alpha_3, \theta$
\item[] {\bf Relative roots:}
$$\begin{matrix}
{\rm res}(\alpha_1)={\rm res}(\alpha_3)\hfill
&:&
x_1-x_2\hfill \quad (l_\beta=2)
\\
{\rm res}(\alpha_1+\alpha_2) ={\rm res}(\alpha_2+\alpha_3)\hfill
&:&
x_1+x_2 \hfill \quad (l_\beta=2)
\\
{\rm res}(\alpha_2)\hfill
&:&
2x_2 \hfill \quad (l_\beta=1)
\\
{\rm res}(\theta)\hfill
&:&
2x_1 \hfill \quad (l_\beta=1)
\end{matrix}$$
\item[] {\bf Affine relative roots:}
$$\begin{matrix}
x_1-x_2+\frac{1}{2}\mathbb Z\hfill
&&
\hfill
\\
x_1+x_2+\frac{1}{2}\mathbb Z\hfill
&&
\hfill
\\
2x_2+\mathbb Z\hfill
&&
\hfill
\\
2x_1+\mathbb Z\hfill
&&
\hfill
\end{matrix}$$
\item[] {\bf Fundamental alcove:} $\frac{1}{2} -x_2\geq x_1\geq x_2\geq 0$
\item[] {\bf Simple affine roots:}
$$\begin{matrix}
x_1-x_2\hfill
&:&
{\rm res}(\alpha_1)= {\rm
res}(\alpha_3)\hfill
\\
2x_2\hfill
&:&
{\rm res}(\alpha_2)\hfill
\\
\frac{1}{2}-x_1-x_2\hfill
&:&
-{\rm res}(\alpha_1+\alpha_2)+\frac{1}{2}\hfill
\end{matrix}$$

\medskip
\item[c)] $X_*(T_{\rm sc})= \{ x\in\mathbb Z^4\mid \Sigma(x)=0\}$
\item[] $X_*(T_{\rm sc})_\Gamma= X_*(T_{\rm sc}) / X_*(T_{\rm sc})\cap X_*(T_{\rm
sc})_{\mathbb Q}^-\to \mathbb Z^2$, via $x\mapsto (x_1-x_4, x_2-x_3)$
\item[] (acts by translation by $\frac{1}{2}$ on apartment). Here $X_*(T_{\rm
sc})_{\mathbb Q}^-$ denotes the $-1$-eigenspace. 
\item[] ${\rm Image} = \{ x\in\mathbb Z^2\mid \Sigma(x)=0(2)\}$.
\item[] $X_*(T_{\rm ad})_\Gamma\buildrel\sim\over\longrightarrow \mathbb Z^2$
surjective.
\end{itemize}

\subsubsection{The case $A_2^{(2)}$ (corresponds to $U_3(F/F_0)$)}

\begin{itemize}
\item[a)] {\bf apartment over $F$:} $(x_1,x_2,x_3)\in \mathbb R^3$
\item[] {\bf fundamental alcove:} $x_3+\frac{1}{2}\geq x_1\geq x_2\geq x_3$
\item[] {\bf Simple affine roots:} $\alpha_1,\alpha_2$, and
$\alpha_0(x)=x_3-x_1+\frac{1}{2}=-\theta+\frac{1}{2}$, with
$\theta=\alpha_1+\alpha_2$.

\medskip
\item[b)] {\bf action of $\tau$ on apartment:} $(x_1,x_2,x_3)\mapsto
(-x_3,-x_2,-x_1)$
\item[] {\bf apartment over $F_0$:} $(x_1, 0, -x_1)$
\item[] {\bf Positive roots:} $\alpha_1,\alpha_2,\alpha_1+\alpha_2$
\item[] {\bf Relative roots:}
$$\begin{matrix}
{\rm res}(\alpha_1)={\rm res}(\alpha_2)\hfill
&:&
x_1\hfill \quad (l_\beta=2)
\\{\rm res}(\alpha_1+\alpha_2)\hfill
&:&
2x_1\hfill \quad (l_\beta=1)
\end{matrix}$$
\item[] {\bf Affine relative roots:}
$$\begin{matrix}
x_1+\frac{1}{2}\mathbb Z\hfill
&&
\hfill
\\
2x_1+\frac{1}{2}+\mathbb Z\hfill
&&
\hfill
\end{matrix}$$
\item[] {\bf Fundamental alcove:} $-x_1+\frac{1}{2}\geq x_1\geq 0$
\item[] {\bf Simple affine roots:}
$$\begin{matrix}
x_1\hfill
&:&
{\rm res}(\alpha_1)\hfill
\\
-2x_1+\frac{1}{2}\hfill
&:&
-{\rm res}(\alpha_1+\alpha_2)+\frac{1}{2}\hfill
\end{matrix}$$

\medskip
\item[c)] $X_*(T_{\rm sc})= \{ x\in \mathbb Z^3\mid \Sigma(x)=0\}$
\item[] $X_*(T_{\rm sc})_\Gamma=\{ x\in\mathbb Z^3\mid \Sigma(x)=0\} /\{ x\mid
\Sigma(x)=0, x_1=x_3 \} \buildrel\sim\over\longrightarrow \mathbb Z$, via $x\mapsto x_1-x_3$.
\item[] $X_*(T_{\rm ad})_\Gamma\buildrel\sim\over\longrightarrow\mathbb Z \ .$
\end{itemize}

\subsubsection{The case $A_4^{(2)}$ (corresponds to $U_5(F/F_0)$)}

\begin{itemize}
\item[a)] {\bf apartment over $F$:} $(x_1,x_2,x_3, x_4, x_5)\in \mathbb R^5$
\item[] {\bf fundamental alcove:} $x_5+\frac{1}{2}\geq x_1\geq x_2\geq \ldots \geq x_5$
\item[] {\bf Simple affine roots:} $\alpha_1,\ldots , \alpha_4$, and
$\alpha_0(x)=x_5-x_1+\frac{1}{2}=-\theta+\frac{1}{2}$, with
$\theta=\alpha_1+\ldots +\alpha_4$.

\medskip
\item[b)] {\bf action of $\tau$ on apartment:} $(x_1,\ldots ,x_5)\mapsto
(-x_5, \ldots ,-x_1)$
\item[] {\bf apartment over $F_0$:} $(x_1, x_2,  0, -x_2,  -x_1)$
\item[] {\bf Positive roots:} $\alpha_1,\ldots, \alpha_4,\alpha_1+\alpha_2, \alpha_2+\alpha_3, \alpha_3+\alpha_4, \alpha_1+\alpha_2+\alpha_3, \alpha_2+\alpha_3+\alpha_4, \theta$
\item[] {\bf Relative roots:}
$$\begin{matrix}
{\rm res}(\alpha_1)={\rm res}(\alpha_4)\hfill
&:&
x_1-x_2\hfill \quad (l_\beta=2)
\\{\rm res}(\alpha_1+\alpha_2+\alpha_3)={\rm res}(\alpha_2+\alpha_3+\alpha_4)\hfill
&:&
x_1+x_2\hfill \quad (l_\beta=2)
\\{\rm res}(\alpha_2)={\rm res}(\alpha_3)\hfill
&:&
x_2\hfill \quad (l_\beta=2)
\\{\rm res}(\alpha_1+\alpha_2)={\rm res}(\alpha_3+\alpha_4)\hfill
&:&
x_1\hfill \quad (l_\beta=2)
\\{\rm res}(\alpha_2+\alpha_3)=2\ {\rm res}(\alpha_2)=2\ {\rm res}(\alpha_3)\hfill
&:&
2x_2\hfill \quad (l_\beta=1)
\\{\rm res}(\theta)=2{\rm res}(\alpha_1+\alpha_2)=2{\rm res}(\alpha_3+\alpha_4)\hfill
&:&
2x_1\hfill \quad (l_\beta=1)
\end{matrix}$$
\item[] {\bf Affine relative roots:}
$$\begin{matrix}
x_1-x_2+\frac{1}{2}\mathbb Z\hfill
&&
({\rm mult}\ 2)\hfill
\\
x_1+x_2+\frac{1}{2}\mathbb Z\hfill
&&
({\rm mult}\ 2)\hfill
\\
x_1+\frac{1}{2}\mathbb Z\hfill
&&
({\rm mult}\ 2)\hfill
\\
x_2+\frac{1}{2}\mathbb Z\hfill
&&
({\rm mult}\ 2)\hfill
\\
2x_1+\frac{1}{2}+\mathbb Z\hfill
&&

({\rm mult}\ 1)\hfill
\\
2x_2+\frac{1}{2}+\mathbb Z\hfill
&&
({\rm mult}\ 1)\hfill

\end{matrix}$$
\item[] {\bf Fundamental alcove:} $-x_1+\frac{1}{2}\geq x_1\geq x_2\geq 0$
\item[] {\bf Simple affine roots:}
$$\begin{matrix}
x_1-x_2\hfill
&:&
{\rm res}(\alpha_1)\hfill
\\
x_2\hfill
&:&
{\rm res}(\alpha_2)\hfill
\\
-2x_1+\frac{1}{2}\hfill
&:&
-{\rm res}(\theta)+\frac{1}{2}\hfill
\end{matrix}$$

\medskip
\item[c)] $X_*(T_{\rm sc})= \{ x\in \mathbb Z^5\mid \Sigma(x)=0\}$
\item[] $X_*(T_{\rm sc})_\Gamma=\{ x\in\mathbb Z^5\mid \Sigma(x)=0\} /\{ x\mid
\Sigma(x)=0, x_1=x_5, x_2=x_4 \} \buildrel\sim\over\longrightarrow \mathbb Z^2$, via $x\mapsto (x_1-x_5, x_2-x_4)$.
\item[] $X_*(T_{\rm ad})_\Gamma\buildrel\sim\over\longrightarrow\mathbb Z^2 \ .$
\end{itemize}


\section{Affine Flag varieties}\label{sec3}
\setcounter{equation}{0}

\subsection{Affine flag varieties and the coherence conjecture }\label{3a}

Let $k$ be an algebraically closed field. We denote by $K=k((t))$ the field of Laurent  power series with coefficients in $k$. Let $G$ be a reductive algebraic group over $K$. 
We will  assume for simplicity that the derived group of $G$ is simply connected and absolutely simple, and also that $G$ splits over a tamely ramified extension of $K$.

Let $S$ be a maximal split torus  in $G$. We fix  an alcove $C$ in the apartment corresponding to $S$, and denote by $B$ the associated Iwahori subgroup of $G(K)$. 
Let ${\bf S}$ be the set of reflections about the walls of $C$. The parahoric subgroups containing $B$ correspond to the non-empty subsets $Y$ of ${\bf S}$. More precisely,
we associate to $Y$ the unique parahoric subgroup containing $B$ such that the  associated subgroup $W^Y= W^{P^Y}$ in the sense of (\ref{relpos}) is equal to the subgroup of $W_a$
generated by the simple reflections for the simple affine roots $\alpha_i$ with $i\notin Y$. In particular, $P^{\bf S}=B$. 

Let $LG$ be the loop group associated to $G$, cf.~[PR3]. This is the ind-group scheme over $k$ which represents the functor $R\mapsto G(R((t)))$ on the category of $k$-algebras. To any parahoric subgroup $P$ there is associated a smooth affine group scheme
with connected fibers over $\O_K$ and with generic fiber equal to $G$. We will again denote this group scheme by $P$ and let $L^+P$ be the associated group scheme over $k$. Recall that $L^+P$ represents the functor $R\mapsto P(R[[t]])$ on the category of $k$-algebras. Let $\F(G)=LG/L^+B$, resp. $\F^Y(G)= LG/L^+P^Y$ be the affine flag variety of $G$, resp. the partial affine flag variety of $G$ corresponding to $Y$. These
fpqc quotients are representable by ind-schemes.

By [PR3], \S 10, there is a canonical isomomorphism
\begin{equation}\label{pic}
{\rm Pic}(\F^Y(G_{\rm der}))\xrightarrow {\ \sim\ }\bigoplus\nolimits_{i\in Y} {\mathbb Z}\cdot \epsilon_i \ ,
\end{equation}
given by sending a line bundle $\L$ to the degrees of its restrictions to the projective lines
${\mathbb P}_{\alpha_i}$ corresponding to the simple affine roots  $\alpha_i$, for $i\in Y$.
Put $\kappa(i)=1$, unless the vector part of the simple affine root $\alpha_i$ is a multipliable root, in which case we set $\kappa(i)=2$ (this last case only arises when $G$ is of type $A_{2m}^{(2)}$). Let $\L'(Y)$ be the line bundle on $\F^Y(G_{\rm der})$ which is sent to $\sum_i \kappa(i) \epsilon_i$ under (\ref{pic}). Then $\L'(Y)$ is ample, i.e., the restriction of $\L'(Y)$ to any Schubert variety $S_w$ is ample, for any $w\in W^Y\backslash W_a/W^Y$. 

Let $\mu$ be a cocharacter of $G_{\rm ad}$. 
Recall from \S\ref{admissible} the subset ${\rm Adm}^Y(\mu)^\circ$ of $ W_a $. More precisely, in the notation of loc. cit. we take $I={\bf S}\setminus Y$ and $I^\circ=\tau^{-1}_{\rm ad}(I)={\bf S}\setminus Y^\circ$, with $Y^\circ=\tau^{-1}_{\rm ad}(Y)$; then ${\rm Adm}^Y(\mu)^\circ={\rm Adm}_I(\mu)^\circ$. We define the subset $\A^Y(\mu)^\circ$ of $\F^{Y^\circ}(G_{\rm der})$ as a reduced union of 
$L^+B$-orbits
\begin{equation}
\A^Y(\mu)^\circ =\ \bigcup\nolimits_{w\in {\rm Adm}^Y(\mu)^\circ} L^+B\cdot n_w \ .
\end{equation}
Note that the union of Schubert varieties in  $\F^Y(G_{\rm ad})$ defined by 
$\A^Y(\mu) =\ \bigcup\nolimits_{w\in {\rm Adm}^Y(\mu)} S_w \ $
is the translation by $\tau$ of $\A^Y(\mu)^\circ$. 

We recall the {\it coherence conjecture} from [PR3], \S 10. It concerns the dimensions $h_Y^{(\mu)}=\dim H^0(\A^Y(\mu)^\circ, \L'(Y)^{\otimes k})$, as $Y$ ranges over the non-empty subsets of ${\bf S}$. For any minuscule coweight $\mu$ of $G_{\rm ad}$ we introduce the polynomial
\begin{equation}
h^{(\mu)}(k)= \dim H^0(X(\mu), \L(\mu)^{\otimes ek}) \ .
\end{equation}
Here $e=[K':K]$ is the degree of the splitting field $K'$ of $G$, and $X(\mu)=G_{K'}/P(\mu)$ is the homogeneous projective variety associated to $\mu$, and $\L(\mu)$ is the ample generator of the Picard group of $X(\mu)$. If $\mu= \mu_1+\ldots +\mu_r$ is a sum of minuscule coweights, we set $h^{(\mu)}= h^{(\mu_1)}\cdot \cdots \cdot h^{(\mu_r)}$. The conjecture then states that
\begin{equation}\label{coherence}
h_Y^{(\mu)}(k)= h^{(\mu)}(|Y|\cdot k) \ ,
\end{equation}
provided that $\mu$ is a sum of minuscule coweights for $G_{\rm ad}$. In [PR3] this conjecture is proved for $G=GL_n$ and for $G=GSp_{2n}$. In the present paper we need the conjecture for $G=GU(V,\phi)$, in which case it is an open question.

\subsection{Unitary affine flag varieties}\label{3b}

Assume ${\rm char}(k)\neq 2$.
We consider  the unitary similitude group described in \S \ref{parahoric} in the special case when $F_0=K=k((t))$ with uniformizer $\pi_0=t$, and $F=K'=k((u))$ with uniformizer $u$ with $u^2=t$. We also have the hermitian vector space $(V, \phi)$ of dimension $n\geq 3$ over $K'$, which we assume split, with distinguished basis $e_1,\dots, e_n$. In this case, we will denote the standard lattice chain by $\lambda_0,\dots , \lambda_{n-1}$. As in \S \ref{unimoduli}, we can complete this into a periodic self-dual lattice chain $\lambda$ in $V=K'^n$. 
We denote by $G=GU(V, \phi)$ the group of unitary similitudes, which is an algebraic group over $K$, of the type considered in \S\ref{3a}. The simultaneous stabilizer of all the lattices in  the lattice chain $\lambda$ is an Iwahori subgroup $B$ as before. As in \S\ref{parahoric} we  set $m=[n/2]$. 

Let $I\subset \{0, \ldots , m\}$ be a non-empty subset. We also require as usual  for $n=2m$ that if $m-1\in I$, then also $m\in I$, cf. \S\ref{parahoric}. Let us write $I=\{i_0<i_1<\cdots <i_k\}$. We  consider the part of  the ``standard" lattice chain
\begin{equation*}\label{unichainSt}
\lambda_{i_0}\subset \lambda_{i_1}\subset \cdots \subset \lambda_{i_k}\subset u^{-1}\lambda_{i_0}\ .
\end{equation*}
We can complete this to a periodic self-dual lattice chain $\lambda_I=\{\lambda_j\}_j$ 
in $V=K'^n$ as in \S \ref{unimoduli}. (Recall that here $j$ is of the form 
$k\cdot n\pm i$, $i\in I$, $k\in \Z$, and $\lambda_{k\cdot n-i}=\pi^{-k}\cdot \hat \lambda_i$, $\lambda_{k\cdot n+i}=\pi^{-k}\cdot \lambda_i$.)

Let us now consider the functor $\F_I$ which to a $k$-algebra $R$ associates
the set of pairs of  $R[[u]]$-lattice chains 
 \begin{equation*}\label{chains1}
L_{i_0}\subset L_{i_1}\subset \cdots \subset L_{i_k}\subset u^{-1}L_{i_0}
\end{equation*}
in $V\hat\otimes_{K'}R=R((u))^n$ together with an invertible  element $\alpha\in   R((t))^{\times}$ which satisfy the following conditions:

a) For any $q\in \{0,\ldots, k\}$, we have
\begin{equation*}
L_{i_q}\subset u^{-1}\alpha^{-1}\hat L_{i_q}\subset u^{-1}L_{i_q}\ ,
\end{equation*}

b) The quotients $L_{i_{q+1}}/L_{{i_q}}$, $u^{-1}\alpha^{-1}\hat L_{i_q}/L_{i_q}$, $u^{-1}L_{i_0}/L_{i_{q+1}}$ are projective $R$-modules
of rank equal to the rank of the corresponding quotients for the standard chain  (when $q=k$, these conditions have to be interpreted in the obvious way).
\medskip

Here, we may think of $\alpha$ as a ``similitude'' that modifies the hermitian form 
$\phi\hat\otimes_{K'}R$ to $\phi_\alpha:=\alpha\cdot (\phi\hat\otimes_{K'}R)$.
Notice that the dual of $L_i$ with respect to the new form $\phi_\alpha$ is   $\alpha^{-1}\hat L_i$. 
The ind-group scheme $LG$ over $k$  acts naturally on $\F_I$
by $g\cdot (\{L_i\}, \alpha)\mapsto (\{g\cdot L_i\}, c(g)^{-1}\cdot \alpha)$.
By following the arguments of [PR], \S 4 (which deals with the case 
of the unitary group, i.e when $\alpha=1$), we 
see that  there is an $LG$-equivariant isomorphism
$ LG/L^+P_I\simeq \F_I $ of sheaves for the fpqc topology. Note, however, that if $n=2m$, $\F_I$ is not always a partial flag variety associated to $G$, since the stabilizer group $P_I$ is not always a parahoric subgroup. In fact, this happens
if and only if $m\notin I$, cf. \S\ref{parahoric}. In this case the parahoric subgroup $P_I^0$ has index $2$ in $P_I$ and the partial flag variety $LG/L^+P_I^0$ is isomorphic to the disjoint sum of two copies of $\F_I$.
 
Recall from \S \ref{oddeven} that we identified the local Dynkin diagram of $SU(V, \phi)$ with $\{0, \ldots, m\}$ when $n=2m+1$ is odd, and with $\{0, \ldots, m-2, m, m'\}$ when $n=2m$ is even. Using the general notation for affine flag varieties introduced in \S\ref{3a}, 
we have $LG/L^+P_I^0=\F^Y$. Here $Y=I$ when $n$ is odd or $n$ is even and $m-1\notin I$. If $n$ is even and $m-1\in I$, then $Y=(I\setminus\{m-1\})\cup \{m'\}$. 

\subsection{Embedding of the special fibers of local models in partial affine flag varieties.}\label{3c}

We now return to the set-up of \S\ref{naive}. In particular, we fix integers $r, s$ with $r+s=n$ and denote by $E$ the reflex field for $(F/F_0, r, s)$. Let $k$ denote an algebraic closure of the residue field of  $\O_E$. Similarly to [PR3] 
we construct a 
natural embedding of the geometric special fibers 
$\bar{M}_I^{\rm naive}= M_I^{\rm
naive}\otimes_{R}k$ into a partial affine flag
variety associated to the unitary similitude group $G$. For this we fix identifications
compatible with the actions of $\pi$ resp.\ $u$,
$$
\Lambda_j\otimes_{\mathcal O_{F_0}}k=
\lambda_j\otimes_{k[[t]]}k
$$ 
which sends the natural bases
of each side to one another. We therefore also obtain a
$k[[u]] / (u^2)$-module chain isomorphism
$$
\Lambda_\bullet\otimes_{\mathcal O_{F_0}}k\simeq
\lambda_\bullet\otimes_{k[[t]]}k\   ,
$$ 
which is in fact an
isomorphism of polarized periodic module chains [RZ]. Let $R$
be a $k$-algebra. For an $R$-valued point of $M^{\rm naive}$ we have
$$
\mathcal F_j\subset \Lambda_j\otimes_{\mathcal O_{F_0}}R=
(\lambda_j\otimes_{k[[t]]}k)\otimes_k R\   .
$$ Let
$ L_j:=L_{\F_j}\subset\lambda_j\otimes_{k[[t]]} R[[t]]$ be the
inverse image of $\mathcal F_j$ under the canonical
projection
$$
\lambda_j\otimes_{k[[t]]}R[[t]]\longrightarrow
\lambda_j\otimes_{k[[t]]}R\  .
$$ 
Also set $\alpha=-t^{-1}$. Then for this choice of $\alpha$,
$$
 L_0\subset  L_{1}\subset\cdots\subset L_m
\subset u^{-1} L_0=  L_{n}
$$
is a lattice chain in $R((u))^n$ which satisfies
conditions a) and b) in \S\ref{3b}. This gives
 a well-defined point in $\F(R)=(LG/L^+B)(R)$. We obtain in this way a morphism
$$
\iota: M^{\rm naive}\otimes_{\mathcal
O_E}k\longrightarrow \mathcal F\ \ ,
$$ 
which is a closed
immersion (of ind-schemes).

In the case of incomplete
lattice chains, one can proceed in a similar way and obtain
an embedding
$$
\iota_I: M_I^{\rm naive}\otimes_{\mathcal
O_E}k\hookrightarrow \mathcal F_I\   .
$$
For
$I'\supset I$, the following diagram is commutative, 
$$
\begin{matrix}
 M_{I'}^{{\rm naive}}\otimes_{\mathcal O_E}k & \hookrightarrow &
\mathcal F_{ I'}
\\
\big\downarrow && \big\downarrow
\\
  M_I^{{\rm naive}}\otimes_{\mathcal O_E}k & \hookrightarrow &
\ \ \mathcal F_{ I}\ .
\end{matrix} 
$$
The horizontal morphisms are equivariant for the actions of
$L^+P_{I'}$ resp.\ $L^+P_I$, in the sense of \cite{P-R2},
section 6.

Note that in [PR3], we have set $L_j:=u^{-1}L_{\F_j}$ to produce an embedding 
of $M_I^{{\rm naive}}\otimes_{\mathcal O_E}k$ into $LU_n/L^+P_I(U_n)$,
where $P_I(U_n)$ is the stabilizer of our lattice chain in the unitary group. Observe that the similitude of scalar multiplication by $u^{-1}$ is our choice of $\alpha=-t^{-1}$ here. Hence, we obtain a commutative diagram
$$
\xymatrix@C=27mm{
M_I^{{\rm naive}}\otimes_{\mathcal O_E}k\ar[r]^{\F_j\mapsto  u^{-1}\cdot L_{\F_j}}
\ar[dr]_{i_I} & \ \ \ \ LU_n/L^+P_I(U_n)\ar[d]^{u\cdot {\rm Id}}  \\
\ & \ \ \ \ \F_I=LG/L^+P_I\ . 
}$$

Let $\mu=\mu_{r, s}$ denote the cocharacter of $G$ given under the isomorphism $G\otimes_KK'\simeq GL_n\times \mathbb G_m$ by
\begin{equation}\label{mu}
\mu_{r, s}(z)=({\rm diag}(\ z^{(s)}, 1^{(r)}\ ), z) \ .
\end{equation}
Let $\A(\mu)=\bigcup_{w\in {\rm Adm}(\mu)} S_w$. This is a closed reduced subset of the full flag variety $\F=LG/L^+B$. Similarly, for $I\subset \{0, \ldots, m\}$, we denote by $\A^I(\mu)$
the image of $\A(\mu)$ in $\F_I$.
\begin{prop}\label{lift}
$\A^I(\mu)$ is contained in the image of $M_I^{{\rm loc}}\otimes_{\mathcal O_E}k$ under $\iota_I$.  
\end{prop}
It obviously suffices to prove this in the case of the full flag variety, i.e., for $I=\{0,\ldots,m\}$. This is done in the next subsection.

\subsection{Lifting of points in $\mu$-strata to the generic fiber}\label{liftmu}

Recall that we are treating the Iwahori case of signature $(r,s)$ with $s\leq r$. 
 Recall our notation of the standard lattices $\Lambda_j$;
these come with distinguished $\O_{F}$-generators. For example, for
$\Lambda_0$ these are $e_1,\ldots, e_n$. For $\Lambda_{m}$, if $n=2m$ is even, these are 
$f_1:=-\pi^{-1}e_1,\ldots f_m:=-\pi^{-1}e_m, f_{m+1}:=e_{m+1},\ldots , f_n:=e_n$.
We have
$$
\langle \pi e_i,  e_{n+1-j}\rangle=\delta_{i, j}
$$
and (for $n=2m$)
$$
(f_i, \pi f_{n+1-j})=\pm\delta_{i, j}
$$
where $\langle \ , \  \rangle$, resp. $(\ ,\ )$ is the alternating, resp. symmetric form associated to the hermitian form $\phi$, as in \S \ref{latt}. In particular, $\Lambda_0$, resp. $\Lambda_m$ is self-dual for $\langle \ , \ \rangle$, resp. $(\ , \ )$.
In general, if 
$$
\Lambda={\rm span}_{\O_F} \{\pi^{-1}e_1,\ldots ,\pi^{-1}e_i, e_{i+1}, \ldots, e_n\},
$$
 we will
set 
$$
f^{\Lambda}_1:= -\pi^{-1}e_1,\ldots , f^{\Lambda}_{i}:=-\pi^{-1}e_{i}, f^{\Lambda}_{i+1}:=e_{i+1},\ldots , f^\Lambda_n:=e_n
$$

If $T$ is a subset of   $\{1,\ldots, n\}=[1,n]$, we will 
set $T^*=\{n+1-a\ |\ a\in T\}$.
We denote by $S$ a subset of $[1,n]$ of cardinality $s$;
then $R=R_S$ will denote the complement of $S^*$. Note that $R$ has cardinality $r=n-s$
and that if $S=[1,s]$, then $R=[1,r]$.  We will consider 
subsets $S$ of cardinality $s$ for which $S\cap S^*=\emptyset$, i.e $S\subset R_S$. 

Set 
$$
\ti \pi=\pi\otimes 1-1\otimes \pi\in \O_{F}\otimes_{\O_{F_0}}\O_F
$$
We will write this as $\ti\pi=\pi-\sqrt\pi_0$.
For such a subset $S$ we will consider 
$$
\F^\Lambda_S=<f^\Lambda_{S^*},\  \pi f^\Lambda_{S^*},\  \ti \pi f^\Lambda_{R\setminus S}>\ \subset  \Lambda\otimes_{\O_{F_0}}\O_F
$$
where $f^\Lambda_T=\{f^\Lambda_t\}_{t\in T}$ and where $<, >$ denotes the $\O_{F_0}$-module generated by these vectors (not to be confused with the notation for the alternating form).
Notice that $R\setminus S=[1,n]\setminus (S\cup S^*)$.
We claim this corresponds to an $\O_{F}$-point of the naive local model for signature $(r,s)$
and the complete lattice chain.
Note the following identity, 
$$
\pi\cdot (\ti \pi v)=-\sqrt\pi_0\cdot (\ti \pi v)\ ,
$$
so that $\ti\pi v$ is an eigenvector for the action of $\pi$ with eigenvalue $-\sqrt{\pi_0}=-1\otimes\pi$.
This shows that $\F^\Lambda_S$ is stable for the action of $\pi=\pi\otimes 1$ and the characteristic polynomial 
is 
$$
\det(T\cdot I-\pi \ |\ \F^\Lambda_S)=(T^2- \pi_0)^{s} (T-\sqrt{\pi_0})^{r-s}=(T+\sqrt{\pi_0})^r(T-\sqrt{\pi_0})^s\ .
$$
Now we explain why $\F^\Lambda_S$ is an isotropic sequence of flags. Since this is a closed condition
it is enough to check it over $F$. Since $\Lambda\otimes_{\O_{F_0}}F=\Lambda_0\otimes_{\O_{F_0}}F$
it is enough to check that 
$$
<\F^{\Lambda_0}_S, \F^{\Lambda_0}_S>=0
$$
for all such $S$. Recall that $f^{\Lambda_0}_j=e_j$ for all $j$.
Since $S\cap S^*=\emptyset$, we have $\langle e_i, \pi e_j\rangle =0$ if $i, j\in S^*$.
We also have
$$
\langle\ti\pi v, \ti\pi w\rangle =\langle \pi v-\sqrt\pi_0 v, \,  \pi w-\sqrt \pi_0 w\rangle =0\,.
$$
Also $\langle e_i, e_j\rangle=0$ for all $i$, $j$ and $\langle e_i, \ti\pi e_j\rangle =0$ if $i\in S^*$, $j\in R\setminus S$.
Similarly $\langle \pi e_i, \ti\pi e_j\rangle =0$
for $i\in S^*$, $j\in R\setminus S$.

This shows that $\F^\Lambda_S$ gives an $\O_F$-point of $M^{\rm naive}$.
Denote by 
$$
\bar \F^{\Lambda}_S\subset \Lambda\otimes_{\O_{F_0}}(\O_{F}/(\pi)\otimes_{\O_{F}/(\pi)}k)=\Lambda\otimes_{\O_{F_0}}k
$$
the reduction of the subspaces $\F^{\Lambda}_S$ modulo the maximal ideal of $\O_F$. Denote by ${L_S}_{\bullet}$ the lattice chain that corresponds to the subspaces $\bar \F^{\Lambda}_S$ as in \S\ref{3b}. Then 
$$
u^2\lambda_{\bullet}  \subset {L_S}_{\bullet} \subset {\lambda_{\bullet}}
$$
(We also think of ${L_S}_{\bullet}$ as giving a point of the affine flag variety $\F=LG/L^+B$ for the unitary similitude group $G$.)

Recall from \S\ref{caseEv}, \ref{caseOdd}, that $W_0=S_m\rtimes \{ \pm 1\}^m$, no matter whether $n$ is even or odd. 
For each $w\in W_0$, and each subset $S\subset [1,n]$ of cardinality $s\leq m$ 
with $S\cap S^*=\emptyset$, $w\cdot S$ is of this type again
and this action is transitive.

\begin{lemma} { The relative position ${\rm inv}(\lambda_{\bullet}, {L_S}_{\bullet})$ in the Iwahori-Weyl group $\ti W$ of $G$ is
$w_S(\lambda)$. Here $\lambda$ is as in \S\ref{admissible} the image of $\mu\in X_*(T)$ in the coinvariants $X_*(T)_I$
and $w_S\in W_0$ has the property that $w_S\cdot [1,s]= S$.}
\end{lemma}

\begin{Proof}
Here we may think of the Iwahori-Weyl group $\ti W$ of 
the unitary similitude group $G\subset {\rm Res}_{F/F_0}(GL_n  \times \Gm) $  as a subgroup of the Iwahori-Weyl group for $GL_n \times \Gm$. This follows easily when we 
 we represent these groups as quotients of normalizer
groups; on the translation part this 
induces the injection given by the norm
\begin{equation}
X_*(T)_I\to X_*(T)=\Z^n\times \Z\ .
\end{equation}
(Note here that $X_*(T)_I$  is torsion-free.) 
Under this map the element $\lambda$ is sent to 
$$
{\rm Norm}(\mu_{r,s})=((2^s, 1^{r-s}, 0^{s}), 2)\, .
$$ 
Set
$$
\tau_{r,s}={\rm diag}(-t^{(s)}, u^{(r-s)}, 1^{(s)})
$$
for a corresponding torus element in the unitary similitude group.

Set $\Lambda=\Lambda_0$ if $n$ is odd and $\Lambda=\Lambda_{m}$ if $n=2m$ is even (this amounts to choosing a special vertex for $G$).
Similarly, set $\lambda$ for $\lambda_0$ or $\lambda_m$,  and $L_S$ for the lattices ${L_{S}}_{0}$ or ${L_{S}}_{m}$ respectively.
Then 
$$
{L_{[1,s]}}_{\bullet}=\tau_{r,s}\cdot \lambda_{\bullet}\ .
$$
This shows that the claim is true for $S=[1,s]$ and $w_{[1,s]}$ the identity; the same is true for any 
other choice of $w_{[1,s]}$ since each such choice centralizes $\tau_{r,s}$. 
Now write $\lambda_j=t_j\cdot \lambda$ where $t_j$ is a
(diagonal) translation element for $GL_n$.
Then we also have
$
(L_S)_j=t_j\cdot L_S
$
for any $j$. Now we think of elements of $W_0$ as permutation matrices in $GL_n$
in the standard basis $f^\Lambda_i$ of $\Lambda$. We can see that the permutation
$w_S$ takes the $k[[u]]$-basis 
$$
\{u^2f^{\Lambda}_1,\ldots , u^2f^{\Lambda}_s, uf^{\Lambda}_{s+1}, \ldots ,uf^{\Lambda}_{r},   f^{\Lambda}_{r+1},\ldots ,  f^{\Lambda}_n\}
$$
of $L_{[1,s]}$ to the basis 
$$
\{f^{\Lambda}_{S^*}\, ,\, uf^{\Lambda}_{R-S}\, ,\, u^2 f^{\Lambda}_{S}\}
$$
of $L_S$. This shows that $L_{S}=w_S\cdot L_{[1,s]}$ which combined with the above gives
\begin{equation}
(L_S)_j=t_j w_S L_{[1,s]}=t_jw_S \tau_{r,s} \lambda=t_jw_S \tau_{r,s} w_S^{-1} \lambda=t_jw_S \tau_{r,s} w_S^{-1} t_j^{-1} \lambda_j\ .
\end{equation}
This amounts to $(L_S)_\bullet=w_S\tau_{r,s}w_S^{-1} \Lambda_\bullet$ which proves the claim.\endproof
\end{Proof}
\medskip

The claim shows that the $\O_F$-points of $M^{\rm loc}$ that we constructed above (one for each subset $S$)
reduce to points lying in the Schubert varieties that correspond to the extreme elements
of the set ${\rm Adm}(\mu)$; all these extreme elements are then obtained by this construction. Since ${\rm Adm}(\mu)$ is closed under the Bruhat order, the set $\A(\mu)$ is closed in $\F$, and is contained in the closed subset $\iota (M^{\rm loc}\otimes_{\O_E}k)$. This proves Proposition \ref{lift}.
\medskip

\section{The structure of local models}
\setcounter{equation}{0}

\subsection{Consequences of the coherence conjecture.}

It is shown in [PR3] \S 11,  that as a consequence of  Proposition \ref{lift},  we have the following statement.
\begin{thm}
Assume the validity of the coherence conjecture for the pair $(G, \mu)= (GU(V,\phi), \mu_{r, s})$. Then the geometric special fiber of $M^{\rm loc}_I$ is reduced and is isomorphic to $\A^I(\mu)$. Its irreducible components are normal and with rational singularities.
\end{thm}

\subsection{Vertex-wise admissibility.}

Recall that we are assuming $s\leq r$.  By the previous section, all the extreme elements of the set ${\rm Adm}^I(\mu_{r,s})$ 
correspond to sequences of subspaces $\F_i$ for which the rank of $(\pi | \F_i)$
is maximal, i.e equal to $s$,  for all $i\in I$. 
\begin{conjecture}\label{vertexwise}
  The converse is true:
given a sequence of subspaces $\F_i$ (in the special fiber of the naive local model $M^{\rm naive}_I$)
for which the rank of $(\pi | \F_i)$ is equal to $s$ for all $i\in I$, then the corresponding stratum 
is admissible and extreme, i.e of the form $w_0(\lambda)$ for some $w_0\in W_0$ (or rather its image in 
$W^{I}\backslash \tilde W/W^{I}$).
\end{conjecture}

 In particular, we conjecture that admissibility 
has a vertex by vertex characterization.

\begin{Remark}
{\rm For ${\rm Res}_{F/F_0}GL_n$ (here again $F/F_0$ is ramified quadratic extension) this statement is true; this follows from [PR1].
In fact, in this case,  an element of the naive local model lies in a Schubert cell for an
element of ${\rm Adm}^I(\mu)$ if and only if the rank of $\pi$ on each subspace is $\leq s$
and the extreme elements are the ones for which the rank is $s$. This stronger statement is false in the case of a ramified unitary group. For example, let $n=2m$ be even and let $I=\{ m\}$. Then the subspace $\F=\pi \Lambda_m/ \pi_0 \Lambda_m$ of $\Lambda_m/ \pi_0 \Lambda_m$ has rank of $\pi | \F$ equal to zero, but does not belong to $\A^I(\mu)$ if $s$ is odd, cf. \ref{caseEv}. 

In the global context of \S\ref{unimoduli} U. G\"ortz and T. Wedhorn ask whether the points in the special fiber of $\A_{C'}$, for which the stronger statement fails, are all supersingular. 
}\end{Remark}

 \begin{prop}
Assume the validity of Conjecture \ref{vertexwise}. Then
$$
\A^I(\mu)=\bigcap\nolimits_{i\in I} \pi^{-1}_{I, \{i\}}(\A^{\{i\}}(\mu)) \ 
$$
(intersection inside $\F_I=LG/L^+P_I$). Assuming in addition the coherence conjecture (for
the pair $(G, \mu)= (GU(V, \phi), \mu_{r, s})$), 
$$
M^{\rm loc}_I=\bigcap\nolimits_{i\in I} \pi^{-1}_{I, \{i\}}(M^{\rm loc}_{\{i\}}) \ 
$$
(intersection inside $M^{\rm naive}_I$). 
\end{prop}
\begin{proof}
In either case, it is clear that the left hand side is contained in the right hand side. 
It suffices to prove  the second identity  for the special fibers, in which case it follows from the first identity since we are assuming the coherence conjecture. For the first identity, it suffices to see that any extreme stratum of the right hand side is contained in the left hand side. This follows from the following identity which is a consequence of conjecture \ref{vertexwise},
\begin{equation*}
{\rm Adm}^I(\mu)=\bigcap\nolimits_{i\in I}\pi_{I, \{i\}}^{-1}({\rm Adm}^{\{i\}}(\mu))
\end{equation*}
(intersection in $W^I\backslash \tilde W/W^I$).
\end{proof}
\medskip

\section{Special parahorics}\label{special}
\setcounter{equation}{0}

In the beginning of this section, we give affine charts for the  local models $M^{\rm naive}_I$ and $M^\wedge_I$
in the cases $I=\{0\}$ when $n$ is odd and $I=\{m\}$ when $n=2m$ is even. 
For simplicity, we will omit the subscript $\{0\}$ or $\{m\}$ and write simply $M^{\rm naive}$, $M^{\wedge}$ etc.

The corresponding parahoric subgroups $P^0_I$  are then special in the sense of Bruhat-Tits ([T]).
As we will see, the reduced special fibers of these charts are the closures of certain nilpotent orbits 
for the ``classical symmetric pairs" $({\mathfrak gl}(n), {\mathfrak o}(n))$, 
resp. $({\mathfrak gl}(n), {\mathfrak sp}(n))$.

\subsection{An affine chart} 

We distinguish two cases:
\smallskip

A. Assume  that $n=2m+1$ is odd; then we take $I=\{0\}$. Recall that $\Lambda:=\Lambda_0$ 
is equipped with the perfect alternating form $\langle\ ,\ \rangle$
which satisfies $\langle \pi e_i,   e_{n+1-j}\rangle=\delta_{ij}$. 
The corresponding $\O_{F_0}$-basis is
$$
\{e_1,\ldots, e_n, \pi e_1,\ldots ,\pi e_n\}
$$
and the matrix of $\langle\ ,\ \rangle$ in this basis is 
$$
J_{2n}=\begin{pmatrix} 0 & -H_n\\ H_n & 0\end{pmatrix} \ ,
$$
where $H_n$ is the unit antidiagonal matrix (of size $n$).
Let $\L=\langle e_1,\ldots , e_n\rangle$ be 
the standard Lagrangian
$\O_{F_0}$-direct summand of $\Lambda$, and denote by $Q$ the Siegel parabolic subgroup of $ Sp(\Lambda, \langle, \rangle)$ that preserves $\L$.  
\smallskip

B. Assume that $n=2m $ is even; then we take $I=\{m\}$. Recall that $\Lambda:=\Lambda_m$ 
is equipped with the perfect symmetric form $(\ ,\ )$. 
Set $f_i=-\pi^{-1}e_i$, for $1\leq i\leq m$,
$f_i=e_i$, for $m+1\leq i\leq n$. Consider the   
$\O_{F_0}$-basis of $\Lambda_m$ given by
$$
\{f_1,\ldots, f_n, -\pi f_1,\ldots ,-\pi f_m, \pi f_{m+1},\ldots, \pi f_n\}
$$
The matrix of $(\ ,\ )$ in this basis is $H_{2n}$.
Here we let $Q$ be the  parabolic subgroup of $ O(\Lambda_m, (, ))$
that preserves the totally isotropic $\O_{F_0}$-direct summand $\L =\langle f_1,\ldots, f_n\rangle$. 
\smallskip

In each case, the Lie algebra ${\rm Lie}(Q)$ consists of the $\O_{F_0}$-endomorphisms
$A$ of $\Lambda$ that satisfy $A(\L)\subset \L$
and \begin{equation}
\langle Av, w\rangle+\langle v, Aw\rangle=0, \quad \hbox{\rm resp. \ } (Av, w)+(v, Aw)=0 .
\end{equation}
For simplicity, denote by $X$ the restriction of $A$ to $\L$.
Now let $N_{r,s}$ be the scheme of endomorphisms
$A$ in the Lie algebra ${\rm Lie}(Q)\otimes_{\O_{F_0}}\O_E$ that satisfy 
the equations 
$$
A^2=\pi_0\cdot I_{2n}\ , \hbox{\rm\ \, and}
$$
$$
{\rm char}_{X}(T)=(T-\sqrt \pi_0)^s(T+\sqrt \pi_0)^r\ .
$$

Also, denote by $N^\wedge_{r,s}$ the closed subscheme of $N_{r,s}$
given by $A$ which when $r\neq s$ satisfy in addition:
\begin{equation}\label{wedgeX}
\wedge^{r+1}(X-\sqrt\pi_0\cdot I)=0,\quad \wedge^{s+1}(X+\sqrt\pi_0\cdot I)=0\ .
\end{equation}

As in [PR1] (1.3), we see that there is a diagram
\begin{equation}
M^{\rm naive}\xleftarrow{\ \pi\ } \widetilde M^{\rm naive}\xrightarrow{\ \phi\ } N_{r,s}\ ,
\end{equation}
where $\pi$ is a $Q$-torsor, and where $\phi$ is a smooth morphism. 
Here 
$$
\widetilde M^{\rm naive}(S)=\{(\F\subset \Lambda\otimes_{\O_{F_0}}\O_S, \alpha)\}
$$
 where $\F$ gives a point of $M^{\rm naive}(S)$ and $\alpha$ is a 
choice of basis $v_1,\ldots, v_{2n}$  of $\Lambda\otimes_{\O_{F_0}}\O_S$ with matrix $J_{2n}$, resp. $H_{2n}$ for the bilinear form on $\Lambda\otimes_{\O_{F_0}}\O_S$ and  such that  the subspace $\F$ is generated by $v_1,\ldots,v_n$.  Here
$\phi((\F, \alpha))$
is given by the endomorphism $\alpha^{-1}\cdot \pi\cdot \alpha$
(which we can express by a matrix $A$ using $\{v_i\}$). Similarly, by restriction along the closed immersions $ M^{\wedge}\subset M^{\rm naive}$, $N^\wedge_{r,s}\subset N_{r,s}$
we also obtain 
\begin{equation}
M^{\wedge}\xleftarrow{\ \pi^{\wedge}\ } \widetilde M^{\wedge} \xrightarrow{\ \phi^{\wedge}\ } N^\wedge_{r,s}\ ,
\end{equation}
with $\pi^\wedge$ a $Q$-torsor and $\phi^\wedge$ smooth again. 
Notice that the morphisms $\phi$, $\phi^\wedge$ are not surjective; 
for example, the closed fiber of their image lies in the open subset where the rank of $A$ 
is exactly $n$.

\subsection{The odd case} \label{odd}

Assume we are in case A. 
For each signature $(r,s)$, we can consider the point of $M^{\wedge}(k)\subset M^{\rm naive}(k)$ given by 
$$
\F_0=\pi \cdot (\Lambda_0\otimes_{\O_{F_0}}k) \subset \Lambda_0\otimes_{\O_{F_0}}k\ .
$$
 To obtain 
an affine chart of the naive local model around this point consider the subspaces $\F$  which are of the form
$$
\F=\{ \pi\cdot v+X\cdot v\mid v\in \O_S^n\}\ ,
$$
with $v$ a column vector in $\O_S^n=\O_S\cdot e_1\oplus\cdots \oplus \O_S\cdot  e_n$. 
Then we can see as in [P1], p. 596, that the corresponding open subscheme $U_{r,s}$ of the local model 
$M^{\rm naive}$ is the scheme of $n\times n$-matrices $X$ which satisfy 
\begin{equation}\label{Xequ1}
X^2=\pi_0\cdot I, \quad X^t=HXH, \quad
 {\rm char}_X(T)=(T-\sqrt \pi_0)^s(T+\sqrt \pi_0)^r\ .
\end{equation}
Here $H=H_n$.
Similarly, the corresponding open subscheme
$U^\wedge_{r,s}$ of the wedge local model 
$M^{\wedge}$ is the scheme of $n\times n$-matrices $X$ which satisfy the equations (\ref{Xequ1})
above together with the equations of minors (\ref{wedgeX}).

To connect this with the previous picture, notice that over $U_{r,s}$ we 
can split the $Q$-torsor $\pi$ by picking the symplectic basis $\{v_i\}$ of $\Lambda_0\otimes_{\O_{F_0}}\O_S$ to be
$\pi e_1+X e_1, \ldots , \pi e_n+X e_n, -e_1, \ldots ,-e_n$. Then $U_{r,s}$ can be  identified 
with a subscheme of the scheme of matrices $A$ of $N_{r,s}$ of the form
\begin{equation*}
A=\left(\begin{matrix} X& -I\\ 0 &-X \end{matrix}\right)\ .
\end{equation*}
Similarly for $U^\wedge_{r,s}$.

\subsection{The even case}\label{even}

 Assume now that we are in case B.  
We can see as above that an open affine chart $U_{r,s}$ of $M^{\rm naive}$ around the point $\F_0=\pi(\Lambda_m\otimes_{\O_{F_0}}k)\subset \Lambda_m\otimes_{\O_{F_0}}k$
is given by the scheme of $n\times n$-matrices $X$ which satisfy 
\begin{equation}\label{Xequ2}
X^2=\pi_0\cdot I, \quad X^t=-JXJ, \quad
 {\rm char}_X(T)=(T-\sqrt \pi_0)^r(T+\sqrt \pi_0)^s\ .
\end{equation}
Here $J=J_{n}$. 
Similarly, the corresponding open subscheme
$U^\wedge_{r,s}$ of  
$M^{\wedge}$ is the scheme of $X$ which satisfy the equations (\ref{Xequ2})
above together with the equations of minors (\ref{wedgeX}).

Suppose $s\neq 0$.
It is then also useful to  consider  affine charts $_1U_{r, s}$, $_1U^{\wedge}_{r,s}$
of $M^{\rm naive}$, $M^{\wedge}_{r,s}$ around the point 
$$
\F_1=\langle f_1, \pi f_1,\pi f_2, \ldots , \pi f_{n-1} \rangle\subset \Lambda_m\otimes_{\O_{F_0}}k
$$
of $M^\wedge(k)\subset M^{\rm naive}(k)$ (especially when $s$ is odd). Set 
$$
\Lambda'={\rm span}_{\O_{F_0}}\langle f_1,\pi f_1, \pi f_2,\ldots , \pi f_{n-1} \rangle,
$$
$$
\Lambda''={\rm span}_{\O_{F_0}}\langle  f_n, \pi f_n, f_{2}, \ldots , f_{n-1}  \rangle 
$$
so that $\Lambda_m=\Lambda'\oplus\Lambda''$. 
Notice that the matrix of the form $(\, ,\, )$ in the basis
\begin{equation}\label{ba}
\{f_1,\pi f_1, \pi f_2,\ldots , \pi f_{n-1} ,  f_n, \pi f_n, f_{2}, \ldots , f_{n-1}\} 
\end{equation}
is 
$$
\begin{pmatrix}0&S\\ S^t&0\\
\end{pmatrix}
$$
where $S$ is the skew matrix of size $n$,
$$
S=\begin{pmatrix}J^t_2&0\\
0&J_{n-2}\\
\end{pmatrix}.
$$
\smallskip

Now we can find an affine chart $_1U_{r,s}$ around the point $\F_1$ by looking at graphs of linear maps
$f: \Lambda'\to \Lambda''$ so that 
$$
\F=\{v+f(v)\mid v\in \Lambda'\otimes\O_S\} \ .
$$
Write the map $f$ by an $n\times n$ matrix 
$$
X= \begin{pmatrix}T&B\\
C&Y\\
\end{pmatrix}
$$
by using the basis (\ref{ba}) above. Here $T$ is a $2\times 2$ matrix
and $Y$ an $(n-2)\times(n-2)$ matrix; the matrices $B$ and $C$ 
are of sizes $2\times (n-2)$ and $(n-2)\times 2$ respectively.
We can see that the condition that $\F$ is isotropic
is given by
\begin{equation}\label{24b}
SX^t= XS\ .
\end{equation}
This translates to $J_2\cdot T^t=T\cdot J_2$, $J_{n-2}\cdot Y^t=Y\cdot J_{n-2}$
and $C\cdot J^t_{2}=J_{n-2}\cdot B^t$. The 
first condition implies that $T$ is diagonal and scalar,  $T={\rm diag}(x,x)$;
the last condition shows that $C$ is determined by $B$.
The condition that $\F$ is
$\pi$-stable translates to
$$
Y^2=\pi_0\cdot I_{n-2}, \quad B_1=B_2\cdot Y
$$
where $B_i$ is the $i$-th row of $B$; therefore $B_1$ is determined by $B_2$.
We can see that the action of $\pi$ on $\F$ is now given by
the block sum of the matrices $\begin{pmatrix}0&\pi_0\\ 1&0\end{pmatrix}$
and $Y$. This allows us to translate the characteristic polynomial and exterior power conditions
to conditions about the matrix $Y$. We obtain
\begin{equation}
{\rm char}_Y(T)=(T -\sqrt{\pi_0})^{s-1}(T+\sqrt{\pi_0})^{r-1},
\end{equation}
and for $M^{\wedge}_{r,s}$
\begin{equation}
\wedge^{s}(Y+\sqrt{\pi_0}\cdot I_{n-2})=0, \
\wedge^{r}(Y-\sqrt{\pi_0} \cdot I_{n-2})=0.
\end{equation}
We conclude that in this case, the affine chart $_1U_{r,s}$ is the product of the affine scheme
$U_{r-1, s-1}$ with an affine space of dimension $n-1$ over $\O_E$ 
(which corresponds to the free coordinates $x$, $B_2=(b_{2\, 1},\ldots b_{2\, n-2})$).
In particular, we can see that if $s$ or $r$ is $1$, then $_1U_{r,s}$ is smooth.

\subsection{Symmetric pairs} Let $\fG$ be a reductive Lie algebra over a field $k$ of odd characteristic. 
Let $\theta$ be a non-trivial involution of $\fG$. Let $\fG=\fF+\fP$ be the Cartan decomposition so that
$\fF=\{X\in \fG: \theta(X)=X\}$, $\fP=\{X\in \fG: \theta(X)=-X\}$. Then $(\fG, \fF)$ form a ``symmetric 
pair"; we will call $\fP$ the associated vector space (we can think of it as an infinitesimal 
version of a symmetric space). Now let $G$ be the adjoint group of $\fG$ and let $H=G^\theta$ be the subgroup
 of elements of $G$ fixed by the involution. The group $H$ acts on $\fP$ (via the adjoint action).
The orbits of $H$ on $\fP$, and in particular the orbits of  elements of $\fP$ which are nilpotent in $\fG$,
have been studied by several people starting with Kostant and Rallis (see [KR1], [KR2], [Oh], etc.).
When $k=\R$, $\fG$ is semi-simple and $\theta$ is a  Cartan involution, then $H(\bf R)$
is compact and one has the so-called Kostant-Sekiguchi correspondence [S2]: this is a bijective
correspondence between the nilpotent $H_\C$ orbits of $\fP_\C$ and the nilpotent orbits
of $H(\R)$ on the real Lie algebra $\fG(\R)$. (This plays no role in what follows).
Here we consider two ``classical cases":

Let $V$ be a  vector space over $k$ of dimension $n$ with a non-degenerate bilinear form 
$h: V\times V\to k$ which is either symmetric or alternating. If $X\in {\rm End}(V)$ the adjoint $X^*$ is defined as usual
by $h(Xv, w)=h(v, X^*w)$. Then $\theta(X)=-X^*$ gives an involution of the 
Lie algebra ${\rm End}(V)\simeq {\mathfrak gl}(n)$.  
\smallskip
 
A. The form  is symmetric given by   the antidiagonal matrix  $H=H_n$. Then $(\fG, \fF)=({\mathfrak gl}(n), {\mathfrak o}(n))$; 
the vector space $\fP$ is $\{X\in {\rm Mat}_{n\times n}(k)\ |\ X^t=HXH\}$ and the orbits
are for the action of the (split) orthogonal group. The nilpotent orbits in $\fP$ are parametrized by partitions $P(n)$ of $n$: For $\lambda\in P(n)$,
the corresponding orbit is $\fP_{\lambda}=\fP\cap \O_\lambda$ with $\O_\lambda$
the $GL_n$-orbit of a nilpotent matrix with Jordan blocks given by $\lambda$. 

\smallskip

B. The form is alternating given by the  skew-symmetric matrix $J=J_{2m}$ 
where $n=2m$.  Then $(\fG, \fF)=({\mathfrak gl}(n), {\mathfrak {sp}}(n))$; the vector space
$\fP$ is $\{X\in {\rm Mat}_{n\times n}(k)\ |\ X^t=-J   X  J\}$. The orbits
are for the action of the symplectic group and are parametrized by partitions of $m$;
if $\lambda=(a_1,a_2,\ldots, a_s)$ is a partition of $m$ then 
$\fP_{\lambda}=\fP\cap \O_{\lambda^{(2)}}$ with $\lambda^{(2)}=(a_1,a_1,a_2, a_2,\ldots ,a_s, a_s)$
a partition of $n=2m$.
\smallskip

(Note that the parametrization of the nilpotent orbits $\fP_{\lambda}$ stated above
is shown  by Sekiguchi [S1], see also [Oh], when $k=\C$. This result extends 
to any $k$ of odd characteristic.)

\subsection{Generic smoothness} Here we show that our results so far allow us to deduce that the special fiber
$\bar M^{\rm loc}$ is irreducible and generically smooth in these special
cases (when $n=2m+1$ and $I=\{0\}$, or $n=2m$ and $I=\{m\}$). 

By \S \ref{3c}, $\bar M^{\rm loc}$ is connected and projective and is a union of left orbits
for $P_I$ in the affine Grassmannian for $GU_n$. We can  see that there is
a unique closed $P_I$-orbit in $\bar M^{\rm loc}$.  This closed orbit has to be contained in $\A^I(\mu)$. From  the description
of ${\rm Adm}_0(\mu_{r, s})$ in \S\ref{caseEv} and \S\ref{caseOdd} we see that the closed orbit is the point $\{\F_0\}$
if $n=2m+1$, or if $n=2m$ and $s$ is even. If $n=2m$ and $s$ is odd, then the closed orbit is the orbit of  
$\F_1={\rm span}_{\O_F}\{f_1,\pi f_1, \pi f_2,\ldots , \pi f_{n-1}\}$. By the above, we can now obtain 
a description of an affine chart of $M^{\rm loc}$ that contains a point from the closed orbit:
In each case, such a chart is given by the flat closure $U_{r,s}^{\rm flat}$ of $U_{r,s}^{\wedge}$,
resp. $_1U_{r,s}^{\rm flat}$ of $_1U_{r,s}^{\wedge}$.
(This is the same as the flat closure of $U_{r,s}$ resp. $_1U_{r,s}$ .) 

Now consider the flat closure $V^{\rm flat}_{r,s}$
of the scheme  $V_{r,s}$ of $n\times n$ matrices $X$ over $\O_E$ that satisfy
$$
X^2=\pi_0\cdot I, \quad {\rm char}_{X}(T)=(T-\sqrt \pi_0)^r(T+\sqrt \pi_0)^s\
$$
$$
\wedge^{s+1}(X-\sqrt\pi_0\cdot I)=0,  \quad \ \wedge^{r+1}(X+\sqrt\pi_0\cdot I)=0\ , \ \hbox{\rm when} \ r\neq s.
$$
By [PR1], we see that the scheme $V^{\rm flat}_{r,s}$  has relative dimension $2rs$.
Its special fiber is reduced, irreducible and 
is the union of the nilpotent orbits $\O_\rho$ for $GL_n$ that 
correspond to partitions $\rho$ of $n$ with $\rho\leq (2^s, 1^r)=(r,s)^\vee$ (the partition dual to $(r,s)$). 
Its smooth locus $V^{\rm sm}_{r,s}$ is the complement of $\O_{(2^{s-1}, 1^{r+1})}$. 
The involution $\sigma$ given by
$X\mapsto H\, X^t\, H$ when $n$ is odd, resp. $X\mapsto -J\, X^t\,J$ when $n$ is even, acts 
on $V_{r,s}$, $V^{\rm flat}_{r,s}$ and on $\O_{(2^s, 1^r)}$. 

Suppose first that either $n$ is odd, or that both $n$  and $s$ are even. 
Then, by definition, $U^{\wedge}_{r,s}$ is the fixed point scheme
$(V_{r,s})^{\sigma}$; over the generic fiber $E$ this scheme has dimension $rs$. 
In this case, by Prop. 1 of [Oh], we can see that there is at least one 
$\sigma$-fixed point  on $\O_{(2^s, 1^r)}$;  since the action is tame,  the fixed point scheme 
$(\O_{(2^s,1^r)})^\sigma$ is smooth. In fact, it follows from [KR2], Prop. 5 and its proof
(they deal with the case $k=\C$ but their proof is valid for every algebraically closed field of odd characteristic) that each component of $(\O_{(2^s,1^r)})^\sigma$ has dimension half of the dimension of the corresponding nilpotent orbit $\O_{(2^s, 1^r)}$.
Hence, this dimension is ${\rm dim}((\O_{(2^s,1^r)})^\sigma)=rs$. 
We can now deduce that $(V_{r,s}^{\rm sm})^\sigma\subset U_{r,s}^{\wedge}$
is smooth over $\O_E$ of relative dimension $rs$ and provides an open 
subscheme of $U^{\rm flat}_{r,s}$ with (smooth) special fiber $(\O_{(2^s,1^r)})^\sigma$. 
In fact,  by Prop. 1 of [Oh],  $(\O_{(2^s,1^r)})^\sigma$ is connected, hence  we can see 
that $(\O_{(2^s,1^r)})^\sigma$ coincides with a nilpotent orbit for the associated symmetric pair, equal to $\fP_{(2^s, 1^r)}$ in the notation of the previous paragraph. It follows that the inclusion $\bar U_{r, s}^{\rm flat}\subset \bar U_{r, s}^\wedge=
\overline{ \fP_{(2^s, 1^r)}}$ is an equality on points. 
 It now follows that the special fiber of
$U^{\rm flat}_{r,s}$ is irreducible and generically smooth.  We also deduce that  $\bar M^{\rm loc}$ is irreducible. 
 
Suppose now that $n$ is even and $s$ is odd. (Note that in this case, $(\O_{(2^r,1^s)})^\sigma=\emptyset$.)
Then we can reduce to the 
smaller case where $n=r+s$ is replaced by $n-2$ partitioned by $(r-1, s-1)$; indeed, by our work in the last part of \S \ref{even}
and the above, we can see that the special fiber of the flat closure 
$_1U_{r,s}^{\rm flat}$ is generically the product of an affine space
with $(\O_{(2^{s-1}, 1^{r-1})})^\sigma$ and hence it is 
generically smooth and irreducible.

 As a corollary of the main result of [PR3] we can now show:
\begin{thm}\label{mainspecial}
Let $I=\{0\}$ if $n$ is odd, and $I=\{m\}$ if $n=2m$ is even. The special fiber of the local model $M_I^{\rm loc}$ is irreducible and reduced and is normal, Frobenius split and with only rational singularities. 
\end{thm}
\begin{proof}
By the previous considerations  the underlying reduced scheme $(\bar M_I^{\rm loc})_{\rm red}$ is a Schubert variety equal to $\A^I(\mu)$. By the main result of [PR3] we deduce that $(\bar M_I^{\rm loc})_{\rm red}$ is normal, Frobenius split and has only rational singularities. On the other hand, we saw that the special fiber $\bar M_I^{\rm loc}$ contains an open dense  subset which is reduced. By Hironaka's Lemma (EGA IV.5.12.8) we deduce that
$\bar M_I^{\rm loc}$ is reduced, which proves the claim, comp. also [PR3], Remark 11.4.
\end{proof}

\begin{conjecture}\label{conj52}
Let $n$ be odd, or both $n$ and $s$ be even. Then the scheme of matrices $U_{r,s}^\wedge$ is flat over $\Spec \O_E$. Equivalently, consider the space of matrices
$X$ in $M_n$ over $\Spec k$ with 
$$
X^2=0, \quad X^t=HXH, \quad
 {\rm char}_X(T)=T^n\ ,\quad \wedge^{s+1}X=0,\quad  \wedge^{r+1}X=0\ ,
$$
if $n$ is odd, resp. 
$$
X^2=0, \quad X^t=-JXJ, \quad
 {\rm char}_X(T)=T^n\ ,\text{ and } \wedge^{s+1}X=0,\quad  \wedge^{r+1}X=0\ , \text{ when } r\neq s \ ,
$$
if $n$ and $s$ are even. 

Then this scheme is reduced (in which case  it is normal, with rational singularities). 

\end{conjecture}

\begin{Remarknumb}\label{nonflat}
{\rm  a) Conjecture \ref{conj52} implies that, under the assumptions 
$n$ is odd, or both $n=2m$ and $s$ even, $M^{\wedge}_{\{0\}}=M^{\rm loc}_{\{0\}}$, resp. 
$M^{\wedge}_{\{m\}}=M^{\rm loc}_{\{m\}}$.

b) Suppose  that $n$ is even and $s$ is odd. Recall that we denote by $U_{r,s}^{\wedge}$ 
the affine chart around the point $\F_0$ of  $M^{\wedge}(k)=M^{\wedge}_{\{m\}}(k)$. By the above,
we have $U_{r,s}^{\wedge}=(V_{r,s}^\wedge)^\sigma$. However, we can see that in this case
the generic fiber $U_{r,s}^{\wedge}\otimes_{\O_{E}}E$ is empty. 
(Indeed, when $\pi_0$ is invertible, if  $X^tJ=JX$ and $X^2=\pi_0$, then the matrix $(\sqrt{\pi_0})^{-1}X$ 
belongs to the symplectic group and hence has determinant $1$. 
Therefore, both eigenvalues $\sqrt{\pi_0}$, $-\sqrt{\pi_0}$ of $X$
have to appear with even multiplicity.) In fact, in this case $(\O_{(r,s)})^\sigma$ is empty and an argument as above shows that ${\rm dim}(U_{r,s}^{\wedge})<rs$.
It follows that the point $\F_0$ does not lift to characteristic zero.
Hence $M^{\wedge}_{\{m\}}$ and also $M^{\rm naive}_{\{m\}}$ are not flat over $\O_E$.}
\end{Remarknumb}

\subsection{Normality of some nilpotent orbits}
In the proof of the previous theorem we used some facts about nilpotent orbits for symmetric pairs. Conversely, we can use the results on local models (which ultimately rely on structure theorems for affine flag varieties) to deduce results on nilpotent orbits. 
\begin{thm}
 a) Suppose that $n$ is odd. Then the Zariski closure of the nilpotent orbit $\fP_{(2^s, 1^r)}$   for the symmetric pair $({\mathfrak gl}(n), {\mathfrak o}(n))$ is normal.

b) Suppose that $n$ and $s$ are both even. Then the Zariski closure of the nilpotent orbit $\fP_{(2^s, 1^r)}$ for the symmetric pair $({\mathfrak gl}(n), {\mathfrak {sp}}(n))$ is normal. 
\end{thm}

(When $k=\C$ this result is a very special case of the results of Ohta [Oh];
his methods are particular to characteristic $0$. Note that Ohta shows that the Zariski closure of any nilpotent orbit for the pair  $({\mathfrak gl}(n), {\mathfrak {sp}}(n))$ is normal, but exhibits examples for  pairs $({\mathfrak gl}(n), {\mathfrak o}(n))$ of non-normal orbit closures.)
\smallskip

\begin{Proof} 
By Theorem \ref{mainspecial} the special fiber $\bar M^{\rm loc}$ is reduced and normal.
By our discussion above, the Zariski closure  of the nilpotent orbit $\fP_{(2^s, 1^r)}$
can be identified with an open affine subscheme of $\bar M^{\rm loc}$. The result follows.\endproof
\end{Proof}
\medskip

\section{The local models of Picard surfaces}
\setcounter{equation}{0}

In this case ($G=GU(2,1)$), there are three conjugacy classes of parahoric subgroups. In what follows we will show Theorem \ref{picard} of the introduction.
\smallskip

A) Let $I=\{0\}$. Then, by [P1], Theorem 4.5 and Remark 4.15., 
we have the following statement. 
\begin{prop}\label{M0}
$M^{\rm naive}_{\{0\}}$ is normal and Cohen-Macaulay. Furthermore,  $M^{\rm naive}_{\{0\}}$ is flat over $ \Spec(\O_F)$  and is smooth outside the special point $\F_0$ of the special fiber. Blowing up this special point yields a semi-stable model with special fiber consisting of two smooth surfaces meeting transversely along a smooth curve.
\end{prop}
We note that in [Kr] it is shown that the blow-up scheme represents a moduli problem analogous to the Demazure resolution of a Schubert variety in the Grassmannian.

B) Let $I=\{1\}$. In this case we have the following statement.
\begin{prop}\label{M1}
$M^{\rm loc}_{\{1\}}$ is smooth over $\Spec(\O_F)$. 
\end{prop}
\begin{proof}
The dual $\hat\Lambda^s_1$ of $\Lambda_1$ with respect to the symmetric bilinear form $(\ ,\ )$ is $\Lambda_2$. More precisely, the matrix of this bilinear form with respect to the  $\O_{F_0}$-basis 
$e_3, \pi^{-1}e_1, e_2, \pi e_3, e_1, \pi e_2$ (in this order)
is equal to 
$$
D= \begin{pmatrix}K&H\\
^t\!H&-\pi_0K\\ 
\end{pmatrix}\ ,
$$
where 
$$
K= \begin{pmatrix}0&0&0\\
0&0&0\\
0&0&1\\
\end{pmatrix}\ , \ \ 
H= \begin{pmatrix}0&1&0\\
-1&0&0\\
0&0&0\\
\end{pmatrix}\ .
$$
We describe an open neighborhood of the point $\F_0$ in the special fiber. 
As in \ref{odd}, we find that it can be given as a subscheme
of the space of matrices $A$ of the form 
$$A= \begin{pmatrix}X\\
I\\
\end{pmatrix}\ ,
$$
where $I$ is the unit matrix of size $3$, and where $X$ is a square matrix of size $3$ with indeterminates as entries. The  following conditions are imposed (the first condition corresponds to the isotropy condition imposed on $\F={\rm Image}(A)$, the second one to the $\pi$-invariance of $\F$):

(i) $^tA\cdot D\cdot A=0  $,

(ii) $X^2=\pi_0\cdot I $.

\noindent 
Now  condition (i) comes to 
\begin{equation}\label{dual}
^t\!XKX+(^t\!H X+^t\!X H)-\pi_0\cdot K=0\ .
\end{equation}
Write $X$ in the form 
$$
X= \begin{pmatrix}X_1&X_2\\
X_3&X_4\\
\end{pmatrix}
$$
where $X_1$ is a square matrix of size $2$ and $X_4$ a square matrix of size $1$.  Also 
let 
$$
J=\begin{pmatrix}
0&1\\
-1&0
\end{pmatrix}\ 
$$
be the left upper corner of $H$. 
Let us write
$$
X_1=
\begin{pmatrix}
a&b\\
c&d\\
\end{pmatrix} \ ,\qquad X_3=(x,y)\ .
$$
 Then equation (\ref{dual}) becomes the matrix equation
\begin{equation}\label{duala}
\begin{pmatrix}
^t\!X_3 X_3&^t\!X_3 X_4\\
 X_4 X_3& X_4^2\\
\end{pmatrix}\ 
+\ 
\begin{pmatrix}
^t\!J X_1+^t\!X_1 J&^t\!JX_2\\
^t\!X_2 J&-\pi_0\\
\end{pmatrix}\ 
=\ 0 \ ,
\end{equation}
where $^t\!X_3X_3$ means $\begin{pmatrix}
x^2&xy\\
xy&y^2
\end{pmatrix}$.

The condition (ii) becomes
\begin{equation}\label{square}
X^2=
\begin{pmatrix}
X_1^2+X_2 X_3&X_1 X_2+X_2 X_4\\
X_3 X_1+ X_4 X_3&X_3 X_2+X_4^2\\
\end{pmatrix}
=\pi_0\cdot I \ .
\end{equation}
From the left lower corner we get $^t\!X_2=X_4\cdot (y, -x)$.   
From the left upper corner of (\ref{duala}) we obtain the identity
\begin{equation}\label{LU}
-\begin{pmatrix}
-2c&a-d\\
a-d&2b
\end{pmatrix}
=
\begin{pmatrix}
x^2&xy\\
xy&y^2
\end{pmatrix} \ .
\end{equation}
On the other hand,  for the characteristic polynomial we are imposing  the condition
$$
{\rm char}_X(T)= (T+\sqrt{\pi_0})^2 (T-\sqrt{\pi_0}) \ .
$$
In particular ${\rm tr}(X)=-\sqrt{\pi_0}$, i.e. $a+d+X_4=-\sqrt{\pi_0}$. 
Taken together with (\ref{LU}),
we see that $X$ is uniquely determined by $(x, y)$ and $X_4$ with $X^2_4=\pi_0$. In the generic fiber of our scheme of matrices $X$
we see that we have $X_4=-\sqrt{\pi_0}$. (Use (\ref{square}), (\ref{LU})
and the characteristic polynomial condition.) Therefore, the 
identity $X_4=-\sqrt{\pi_0}$ persists in the flat closure:
  We conclude that the intersection of the flat closure 
$M^{\rm loc}_{\{1\}}$ with an  open
neighborhood of $\F_0$ can be identified with a closed subscheme of the affine space
in $(x, y)$ over $\Spec(\O_F)$. It therefore has to coincide with this affine space. Hence
 $M^{\rm loc}_{\{1\}}$ is smooth in a neighborhood of $\F_0$. But then the singular locus is empty since otherwise, as a closed subset in the special fiber invariant under the parahoric subgroup, it would have to contain $\F_0$.  
  \end{proof}
 \begin{Remarks}{\rm
 (i) One can see that $M^{\rm naive}_{\{1\}}$ and even $M^{\wedge}_{\{1\}}$ are not flat; computer calculations indicate that  they  have non-reduced special fibers, cf. Remark \ref{computer}, (iv).
 
   (ii) The special fiber $\overline M$ of $M^{\rm loc}_{\{1\}}$ is isomorphic to ${\bf P}^2$. To see this, let $\bar V=V/R$, where $V=\Lambda_1\otimes_{\O_{F_0}}k$, and where $R=(\pi e_2)$ is the one-dimensional radical of the symmetric bilinear form on $V$. Then $\bar V$ is a non-degenerate quadratic space of dimension $5$, and the image of $\pi \Lambda$ in $\bar V$ is a two-dimensional isotropic subspace $\bar V_2$. Now any $\F$ corresponding to a point in  $\overline M$ contains $R$, and hence defines a two-dimensional isotropic subspace $\bar \F$ in 
   $\bar V$. The fact that $\F$ is $\pi$-stable implies that the intersection $\bar \F\cap\bar V_2$ is non-trivial. We see that $\overline M$ is contained in the Schubert variety of isotropic planes in $\bar V$ which have a non-trivial intersection with the fixed isotropic plane $\bar V_2$. However, the Grassmannian of isotropic planes in $\bar V$ is isomorphic to the Grassmannian of (isotropic) lines in a $4$-dimensional symplectic space $(W, \langle\ ,\  \rangle)$. Indeed, $\bar V$ may be identified with a natural subspace of $\wedge^2 W$, and a line $l$ in $W$ is mapped to $l\wedge l^\perp$. If $L$ is mapped to $\bar V_2$ under this map, then the Schubert variety in question is identified with the set of lines in $W$ which are contained in the $3$-dimensional space $L^\perp$. Hence this Schubert variety is isomorphic to ${\bf P}^2$. It follows that $\overline M_{\rm red}$ is contained in ${\bf P}^2$, and therefore by Proposition  \ref{M1}, $M^{\rm loc}_{\{1\}}\otimes_{\O_F}k$ is isomorphic to ${\bf P}^2$.
   }
 \end{Remarks}
 
 C) Let $I=\{0, 1\}$. In this case we have the following statement.
 \begin{prop}
 The local model $M^{\rm loc}_{\{0, 1\}}$ is normal and Cohen-Macaulay and has reduced special fiber. It has two irreducible components which are normal and with only rational singularities. These two irreducible 
 components meet along two smooth curves which intersect transversally in a single point.
 \end{prop}

 \begin{proof}
 For the special fiber $\bar M^{\rm loc}_{\{0, 1\}}$ there is a chain of closed immersions, where $\mu= \mu_{(2, 1)}$,   
 \begin{equation}\label{inters}
 \A^{\{0, 1\}}(\mu)\subset \bar M^{\rm loc}_{\{0, 1\}}\subset \pi_{\{0\}}^{-1}(\bar M^{\rm loc}_{\{0\}})\cap  \pi_{\{1\}}^{-1}(\bar M^{\rm loc}_{\{1\}})= \pi_{\{0\}}^{-1}(\A^{\{0\}}(\mu))\cap \pi_{\{1\}}^{-1}(\A^{\{1\}}(\mu))\ .
 \end{equation}
 For the $\mu$-admissible set we have
 $${\rm Adm}^{\{0, 1\}}(\mu)=\pi_{\{0\}}^{-1}({\rm Adm}^{\{0\}}(\mu))\cap \pi_{\{1\}}^{-1}({\rm Adm}^{\{ 1\}}(\mu)) \ .
 $$

 This is shown by Figure 3. Indeed, in the picture, the extreme simplices of ${\rm Adm}^{\{0, 1\}}(\mu)$ 
 (w.r.t. the Bruhat order) are translates of the base simplex  by the translation elements $\{\pm2\}$ which have length $2$, and the other simplices in  ${\rm Adm}^{\{0, 1\}}(\mu)$ correspond to the
 elements smaller than either of these two elements in the affine Weyl group. 
  On the other hand, $\pi_{\{0\}}^{-1}({\rm Adm}^{\{0\}}(\mu))$
 is the set of simplices $\Delta$ such that the distance from the vertex labeled $\{0\}$
 in the picture to the vertex of the same type  in $\Delta$ is $\leq 2$,
 and $\pi_{\{1\}}^{-1}({\rm Adm}^{\{1\}}(\mu))$ is the set of simplices $\Delta$ such that the 
 distance from the vertex labeled $\{1\}$
 in the picture to the vertex of the same type  in $\Delta$ is $\leq 2$. It is visibly true
  that the set of simplices $\Delta$ satisfying both conditions is the set of simplices marked in the picture, i.e. is equal to ${\rm Adm}^{\{0, 1\}}(\mu)$.  
 
\begin{figure}\label{U21}
\includegraphics[width=12cm]{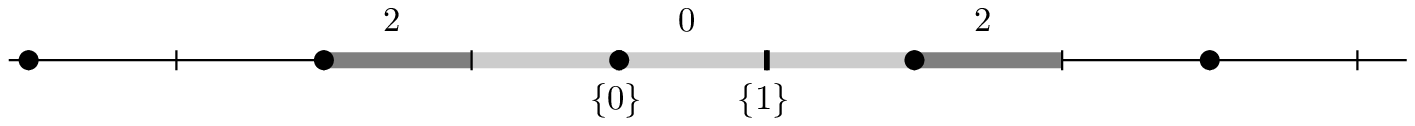}
\caption{The admissible set for $U(2,1)$.}
\end{figure}

 Hence the above inclusions are bijections on the sets of closed points. By Frobenius splitness of Schubert varieties in the flag variety of the unitary group proved in [PR3] ,  the intersection on the RHS of (\ref{inters}) is reduced. Hence all inclusions in (\ref{inters}) are equalities. It follows that   $\bar M^{\rm loc}_{\{0, 1\}}$ is reduced and that all its irreducible components are Schubert varieties,  hence by  [PR3] are normal and with rational singularities.  From Figure 3 we deduce that there are two irreducible components corresponding to the extreme elements of ${\rm Adm}(\mu)$. They intersect in the two Schubert varieties corresponding to the elements of length one in ${\rm Adm}(\mu)$ and these are normal curves which intersect in a single reduced point corresponding to the element of length zero in ${\rm Adm}(\mu)$. We still need to prove that $\bar M^{\rm loc}_{\{0, 1\}}$ is Cohen-Macaulay. However, both irreducible components $X$ and $Y$ are Cohen-Macaulay, and the intersection $X\cap Y$ is Cohen-Macaulay and of codimension $1$ in both $X$ and $Y$. Hence $X\cup Y$ is Cohen-Macaulay, comp. [G\"o1], Lemma 4.22. 
\end{proof}
\medskip

\section{Functor description; the spin condition}
\setcounter{equation}{0}

\subsection{Generalities on symmetric forms of even dimension}\label{symmspin}
 
Let $F$ be a field of characteristic $\neq 2$.  Fix a separable closure $F^{\rm sep}$ of $F$.
Let $V$ be a vector space over $F$ of even dimension $2n$.
We will denote $F$-duals by ${}^*$, so that
$V^*:=\Hom_F(V,F)$ is the $F$-dual of $V$.
There is a perfect, $F$-bilinear pairing
\begin{equation} \label{form1}
\wedge^nV\times\wedge^nV\to \wedge^{2n}V
\end{equation}
defined by
 \ $
(v_1\wedge\cdots\wedge v_n, v'_1\wedge\cdots\wedge v'_n)\mapsto v_1\wedge\cdots\wedge v_n\wedge v'_1\wedge\cdots\wedge v'_n
$.
This gives an isomorphism
\begin{equation}
c:  \wedge^{2n}V\otimes_F(\wedge^nV)^*\to \wedge^nV\ 
\end{equation}
with inverse given by 
$$
v'_1\wedge\cdots\wedge v'_n\mapsto \left(v_1\wedge\cdots\wedge v_n\mapsto 
v_1\wedge\cdots\wedge v_n\wedge v'_1\wedge\cdots\wedge v'_n\right)\ .
$$

\subsubsection{} Now suppose that $V$
supports
a perfect symmetric $F$-bilinear form
\begin{equation*}
h: V\times V\to F\ .
\end{equation*}
The form $h$ induces an $F$-linear isomorphism $b=b_h: V\to V^*$.
For simplicity, let us set $\DD$ for the determinant $F$-line $\wedge^{2n}V$. Then $b$ induces an isomorphism
(``the discriminant") 
\begin{equation*}
d: \DD^{\otimes 2}\xrightarrow{\ \sim\ } F\ .
\end{equation*}
Notice that if $e=e_1\wedge\cdots \wedge e_{2n}\in \Lfr$ with  $\{e_1,\cdots ,e_{2n}\}$ 
an $F$-basis of $V$ then 
\begin{equation}\label{Discr}
d(e\otimes e)=D:=(-1)^n\det((h(e_i,e_j))_{1\leq i,j\leq 2n})
\end{equation}
is the discriminant of the form $h$ in the basis $\{e_i\}_i$.

Define an $F$-algebra structure on $A:=F\oplus {\DD}$ by setting 
$$
(a, e)\cdot (a', e')=(aa'+d(e\otimes e'), ae'+a'e)\ .
$$
 We call $A$ the discriminant 
algebra of $h$. A choice of a non-zero element $e$ in the $F$-line $\DD$
gives an isomorphism 
\begin{equation}
F[x]/(x^2-d(e\otimes e))\xrightarrow{\sim} A\ ,
\end{equation}
and $A$ is either isomorphic to $F\times F$, resp. a quadratic extension of $F$, if 
the discriminant is a square, resp. not a square, in $F$.

\subsubsection{} We will consider the composition
\begin{equation}
a:=c\cdot(id \otimes_F \wedge^n b)\ :\ \DD\otimes_F \wedge^n V \to \wedge^nV .
\end{equation}
where $\wedge^nb: \wedge^nV\to \wedge^n(V^*)=(\wedge^nV)^*$ is given by $b=b_h$.
A choice of a non-zero element $e\in \DD$
provides us with an $F$-linear map:
\begin{equation}
a_e=a(e\otimes -) : \wedge^nV\to \wedge^nV\ .
\end{equation}
For $\lambda\in F^{\times}$, we have $a_{\lambda\cdot e}=\lambda a_e$.
 

\begin{prop} \label{wedge1}
The  map $a: \DD\otimes_F \wedge^n V \to \wedge^nV$ 
provides an $A=F\oplus\DD $-module structure 
which extends the 
$F$-vector space structure on $\wedge^nV$.
\end{prop}

\begin{Proof}
In view of the above, it is enough to show that if 
$e=e_1\wedge\cdots \wedge e_{2n}$, with  $\{e_1,\cdots ,e_{2n}\}$ 
an orthonormal $F$-basis of $V$,  then 
$$
a_e^2=D\cdot {\rm Id}_{\wedge^nV}\ ,
$$
with
$$
D=(-1)^n\prod_{i=1}^{2m}h(e_i,e_i)\ .
$$
For a subset $S=\{i_1,\ldots, i_n\}\subset \{1,\ldots, 2n\}$
with $i_1<\cdots <i_n$ set
$$
e_S=e_{i_1}\wedge\cdots \wedge e_{i_n}\in \wedge^nV\ ,\quad D_S=h(e_{i_1},e_{i_1})\cdots h(e_{i_n},e_{i_n})\in F^{\times} \ .
$$
Denote by  $h^{\wedge^n}$ the symmetric form on $\wedge^nV$ given by the map
$\wedge^n b: \wedge^n V\to (\wedge^n V)^*=\wedge^n V^*$. 
If $S\neq T$, then $h^{\wedge n}(e_S,e_T)=0$.
If $S=T$, then $h^{\wedge n}(e_S,e_S)=D_S$. Therefore, as $S$ runs over all
subsets of $\{1,\ldots, 2n\}$ of order $n$, the elements $e_S$ form an orthogonal
basis of $\wedge^nV$ for the form $h^{\wedge^n}$ and the map
$\wedge^n b $ is given by
$$
e_S=e_{i_1}\wedge\cdots \wedge e_{i_n}\mapsto D_S\cdot e_S^*=e^*_{i_1}\wedge\cdots \wedge e^*_{i_n}\ .
$$
On the other hand, the map $c$ sends $e_S^*$ to $\eta_S\cdot e_{S^c}$
where $S^c=\{1,\ldots, 2n\}\setminus S$ is the complement of $S$ and $\eta_S=\pm 1$ is such that
$$
e_{S}\wedge e_{S^c}=\eta_S\cdot e\ .
$$
Note that $\eta_S\cdot \eta_{S^c}=(-1)^n$ (indeed, the form (\ref{form1}) is
symmetric when $n$ is even, and alternating when $n$ is odd).
We conclude that the map $a_e$ sends $e_S$ to $D_S\eta_S\cdot e_{S^c}$.
Hence, the square $a^2_e$ maps $e_S$ to $D_SD_{S^c}\eta_S\eta_{S^c}\cdot e_S$
and the result follows.\endproof
\end{Proof}
\medskip

Now suppose that $R$ is an $A$-algebra. Then the base change $\wedge^nV\otimes_FR=\wedge^n(V\otimes_FR)$ supports an $A\otimes_FR$--module
structure where the first factor acts via the $A$-module structure of Proposition \ref{wedge1}. Hence, it splits into a direct 
sum of two $R$-modules
\begin{equation}\label{split2}
\wedge^nV\otimes_FR=(\wedge^nV\otimes_FR)_+\oplus (\wedge^nV\otimes_FR)_-
\end{equation}
where the $+$, resp. $-$, part is where $a\otimes  1$ acts as $1\otimes  a$, resp. as $1\otimes  \tau(a)$,
and $\tau:A\to A$ is the non-trivial involution over $F$. 

\subsubsection{} \label{split}

Suppose now that the form $h$ is split; i.e the vector space $V$ has a  basis $\{e_1,\ldots , e_{2n}\}$
with $h(e_i, e_{2n+1-j})=\delta_{ij}$. 
In this case, $D=1$, $A=F\times F$, and
Proposition \ref{wedge1} implies that for $e=e_1\wedge\cdots \wedge e_{2n}$, we have
\begin{equation}
a^2_e={\rm Id}_{\wedge^nV}\ .
\end{equation}
The element $e$ gives the structure of an $A$-algebra to $F$ by sending $e$ to $1$ and we have
\begin{equation*}
\wedge^n V= (\wedge^nV)_+\oplus (\wedge^nV)_-\ 
\end{equation*}
where
$(\wedge^nV)_+$, resp. $(\wedge^nV)_-$ is the eigenspace
of $\wedge^mV$ where $a_e$ has eigenvalue $+1$, resp. $-1$. 
The subspaces
$(\wedge^nV)_+$,  $(\wedge^nV)_-$, depend on the choice of (split) basis 
but the two element set  $\{(\wedge^nV)_+ , (\wedge^nV)_-\}$ does not.
We can see that the action of the special othogonal group $SO(V,h)$ on $\wedge^nV$ preserves the subspaces 
$(\wedge^nV)_+$, $(\wedge^nV)_-$. The element of determinant $-1$ in $O(V,h)$ permutes 
the two subspaces.

In this case, we can easily obtain 
an explicit description of the subspaces $(\wedge^nV)_\pm$,
as follows: For each subset $S$ of $\{1,\ldots, 2n\}$ of order $n$
we set 
$$
2n+1-S:=\{t\ |\ t=2n+1-s, s\in S\}\ .
$$ 
Define a permutation 
$$
\sigma_{S}: \{1,\ldots, 2n\}\to \{1,\ldots, 2n\}
$$
by sending $\{1,\ldots, n\}$ to the
elements of  $ 2n+1-S $ in {\sl decreasing} order and sending the remaining
$\{n+1,\ldots, 2n\}$ to the elements of  the complement
$(2n+1-S)^c$ in {\sl increasing} order.

\begin{lemma}\label{explicit} The subspace $(\wedge^nV)_\pm$ of $\wedge^nV$
is generated by the elements
$$
e_S\pm{\rm sign}(\sigma_S)e_{(2n+1-S)^c}\ .
$$
for $S$ running over all subsets of $\{1,\ldots, 2n\}$ of order $n$.
\end{lemma}

\begin{Proof}
A straightfoward calculation from the definitions gives
that
$$
a_e(e_S)={\rm sign}(\sigma_S)e_{(2n+1-S)^c}\ .
$$
Hence, since $a_e^2={\rm Id}_{\wedge^nV}$, we have
${\rm sign}(\sigma_S)={\rm sign}(\sigma_{(2n+1-S)^c})$.
The result now follows.
\endproof
\end{Proof}

\subsubsection{} We continue to assume  that the form $h$ on $V$ is split with a basis
$\{e_1,\ldots, e_{2n}\}$ such that $h(e_i, e_{2n+1-j})=\delta_{ij}$. 
If $W$ is an totally isotropic subspace of $V$ of dimension $n$, then 
the line $\wedge^nW$ is either contained in $(\wedge^nV)_+$ or in  $(\wedge^nV)_-$.
From the above we can see that   $  e_1\wedge \cdots\wedge e_n   \subset (\wedge^nV)_+$
while   $ e_1\wedge \cdots\wedge e_{n-1}\wedge e_{n+1}  \subset (\wedge^nV)_-$. Now  any basis of an $n$-dimensional isotropic subspace of $V$ can be completed to form 
a split basis of $V$.  Hence, using induction, we can conclude that if
$W$, $W'$ are arbitrary isotropic $n$-dimensional subspaces of $V$ then
$\wedge^nW$ and $\wedge^nW'$ are contained in the same eigenspace $(\wedge^nV)_\pm$ if and only
if 
\begin{equation}\label{interparity}
{\rm dim}(W\cap W')\equiv n\ {\rm mod}\, 2\ .
\end{equation}

\subsection{The spin correction} We now return to the notation of \S \ref{parahoric}. In particular, $F/F_0$ is a 
tamely ramified quadratic extension with automorphism $a\mapsto \bar a$, 
and $\phi: V\times V\to F$ is a  perfect hermitian form on the $F$-vector space
$V$ of dimension $n\geq 3$. 
Suppose $\{e_1,\ldots , e_n\}$ is a basis of $V$ with $\phi(e_i, e_{n+1-j})=\delta_{ij}$.
We will apply the constructions of \S \ref{symmspin} to $V$ considered as a $2n$-dimensional $F_0$-vector space with the symmetric form
$$
h(v,w)=\frac{1}{2}{\rm Tr}_{F/F_0}(\phi(v,w))\ .
$$
Write $m=[n/2]$ and consider the $F_0$-basis
$$
\{-\pi^{-1}e_1, \ldots , -\pi^{-1}e_m, e_{m+1},\ldots, e_{n}, e_1,\ldots , e_m,  \pi e_{m+1},\ldots, \pi e_{n}\}
$$
of $V$. The  $\O_{F_0}$-lattice spanned by this 
basis is invariant for the action of $\O_F=\O_{F_0}[\pi]$ and is equal to 
$\Lambda_{m}$.

There are two cases:

(I) $n=2m+1$ is odd. Then the discriminant $D$ of $h$ in this basis is equal to $\pi_0$
and the discriminant algebra $A$ is isomorphic to $F$. Hence, the form $h$ is not split over $F_0$.
However, after base changing to $F$ we can replace in this list the   vectors $e_{m+1}$, $ \pi e_{m+1}$,
by 
$$
e_{m+1}-\frac{\pi e_{m+1}}{\sqrt{\pi_0}}\ , \ \frac{1}{2}\left(e_{m+1}+\frac{\pi e_{m+1}}{\sqrt{\pi_0}}\right)\ .
$$
The resulting new (ordered) basis splits the form $h\otimes_{F_0}F$. 
We can consider the corresponding eigenspaces
$(\wedge^nV\otimes_{F_0}F)_\pm$ of $\wedge^n V\otimes_{F_0}F$. 

(II) $n=2m$ is even. Then the above basis  splits the form $h$ and we can consider 
the eigenspaces $(\wedge^nV)_\pm$ of $\wedge^n V$. 
\smallskip

 Now consider a choice of $(r,s)$ with $n=r+s$.
 Recall our definition 
of the reflex field $E$. We have $E=F_0$, if $r=s$, and $E=F$ if $r\neq s$. Note that if $r=s$, then $n$ is even
and so we can always make sense of the eigenspaces $(\wedge^n V\otimes_{F_0}E)_\pm$
of $\wedge^nV$. (If $n$ is even and $r\neq s$, we set $(\wedge^n V\otimes_{F_0}E)_\pm=
(\wedge^n V)_{\pm }\otimes_{F_0}F$.)  If $\Lambda$ is one of the $\O_F$-lattices $\Lambda_j$
of $V$ defined in \S \ref{parahoric}, we set
\begin{equation}
(\wedge^n\Lambda\otimes_{\O_{F_0}}\O_E)_\pm :=(\wedge^n\Lambda\otimes_{\O_{F_0}}\O_E)\cap (\wedge^n V\otimes_{F_0}E)_\pm \ .
\end{equation}

\subsubsection{} Now let $I$ be a subset of $\{0,\ldots, m\}$ as in \S \ref{unimoduli} and recall the definition
of the  ``naive" unitary local models $M^{\rm naive}_I$ of \S \ref{naive} corresponding to $I$ and $(r,s)$. 
Recall also the definition of the closed subscheme $M^{\wedge}_I$ of $M^{\rm naive}_I$ defined by requiring the additional 
exterior power condition (e) of 
(\ref{1e4}). 

We define a subfunctor 
$M_I$ of $M^{\rm naive}_I$ by specifying that a point of $M_I$ with values in an $\O_E$-scheme $S$ is given by
an $\O_F\otimes_{\O_{F_0}}\O_S$-submodule $\F_j\subset \Lambda_j\otimes _{\O_{F_0}}\O_S$ for each $j=k\cdot n\pm i$, $i\in I$,
that in addition to the conditions (a)-(e) also satisfies
\smallskip

f) ({\sl Spin condition}) For each $j=k\cdot n\pm i$, $i\in I$, the line  $\wedge^n \F_j\subset   \wedge^n(\Lambda_j\otimes_{\O_{F_0}}\O_S)$ is contained in 
the subspace $(\wedge^n\Lambda_j\otimes_{\O_{F_0}}\O_E)_{\pm}\otimes_{\O_E}\O_S$ with $\pm=(-1)^s$.
\smallskip

It follows from the definition that $M_I$ is represented by a closed subscheme of $M^{\rm naive}_I$
which is contained in $M^{\wedge}_I$; we also denote this subscheme  by $M_I$. 

\subsubsection{} 
We claim that the generic fibers of $M_I$ and $M^{\rm naive}_I$ are equal. Suppose that $S$ is an $E$-scheme.
Then, we can see that condition (f) is open and closed on $S$.
Recall that the generic fiber  $M^{\rm naive}_I\otimes_{\O_E}F$ is isomorphic to the Grassmannian ${\rm Gr}(r, n)_F$. 
Hence, $M^{\rm naive}_I\otimes_{\O_E}E$ is connected and so it is enough to show that (f) is satisfied at some point of the generic fiber. 
It is enough to check this for one choice of signature:
Indeed, we can easily find subspaces $\F_j$, $\F'_j$,
for signatures $(r,s)$ and $(r-1, s+1)$ whose intersection has dimension 
$n-1$; using   (\ref{interparity}) we can inductively reduce to the case
of a single signature choice. Let us take this  to be $(m+1, m)$ if $n=2m+1$ and $(m, m)$ if $n=2m$.
The subspaces given by 
$$
\F^{\rm odd}=\langle\, \pi^{-1}e_1, \ldots , \pi^{-1}e_m, e_1,\ldots , e_m, 
\pi e_{m+1}- {\sqrt{\pi_0}e_{m+1}}  
\, \rangle\ \ ,\hbox{\rm \ and,}
$$
$$
\F^{\rm even}=\langle\, \pi^{-1}e_1, \ldots , \pi^{-1}e_m, e_1,\ldots , e_m 
\, \rangle\ ,
$$
respectively, give points of $M^{\rm naive}_I\otimes_{\O_E}E$ for these 
signatures. We can now see, by calculating the 
dimension of the intersection of $\F^{\rm odd}$, resp. $\F^{\rm even}$, with the standard isotropic subspace
spanned by the first $n$-basis vectors of our chosen split basis of $V\otimes_{F_0}F$, that 
these satisfy condition (f). To recap, we have closed immersions
\begin{equation}
M^{\rm loc}_I\subset M_I\subset M^{\wedge}_I\subset M^{\rm naive}_I
\end{equation}
which are isomorphisms on the generic fiber.

\begin{conjecture}  The scheme $M_I$ is flat over $\O_E$. Equivalently, we have $M^{\rm loc}_I=M_I$.
\end{conjecture}

\begin{Remarks}\label{computer}
{\rm  i) The definition of $M_I$ and the above conjecture are also motivated by a similar 
construction that we have found to be effective in the case of even orthogonal groups (see
\S \ref{ortho}).

ii) This conjecture is supported by some computational evidence that we obtained with the help of Macaulay 2.
In particular, we verified the conjecture for the local models of unitary groups in $3$ and $4$ variables
when $F_0={\bf F}_p((t))$, $F={\bf F}_p((u))$ with $u^2=t$, for various (small) primes $p>2$.

iii) Recall from Remark \ref{nonflat} that the moduli scheme $M^\wedge_I$,  
in which  we omit condition (f),  is not flat   over $\O_E$ in general: 
Indeed, for $n=2m$, $I=\{m\}$ and   signature $(r,s)$
with $s$ odd the subspace $\pi(\Lambda_m\otimes_{\O_{F_0}}k)\subset \Lambda_m\otimes_{\O_{F_0}}k$ gives
a point of the special fiber of   $M^{\wedge}_{\{m\}}$
which does not lift to the generic fiber. We can see that
this point does not satisfy the spin condition (f).

iv) Another interesting example is provided by the case $n=3$, $(r,s)=(2,1)$.
Then we can see that both $M^\wedge_{\{1\}}$ and $M^\wedge_{\{0,1\}}$ are not flat
while our computer calculations 
suggest that $M_{\{1\}}$ and $M_{\{0,1\}}$ are flat.


}
 \end{Remarks}
\medskip

\section{Remarks on the case of (even) orthogonal groups}\label{ortho}
\setcounter{equation}{0}

Here we will discuss a particular type of orthogonal local model
that controls the singularities of certain PEL Shimura varieties;
we will concentrate on the algebraic definition of these models 
and leave the connection to group theory and loop groups for another occasion. Also for the 
connection to the reduction of Shimura varieties we
will refer to [RZ].

We continue with the notations of \S \ref{symmspin}. In addition, assume
that $F$ is local with ring of integers $\O_F$, uniformizer $\pi$
and residue field $k$ of characteristic $\neq 2$.
If $\Lambda$ is an $\O_F$-lattice
in the $F$-vector space $V$, we denote by $\hat\Lambda$
its dual lattice, i.e the image
of $\Lambda^*:={\rm Hom}_{\O_F}(\Lambda, \O_F)$
under the map
$$
\Lambda^*\to \Lambda^*\otimes_{\O_F}F=V^*\buildrel b^{-1}\over \longrightarrow V\ .
$$
The restriction of the form $h$ gives a perfect $\O_F$-bilinear pairing
$$
h_\Lambda: \Lambda\times\hat\Lambda\to \O_F\ .
$$
In what follows we will adhere to the terminology
of [RZ] (note however, that 
in their set-up $F=\Q_p$). Let $\L=\{\Lambda\}$ be a self-dual periodic chain of $\O_F$-lattices
in $V$. Recall that a lattice chain in $V$ is by definition a (non-empty) collection
of lattices such that if $\Lambda$, $\Lambda'$ are in the chain then either
$\Lambda\subset \Lambda'$ or $\Lambda'\subset\Lambda$. A lattice chain
defines a category with objects the lattices and morphisms given by inclusions.
``Periodic" means that $\Lambda\in\L$ implies that $a\Lambda\in\L$ for 
every $a\in F^{\times}$; ``self-dual"
means that $\Lambda\in \L$ implies that $\hat\Lambda\in\L$.

\subsection{The naive local models}\label{orthNaive}

Let us formulate a moduli problem $M^{\rm naive}$ on the category of $\O_F$-schemes. A point of $M^{\rm naive}$ with values 
in an $\O_F$-scheme $S$ is given by functors from the category given by 
the lattice chain $\L$ to the category of $\O_S$-modules 
\begin{equation}
\Lambda\mapsto \F_\Lambda
\end{equation}
together with a morphism of functors 
given by injections
\begin{equation}
j_\Lambda: \F_\Lambda\to \Lambda\otimes_{\O_F}\O_S
\end{equation}
such that:

a) via $j_\Lambda$, $\F_\Lambda$ identifies
with an $\O_S$-locally direct summand of 
rank $n$ of $\Lambda\otimes_{\O_F}\O_S$;

b) if $\Lambda'=\pi\Lambda$, the isomorphism $\Lambda\otimes_{\O_F}\O_S\to
\pi\Lambda\otimes_{\O_F}\O_S=\Lambda'\otimes_{\O_F}\O_S$ 
given by multiplication by the uniformizer $\pi$ induces an isomorphism 
between $\F_{\Lambda}$ and $\F_{\Lambda'}$;
 
c) we have $\F_{\hat\Lambda}=\F_{\Lambda}^{\perp}$, where $\F_{\Lambda}^{\perp}$
is the orthogonal complement of $\F_\Lambda\subset \Lambda\otimes_{\O_F}\O_S$
under the perfect pairing
\begin{equation}
(\hat\Lambda\otimes_{\O_F}\O_S)\,\times\, (\Lambda\otimes_{\O_F}\O_S)\to \O_S
\end{equation}
given by $h_\Lambda\otimes_{\O_F}\O_S$.
\smallskip

We can see that the functor $M^{\rm naive}$ is represented by a projective 
scheme over $\Spec\O_F$ which we will   also denote 
by $M^{\rm naive}$; this is the naive local model defined in [RZ].
If $S$ is an $F$-scheme, then the conditions above
imply that $W:=\F_\Lambda\subset V\otimes_F\O_S$
is independent of the choice of $\Lambda\in\L$.
We can see that generic fiber $M^{\rm naive}_F$ of the naive local
model is the orthogonal Grassmannian of isotropic
subspaces $W$ of $V$ of dimension $n$.

\subsection{The corrected local models}
It was observed by Genestier ([Ge2]),
that the naive local models are not flat
in general. In this paragraph, we suggest a correction in their definition
that should address this problem.
We continue with the notations of the previous paragraph.

Recall $A=F\oplus \DD$ is the discriminant algebra of $h$. By (\ref{split2}) we have
$$
\wedge^nV\otimes_FA=(\wedge^nV\otimes_FA)_+\oplus (\wedge^nV\otimes_FA)_- \ ,
$$
where the $+$, resp. $-$, part is where $a\otimes  1$ acts as $1\otimes  a$, resp. as $1\otimes  \tau(a)$,
and $\tau:A\to A$ is the non-trivial involution over $F$.

Now let $\O_A$ denote the integral closure of $\O_F$ in $A$. For any $\O_F$-lattice $\Lambda$ as above consider the 
$\O_A$-lattice $\Lambda_{\O_A}:=\Lambda\otimes_{\O_F}\O_A$. We have
$
\wedge^n(\Lambda_{\O_A})\subset \wedge^n(V\otimes_FA)
$.
Let us set 
$$
(\wedge^n\Lambda_{\O_A})_\pm =\wedge^n(\Lambda_{\O_A})\cap (\wedge^nV\otimes_FA)_\pm.
$$
This is an $\O_A$-lattice in $(\wedge^nV\otimes_FA)_\pm$. 

Denote by $M$ the closed subscheme of the base change $M^{\rm naive}\otimes_{\O_F}\O_A$
whose $S$-points for an $\O_A$-scheme $S$ 
are given by $S$-points $\F_\Lambda \subset \Lambda\otimes_{\O_F}\O_S$ of $M^{\rm naive}\otimes_{\O_F}\O_A$  
which satisfy
$\wedge^n_{\O_S}(\F_\Lambda)\subset (\wedge^n\Lambda_{\O_A})_+\otimes_{\O_A}\O_S$
(for all $\Lambda$). The corrected local model is by definition $M$ regarded as a scheme over $\Spec(\O_F)$.
By its construction $M$ comes with a morphism
\begin{equation}\label{stein}
M\to \Spec(\O_A)\to \Spec(\O_F).
\end{equation}
There is a natural morphism of $\O_F$-schemes
$q: M \to M^{\rm naive}$.
We will see below that $q$
induces an isomorphism on the generic fibers. However,  $q$ is not  
a closed immersion in general.
\begin{conjecture}
The scheme $M$ is flat over $\Spec\O_F$.
\end{conjecture}

In addition to the examples of the following paragraphs, there is a fair amount
of computational evidence for this (when $n\leq 8$, $F=\fp((t))$,
with $p$ a small prime.) 

\subsubsection{}
To understand this construction a little better,  let us 
consider the situation over $F$, i.e the generic fibers.
Suppose that $W$ is an $S$-valued point 
of $M^{\rm naive}_F$ with $S$ an $F$-scheme (i.e a totally isotropic subbundle of $V\otimes_F\O_S$); the short exact sequence
$$
0\to W\to V\otimes_F\O_S\to (V\otimes_F\O_S)/W\buildrel {b^{-1}}\over {\simeq} W^*\to 0
$$ 
induces an isomorphism $\O_S\simeq  \wedge^nW\otimes\wedge^n(W^*)\simeq \det(V)\otimes_F\O_S\simeq \DD\otimes_F\O_S$. We can see that
this in turns equips $\O_S$ with an $A$-algebra structure. 
We deduce that $M^{\rm naive}_F$ has a natural structure of an $A$-scheme: This structure has the property that if 
$W$ gives an $S$-valued point of $M^{\rm naive}_F$, then the homomorphism $A\otimes_F \wedge^nW\to \wedge^nW$ induced
by $S\to \Spec(A)$ is obtained by restricting the homomorphism 
$A\otimes_F \wedge^nV\to \wedge^nV$
of Proposition \ref{wedge1}. Therefore, we have
$$
M^{\rm naive}_F=M_F\ .
$$
Now denote by $F^{\rm sep}$ a separable closure of $F$; we fix an $F$-homomorphism $A\to F^{\rm sep}$. 
The Pl\"ucker map $W\mapsto \wedge^nW$ gives an embedding
in projective space
\begin{equation*}
M_{F^{\rm sep}}\to {\bf P}(\wedge^nV_{F^{\rm sep}})\ .
\end{equation*}
We can deduce that the image of $M_{F^{\rm sep}}$  under the Pl\"ucker
embedding
\begin{equation*}
M_{F^{\rm sep}}\to {\bf P}(\wedge^nV_{F^{\rm sep}})={\bf P}((\wedge^nV_{F^{\rm sep}})_{ +}\oplus (\wedge^nV_{F^{\rm sep}})_{-}) 
\end{equation*}
lies in the (disjoint) union of the two linear subspaces
${\bf P}((\wedge^nV_{F^{\rm sep}})_{ +})$ and 
${\bf P}((\wedge^nV_{F^{\rm sep}})_{-})$
where the decomposition is obtained as above.
Since the form $h$ splits over $F^{\rm sep}$ we can choose a basis for $V_{F^{\rm sep}}$
as in \ref{split}. 
Since ${\rm SO}(V_{F^{\rm sep}}, h_{F^{\rm sep}})$ is connected,
from \ref{split} and the above, we conclude that $M_{F^{\rm sep}}$ has exactly two
connected components each isomorphic to the special orthogonal Grassmannian ${\rm SO}(V_{F^{\rm sep}}, h_{F^{\rm sep}})/P$;
they are separated by asking that the top exterior power of the
isotropic subspace is contained in $(\wedge^nV_{F^{\rm sep}})_+$,
or $(\wedge^nV_{F^{\rm sep}})_-$ respectively. In fact, it is well-known
that both of these components are of dimension $\binom{n}{2}$. 

\subsection{The split case; examples} Assume that the form $h$ is split
and choose a basis $\{e_1,\ldots, e_{2n}\}$ as in 
\S \ref{split}. Set $\Lambda_0=\la e_1,\ldots,e_{2n}\ra$
and let $\Pi: V\to V$ be the linear map defined by $\Pi(e_1)=\pi\, e_{2n}$,
$\Pi(e_i)=e_{i-1}$, $1\leq i\leq 2n$.
We have $\Pi^{2n}=\pi$. Now let us set $\Lambda_{-i}=\Pi^i\Lambda_0$; 
We can see that $\hat\Lambda_i=\Lambda_{-i}$ so that 
the form induces  perfect $\O_F$-bilinear pairings
\begin{equation}
\Lambda_{-i}\times\Lambda_i\to \O_F\ .
\end{equation}
Now if $I\subset \{0,\ldots, n\}$ is a non-empty subset,
$$
\L_I=\{\pi^k\Lambda_{\pm i}\ |\ k\in\Z\, ,\ i\in I\}.
$$
is a self-dual periodic lattice chain in $V$. We will denote
by $M^{\rm naive}_I$ the ``naive" local model associated to $V$, $h$ 
and the lattice chain $\L_I$. As it was observed by Genestier ([Ge2]),
when $\{0,n\}\subset I$ the scheme $M^{\rm naive}_I\to \Spec\O_F$ is not flat;
he pointed out that the problem is created by the existence 
of two connected components in the (isomorphic) generic fibers   $(M^{\rm naive}_{\{i\}})_F\simeq {\rm OGr}(n, 2n)_F$, $i\in I$. Roughly speaking, the closures of 
several pairs of these components appear in the special fiber 
of the scheme 
$M^{\rm naive}_{I}\subset \prod_{j\in I\cup{-I}}{\rm Gr}(n, 2n)_{\O_F}$. (See  the examples below). His observation 
motivated our work in these sections.
\medskip

\noindent{\bf Example 1:} Let us consider the case of the orthogonal
group ${\rm O}_2$ over $\Q_p$, i.e take $n=1$, $F=\Q_p$, $\pi=p$, and $I=\{0, 1\}$.
In this case we have $\Lambda_0=\la e_1, e_2\ra$, $\Lambda_{-1}=\Pi\Lambda_0=\la e_1, pe_2\ra $. 
For a $\Z_p$-scheme $S$ the $S$-points
of the naive local model $M^{\rm naive}$ are given by diagrams
\begin{equation}
\begin{matrix}
\Lambda_{-1}\otimes_{\Z_p}\O_S&\buildrel \phi\over\to& \Lambda_0\otimes_{\Z_p}\O_S&\buildrel \psi\over\to&
\Lambda_{-1}\otimes_{\Z_p}\O_S&\\
\cup&&\cup&&\cup\\
\F_{-1}&\to&\F_0&\to &\F_{-1}\\ 
\end{matrix}
\end{equation}
such that $\F_0$, $\F_{-1}$ are locally direct summands of rank 
$1$, $\F_0$ is isotropic for $h$, and $\F_{-1}$ is isotropic
for the perfect form on $\Lambda_{-1}\otimes_{\Z_p}\O_S$ given by $p^{-1}h$. 
In these diagrams, the map $\phi$ is given by base changing the inclusion $\Lambda_{-1}\subset \Lambda_0$
while the map $\psi$ is given by base changing the map $\Lambda_0\to \Lambda_{-1}$
given by multiplication by $p$. 

Setting
$$
\F_0= \la ae_1+be_2\ra\ ,\ \ \  \F_{-1}=\la ce_1+dpe_2\ra \ ,
$$
we see that the orthogonality conditions give the (homogeneous) equations:
$$
ab=0\ , \ \ \ cd=0\ .
$$
If $p$ is invertible in $\O_S$, the submodule $\F_{-1}$ is determined
by $\F_0$. We can see that there are two possibilities for $\F_0$:
either $\F_0=\la e_1\ra$, or $\F_0=\la e_2\ra$
and so the generic fiber $M^{\rm naive}_{\Q_p}$ consists of a disjoint union
of two copies of $\Spec\Q_p$.
Let us now consider the special fiber, i.e assume that $p=0$ in $\O_S$.
Taking into account that $\phi(\F_{-1})=\la ce_1\ra  \subset \F_0=\la ae_1+be_2\ra$, $\psi(\F_0)=\la bpe_2\ra \subset \F_{-1}=\la ce_1+dpe_2\ra$,
we see that there are three  possibilities:
$$
\F_{-1}=\la e_1\ra\ , \qquad \F_0=\la e_1\ra\ , \ \  {\rm or,}
$$
$$\F_{-1}=\la pe_2\ra\ , \qquad 
\F_0=\la e_2\ra\ ,  \ \  {\rm or,}
$$
$$
  \F_{-1}=\la pe_2\ra\ , \qquad \F_0=\la e_1\ra\ .\ \ \ \ \ \ \ 
$$
Therefore, the special fiber consists of three copies of $\Spec\fp$
and so $M^{\rm naive}_I$ is not flat over $\Spec\Z_p$. The point
corresponding to the third possibility above does not lift 
to the generic fiber. Notice that in this case
$$
(\wedge^nV)_+=\la e_1\ra\ , \qquad (\wedge^nV)_-=\la e_2\ra\ .
$$
We can now readily see that
$M=\Spec\Z_p\,\sqcup\,\Spec\Z_p$ which is of course flat.
\medskip

\noindent{\bf Example 2:} In this example, we 
 take $n=2$ (which corresponds to $O_4$), $F=\Q_p$, $\pi=p$, and $I=\{1\}$. We have
$$
(\wedge^2V)_+=\la e_1\wedge e_2, e_3\wedge e_4, e_1\wedge e_4+e_2\wedge e_3\ra 
$$
The lattices are $\Lambda_{-1}=\Pi\Lambda_0=\la e_1, e_2, e_3, pe_4\ra$, $\Lambda_1=\Pi^{-1}\Lambda_0=\la e_1/p, e_2, e_3, e_4\ra$,
$p^{-1}\Lambda_{-1}=p^{-1}\Pi\Lambda_0=\la e_1/p, e_2/p, e_3/p, e_4\ra$ and 
the naive local model is given by diagrams
$$
\begin{matrix} \la e_1, e_2, e_3, pe_4\ra_S&\buildrel \phi\over\to& \la e_1/p, e_2, e_3, e_4\ra_S&\buildrel \psi\over\to&
 \la e_1/p, e_2/p, e_3/p, e_4\ra_S\cr
\cup&&\cup&&\cup\cr
\F_{-1}&\to&\F_1&\to &\F_{-1}^{(1/p)}\cr
\end{matrix}
$$
with $\F_i$ of rank $2$ such that $\F_{1}=(\F_{-1})^\perp$. 

In this case
$$
(\wedge^2\Lambda_{-1})_+=\la e_1\wedge e_2,\ e_3\wedge pe_4,\  e_1\wedge pe_4+p\cdot (e_2\wedge e_3)\ra\ ,
$$
$$
\ (\wedge^2\Lambda_{1})_+=\la e_1/p\wedge e_2,\ e_3\wedge pe_4,\ p\cdot (e_1/p\wedge e_4)+e_2\wedge e_3\ra\ .
$$
Our additional spin conditions amount to
\begin{equation}\label{sp1}
\wedge^2\F_{-1}\ \subset\ \la e_1\wedge e_2,\ e_3\wedge pe_4,\ e_1\wedge pe_4+p\cdot (e_2\wedge e_3)\ra_S\ ,
\end{equation}
\begin{equation}\label{sp2}
\wedge^2\F_1\ \subset\ \la e_1/p\wedge e_2,\ e_3\wedge e_4,\ p\cdot (e_1/p\wedge e_4)+e_2\wedge e_3\ra _S\ .
\end{equation}

We will  examine the special fiber; thus we assume $p=0$ in $\O_S$.
If $(\F_{-1}, \F_1)$ corresponds to a $\bar{\bf F}_p$-valued point
of $M$, 
we can see that either $e_2$ or $e_3$ is in $\F_1$.
When $e_2\in \F_1$, then
$\F_1=\la e_2, ae_1/p+be_3\ra $, $\F_{-1}=\la e_1, -ae_2+bpe_4\ra$, where $(a;b)\in {\bf P}^1$. If 
$e_3\in \F_1$, then $\F_1=\la e_3, ce_2+de_4\ra$, $\F_{-1}=\la p\,e_4, -ce_1+de_3\ra$, with $(c;d)\in {\bf P}^1$.
Those two projective lines intersect at the point $\F_1=\la e_2, e_3\ra$
(at which $\F_{-1}=\la e_1, p\,e_4\ra$). This describes the (reduced) special fiber of the corrected local model. 
It is not hard to see that both irreducible components
lift to characteristic zero. For example, the point
$\F_1=\la e_3, e_4\ra$, $\F_{-1}=\la e_3, p\,e_4\ra$ is a smooth point
on the second component  that lifts; similarly the point
$\F_1=\la e_1/p, e_2\ra$, $\F_{-1}=\la e_1, e_2\ra$ is a smooth 
point of the first component that lifts.

Now on to a ``scheme theoretic" calculation. We will calculate an affine chart
of the special fiber at the singular point $\F_{-1}=\la e_1, p\,e_4\ra$, 
$\F_1=\la e_2,e_3\ra$. Set
$$
\F_{-1}=\la e_1+x_2e_2+x_3e_3,\ p\,e_4+y_2e_2+y_3e_3\ra \ ,
$$
Then we can see that
$$
\F_1=(\F_{-1})^\perp=\la e_2-y_3\frac{e_1}{p}-x_3e_4,\  e_3-y_2\frac{e_1}{p}-x_2e_4\ra\ .
$$
The condition $\phi(\F_{-1})\subset \F_1$ translates to
\begin{equation}\label{incl1}
x_2x_3=y_2y_3=x_2y_3+x_3y_2=0\ ,
\end{equation}
while $\psi(\F_{1})\subset \F_{-1}^{(1/p)}$ translates  to
\begin{equation}\label{incl2}
x_2y_2=x_3y_3=x_2y_3+x_3y_2=0\ .
\end{equation}
The spin conditions (\ref{sp1}), (\ref{sp2}) give 
$$
y_3=x_2=0\ ,\ \ x_2y_3-y_2x_3=0\ .
$$
These equations together with (\ref{incl1}), (\ref{incl2}) amount to  
$$
x_2=y_3=0, \ \ y_2x_3=0\ .
$$
 This together
with the ``set-theoretic" arguments above imply that 
the special fiber is reduced and 
is the union of two ${\bf P}^1$'s intersecting transversely
at a point. Since, both of these ${\bf P}^1$'s lift, we also conclude that $M_{I}$ is flat 
over $\Z_p$. 

Notice that if we omit the spin conditions,
we obtain
$$
x_2x_3=y_2y_3=x_2y_2=x_3y_3=x_2y_3+x_3y_2=0\ .
$$
These are then the equations for an affine chart of the
special fiber of the naive local model; this scheme 
has four irreducible components and is non-reduced
at the origin. We can see that $M_I^{\rm naive}$ 
cannot be flat over $\Z_p$.

\end{document}